\documentclass[12pt]{amsart}
\usepackage[utf8]{inputenc}
\usepackage[margin=1in]{geometry}
\usepackage[dvipsnames]{xcolor}
\usepackage{bbm}
\usepackage{tikz}
\usetikzlibrary{knots}
\usetikzlibrary{arrows}
\usepackage{tikz-cd}
\usepackage{pgffor}
\usepackage{amssymb,amsmath,amsfonts,mathrsfs,amsthm}
\usepackage{enumitem}
\usepackage{adjustbox}
\usepackage{hyperref}
\usepackage{relsize}
\usepackage{graphicx}
\usepackage[font={small}]{caption}
\usepackage{float}
\usepackage{subcaption}
\usepackage{changepage}
\usepackage{changes}
\usepackage[toc]{appendix}
\usepackage{pgfplots}
\pgfplotsset{compat=1.7}
\usepackage{adjustbox}
\usepackage{varwidth}
\usepackage{chngcntr}
\usepackage{tcolorbox}
\usepackage{mathtools}
\usepackage{changes}
\usepackage{chngcntr}
\usepackage{standalone}
\usetikzlibrary{calc}

\newcommand{\R}{\mathbb{R}}

\newcommand{\Z}{\mathbb{Z}}
\newcommand{\C}{\mathbb{C}}

\newcommand{\Hom}{\operatorname{Hom}}

\newcommand{\id}{\operatorname{id}}

\renewcommand{\l}{\ell}

\renewcommand{\sl}{\mathfrak{sl}}

\newcommand{\Cs}{\mathscr{C}}

\newcommand{\q}{\mathfrak{q}}

\newcommand{\Tot}{\operatorname{Tot}}
\newcommand{\BN}{\mathcal{BN}}

\newcommand{\A}{\mathbb{A}}

\newcommand{\F}{\mathcal{F}}

\newcommand{\gmod}{{\rm -}\operatorname{gmod}}
\newcommand{\ggmod}{{\rm -}\operatorname{ggmod}}

\newcommand{\qdeg}{\operatorname{qdeg}}
\newcommand{\adeg}{\operatorname{adeg}}
\newcommand{\lr}[1]{\vert {#1} \vert}

\newcommand{\gRep}{\operatorname{gRep}}

\newcommand{\subcat}{\mathcal{Z}}
\newcommand{\J}{\mathcal{J}}

\newcommand{\size}[1]{\left| #1 \right|}
\newcommand{\cube}{\mathcal{Q}}
\newcommand{\subcube}[3] { #3 \leq_{#1} #2 }
\newcommand{\scube}{ \subcube{m}{v}{u} } 
\newcommand{\ModSp} {\mathcal{M}}
\newcommand{\cubeMod}[2] { \ModSp_{\cube} ( #1 , #2 ) }
\newcommand{\cMod} { \cubeMod{v}{u} } 
\newcommand{\AmbSp}[2] { \mathcal{E} ( #1 , #2 ) }
\newcommand{\AmbS} { \AmbSp{v}{u} } 
\newcommand{\EucSp}[2] {\mathbb{E} ( #1 , #2 )} 
\newcommand{\EucS}{ \EucSp{v}{u} } 
\newcommand{\cubeEuc}[2]{ \mathbb{E}_\cube ( #1 , #2 ) }
\newcommand{\cEuc}{\cubeEuc{v}{u}} 

\newcommand{\KG}[1] { \Gamma_{Kh}(#1) }
\newcommand{\de} {\partial}
\newcommand{\pathPts}[2] { \de_{#1} #2 }
\newcommand{\pathP} {\pathPts{P}{\MspJxy}} 
\newcommand{\pathAmb}[1] {\mathcal{E}(#1)}
\newcommand{\pAmb}{\pathAmb{P}} 

\newcommand{\sqMspE}[4] {\de_{P^{#1} P^{#2}} \MspE{#3}{#4} }
\newcommand{\sqMspExy} {\sqMspE{b}{a}{y}{x}}
\newcommand{\MspCC}[4][] { \ModSp_{#1}^{#2} (#3|_{#2},#4|_{#2}) }
\newcommand{\MspECCxy} {\MspCC[E]{\Cs}{y}{x}} 
\newcommand{\MspCCpxy} {\MspCC{\Cs'}{y}{x}} 

\newcommand{\ConnComp}{\Cs}
\newcommand{\CC}{\ConnComp} 

\newcommand{\X} {\mathcal{X}} 
\newcommand{\ffc} {\mathfrak{F}} 
\newcommand{\ffcJ} {\ffc_J} 
\newcommand{\MspJ}[2] { \ModSp_J ( #1 , #2 ) }
\newcommand{\MspJxy} {\MspJ{y}{x}} 
\newcommand{\MspE}[2] { \ModSp_E ( #1 , #2 ) }
\newcommand{\MspExy} {\MspE{y}{x}} 
\newcommand{\MspH}[2] { \ModSp_H ( #1 , #2 ) }
\newcommand{\MspHxy} {\MspH{y}{x}} 
\newcommand{\Cone}{\mathrm{Cone}}
\newcommand{\Obj}{\mathrm{Ob}}
\newcommand{\td}[1]{\widetilde{#1}}

\newcommand{\rar}[1]{\xrightarrow{#1}}
\newcommand{\bv}{\boldsymbol{v}}

\renewcommand{\o}{\otimes}

\renewcommand{\:}{\colon}

\newcommand{\FKh}{\mathcal{F}_{Kh}}
\newcommand{\Fa}{\mathcal{F}_{\mathbb{A}}} 
\newcommand{\vm}{v_-}
\newcommand{\vp}{v_+}

\newcommand{\thickemb}{
    \lhook\joinrel\relbar\mspace{-10mu}\hookrightarrow
} 

\newcommand{\xthickemb}[1]{
    \lhook\joinrel\relbar\mspace{-10mu}\xhookrightarrow{#1}
} 

\newcommand{\thickembtikz}[2]{
\draw[right hook->] (#1) -- (#2);
\draw[right hook->] (#1.east)+(4pt,0) -- (#2);
}

\newcommand{\xthickembtikz}[3]{
\draw[right hook->] (#1) -- (#2);
\draw[right hook->] (#1.east)+(4pt,0) -- node[above]{$#3$} (#2);
}

\newcommand\Label[1]{&\refstepcounter{equation}(\theequation)\ltx@label{#1}&}

\newtheorem{theorem}{Theorem}[section]
\newtheorem{lemma}[theorem]{Lemma}
\newtheorem{proposition}[theorem]{Proposition}
\newtheorem{corollary}[theorem]{Corollary}

\theoremstyle{definition}

\theoremstyle{definition}

\theoremstyle{remark}
\newtheorem{example}[theorem]{Example}

\theoremstyle{remark}
\newtheorem{remark}[theorem]{Remark}

\theoremstyle{remark}

\theoremstyle{definition}
\newtheorem{definition}[theorem]{Definition}

\theoremstyle{definition}

\begin{document}

\title{Towards an $\sl_2$ action on the annular Khovanov spectrum}

\author[R. Akhmechet]{Rostislav Akhmechet}
\address{Department of Mathematics, University of Virginia, Charlottesville VA 22904-4137}
\email{\href{mailto:ra5aq@virginia.edu}{ra5aq@virginia.edu}}

\author[V. Krushkal]{Vyacheslav Krushkal}
\address{Department of Mathematics, University of Virginia, Charlottesville VA 22904-4137}
\email{\href{mailto:krushkal@virginia.edu}{krushkal@virginia.edu}}

\author[M. Willis]{Michael Willis}
\address{Department of Mathematics, University of California, Los Angeles, CA 90095}
\email{\href{mailto:msw188@ucla.edu}{msw188@ucla.edu}}

\begin{abstract}
Given a link in the thickened annulus, its annular Khovanov homology carries an action of the Lie algebra $\sl_2$, which is natural with respect to annular link cobordisms. We consider the problem of lifting this action to the stable homotopy refinement of the annular homology. As part of this program, the actions of the standard generators of $\sl_2$ are lifted to maps of spectra.  In particular, it follows that the $\sl_2$ action on homology commutes with the action of the Steenrod algebra.  The main new technical ingredients developed in this paper, which may be of independent interest, concern certain types of cancellations in the cube of resolutions and the resulting more intricate structure of the moduli spaces in the framed flow category. 
\end{abstract}

\maketitle
\tableofcontents

\section{Introduction}
The structure of link homology theories is closely related to the representation theory of Lie algebras, and to categorification of quantum groups. This paper initiates the study of lifting these relations to the level of stable homotopy types.  In the process of doing so, we provide the first example (to our knowledge) of constructing a stable homotopy refinement via framed flow categories of a homological invariant involving signs and cancellations. 

We work in the setting of the $\sl_2$ homology theory $Kh_{\A}(L)$ for links $L$ in the thickened annulus, sometimes referred to as sutured annular Khovanov homology. Following constructions in \cite{APS, BN2, Ro}, this triply graded theory may be obtained from the usual Khovanov chain complex \cite{Kh1} of a link diagram $D$ by taking the annular degree zero part of the differential. It was shown in \cite{BPW} that given a tangle $T$, the homology $Kh_{\A}(\widehat T)$ of its annular closure may also be obtained as the Hochschild homology of the complex of bimodules over the Chen-Khovanov algebra \cite{CK} associated to $T$.

Motivated in part by results of Auroux, Grigsby and Wehrli \cite{AGW} and by the work of Lauda \cite{Lauda} on categorified $\sl_2$, Grigsby, Licata and Wehrli showed \cite{GLW} that the annular Khovanov homology $Kh_{\A}(L)$ carries an action of $\sl_2$. Analogous actions of $\sl_n$ on annular Khovanov-Rozansky homology were defined in \cite{QR}.  In fact, the $\sl_2$ action is constructed on the chain complex level, and it is natural with respect to annular link cobordisms. Here a trivial simple closed curve in the annulus is assigned the trivial representation of $\sl_2$, and essential curves are assigned the fundamental representation or its dual; \cite{GLW} extended the action to $Kh_{\A}(L)$ for any annular link $L$ and showed that the action is a link invariant; see Section \ref{sec:sl2 action} for a detailed discussion.

Our results concern the stable homotopy refinement of Khovanov homology, introduced by Lipshitz and Sarkar \cite{LS}. This construction, building on the work of Cohen, Jones and Segal \cite{CJS}, associates to a diagram of an oriented link a framed flow category. The resulting suspension spectrum is a stable homotopy invariant of the link $L$; its cohomology is isomorphic to the Khovanov homology of $L$. An alternative, more combinatorial construction using the Burnside category was given by Lawson, Lipshitz and Sarkar in \cite{LLS}.
The stable homotopy refinement $\X_\A(L)$ of an annular link $L$ may be constructed using an analogous method, cf. \cite[Section 4.3]{AKW}. Alternatively, lifting the homology isomorphism with the Hochschild homology of Chen-Khovanov bimodules,  $\X_\A(L)$ may be defined as the topological Hochschild homology of the tangle invariants introduced by Lawson-Lipshitz-Sarkar in \cite{LLSCK}. 

These new links invariants valued in spectra can be studied using tools of stable homotopy theory. For example, they admit an action of the Steenrod algebra, and they are known to be stronger invariants than the underlying link homology, cf. \cite{Seed}. It is an interesting question what features of link homology admit a lift to the level of spectra. We address this question, giving a new structure on the annular Khovanov spectrum and thus providing a new invariant of annular links and cobordisms.  To this end, our main result is the following.  Let $J$ denote one of the standard generators $E, F, H$ of $\sl_2$ (see Section \ref{sec:sl2 action}). 

\begin{theorem}\label{main thm} Given an annular diagram $D$  of a link $L$ in the thickened annulus, there exists a map of spectra $\J:\X_\A(D) \rightarrow \X_\A(D)$, whose induced map $\J^*$ on cohomology is equal to the action of $J$ on $Kh_{\A}(L)$. The homotopy class of $\J$ is an invariant of the link $L$.
\end{theorem}

As a consequence, we have the following statement about the $\sl_2$ action on annular Khovanov homology.

\begin{corollary}
The action of $\sl_2$ on the annular homology $Kh_{\A}(L)$ of a link $L$ commutes with the Steenrod cohomology operations.
\end{corollary}

We carry out the construction in the setting of framed flow categories, recalled in Section \ref{sec:ffc review}.  (An analogous construction likely can be done using a suitably defined ``framed Burnside category''; compare the discussion in \cite[Section 2.5]{LSsurvey}.) The overall strategy is to define a stable homotopy refinement of the cone of the chain map $J$, see Section \ref{sec:constructing a map in general}. 

An important problem that immediately presents itself is the existence of {\em cancellations}. That is, within the cube-shaped chain complex ${\rm Cone}(J)$, there exist pairs of consecutive edges so that a single generator is sent to a pair of generators, which then cancel by the second edge differential. This does not happen in the usual Khovanov chain complex, since the coefficients of the differential with respect to the canonical generators are either $0$ or $1$. There are signs but still no cancellations in the setting of odd Khovanov homology, whose stable homotopy refinement was given in \cite{SSS}. 

A specific example of cancellations in ${\rm Cone}(J)$ is discussed in Section \ref{sec:cubic illustration}. Its manifestation in the context of framed flow categories is the presence of a pair of oppositely framed points, which cobound a framed interval in the $1$-dimensional moduli space. In general this leads to moduli spaces of dimensions $2$ and higher with non-trivial topology, see Sections \ref{sec:cubic illustration}, \ref{sec:m-1 to m} and  Figure \ref{fig:Hexagon times I}.

The construction of the Khovanov homotopy type \cite{LS} and other recent results in this setting involved framed flow categories whose $n$-dimensional moduli spaces are trivial covers of permutohedra, so topologically they are (disjoint unions of) $n$-dimensional balls. In general, if one attempts to build a framed flow category and the moduli spaces have non-trivial topology in some dimension, the existence of higher-dimensional moduli spaces is not assured. Indeed, there are highly restrictive compatibility conditions on how the moduli spaces of various dimensions fit together in a framed flow category (see Section \ref{sec:ffc review}), and there are choices involved in their construction. The union of moduli spaces up to a dimension $n-1$ are required to form the boundary of the moduli spaces of dimension $n$, forming manifolds with corners. Thus extending the construction to dimension $n$ amounts to finding a collection of null-cobordisms, a problem which a priori might not have a solution.

A key point in our definition of a framed flow category for ${\rm Cone}(J)$ in the proof of Theorem \ref{main thm} is that the moduli spaces may be constructed as (trivial covers of) codimension-1 submanifolds of some canonical spaces. The bulk of the paper, Sections \ref{0 dim sec} - \ref{sec:Path moduli spaces}, deals with setting up the base cases: $0$- and $1$-dimensional moduli spaces. The inductive step,  using the Pontryagin-Thom construction and propagating the codimension-1 condition, is given in Section \ref{sec:m-1 to m}.   

Another important feature in our work, as compared to prior constructions, is the necessity of performing a global combinatorial analysis of configurations of curves and surgery arcs in a link diagram.  Recall that the structure of the 1-dimensional moduli spaces used in the construction of the Khovanov homotopy type is analyzed by hand only for 2- and 3-dimensional faces of the Khovanov cube \cite[Sections 5.4,5.5]{LS}.  These analyses are local in the sense that no more than three surgery arcs need to be considered.  The triviality of the permutohedral covers then ensures that all higher dimensional moduli spaces can be built.  In contrast, the more intricate structure of moduli spaces in our setting requires a detailed analysis of configurations of an arbitrary number of surgery arcs; see Section \ref{sec:Path moduli spaces}.

The construction of the map $\J$ in Theorem \ref{main thm} is stated as Corollary \ref{cor:ffcj refines Cone J}; the main technical result underlying its definition is Theorem \ref{thm:ffcJ}.  The invariance of the homotopy class of $\J$ with respect to various choices made in the construction, as well as invariance under Reidemeister moves, is established in Section \ref{sec:invariance}.

Recall that the $\sl_2$ action on the annular complex is natural with respect to link cobordisms.  However, link cobordism maps on spectra are only defined \cite{LS2} for link cobordisms which are presented as a sequence of Morse moves and Reidemeister moves.  These maps are not known to be well-defined in general, although they are conjectured to be so.

\begin{theorem}
\label{thm:J commutes with cobs}
The map $\J$ constructed in Theorem \ref{main thm} is natural with respect to link cobordisms presented as sequences of Morse moves and Reidemeister moves.
\end{theorem}

The results of this paper concern the action of generators of $\sl_2$. The problem of lifting the relations of $\sl_2$ is outside the scope of our present work; we hope to address it in a future paper. Another interesting problem concerns quantum annular homology $Kh_{\A_\q}(L)$. This theory was defined by Beliakova, Putyra and Wehrli in \cite{BPW}. Extending the work of \cite{GLW}, they showed that $Kh_{\A_\q}(L)$ admits an action of $U_{\q}(\sl_2)$. Our earlier work \cite{AKW} defined a stable homotopy refinement of quantum annular homology, which takes the form of a $\Z/r\Z$-equivariant spectrum  ${\mathcal X}^r_{\A_\q}(L)$, where $r\geq 2$.  Conjecture 1.4 in \cite{AKW} states that the action of $U_{\q}(\sl_2)$ on $Kh_{\A_\q}(L)$  can be lifted to an action on ${\mathcal X}^r_{\A_\q}(L)$. The methods of the present paper do not immediately extend to the setting of ${\mathcal X}^r_{\A_\q}(L)$ which was constructed using an equivariant version of the Burnside category. Nevertheless, our results are consistent with that conjecture, and a reformulation in terms of a ``framed Burnside category", briefly discussed above, is likely to lead to a lift of the action of the generators of $U_{\q}(\sl_2)$.

{\bf  Acknowledgements.} We are grateful to Sucharit Sarkar for many very helpful discussions on framed flow categories and the Burnside category.  MW is grateful to Matt Stoffregen for several helpful conversations.

RA was supported by NSF RTG Grant DMS-1839968 and the Jefferson Scholars Foundation. 
VK was supported in part by the Miller Institute for Basic Research in Science at UC Berkeley, Simons Foundation fellowship 608604, and NSF Grant DMS-1612159.  MW was supported in part by NSF FRG Grant DMS-1563615.  

\section{Background and conventions}

\subsection{Khovanov homology}\label{sec:Khovanov homology}

We begin with a brief overview of the Bar-Natan category $\BN$ and the construction of the formal chain complex $[[D]]$ assigned to a link diagram $D$; for a complete treatment see \cite{BN1}. 

\begin{figure}
\centering
\begin{subfigure}[b]{.3\textwidth}
\begin{center}
\includegraphics{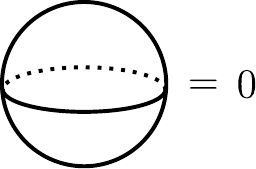}
\end{center}
\caption{Sphere}\label{fig:sphere}
\end{subfigure}%
\begin{subfigure}[b]{.4\textwidth}
\begin{center}
\includegraphics{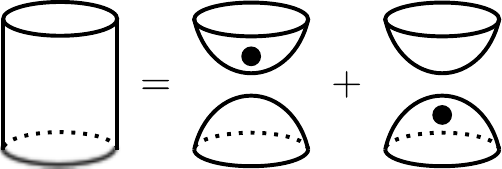}
\end{center}
\caption{Neck-cutting}\label{fig:neck-cutting}
\end{subfigure}\\
\vskip1em
\begin{subfigure}[b]{.3\textwidth}
\begin{center}
\includegraphics{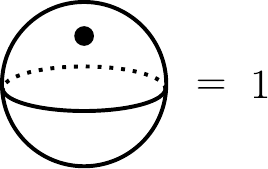}
\end{center}
\caption{Dotted sphere}\label{fig:dotted sphere}
\end{subfigure}%
\begin{subfigure}[b]{.35\textwidth}
\begin{center}
\includegraphics{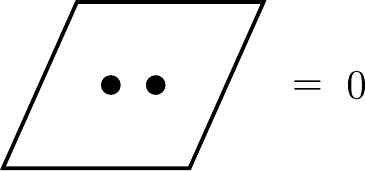}
\end{center}
\caption{Two dots}\label{fig:two dot}
\end{subfigure}
\caption{Bar-Natan relations}\label{fig:BN relations}
\end{figure}

Recall the (dotted) Bar-Natan category $\BN=\BN(\R^2)$. Let $I :=[0,1]$ denote the unit interval. Objects of $\BN$ are formal direct sums of formally graded collections of simple closed curves in the plane $\R^2$. Morphisms are matrices whose entries are formal $\Z$-linear combinations of dotted cobordisms properly embedded in $\R^2 \times I$. Cobordisms are considered up to isotopy relative to the boundary, and they are subject to the local relations shown in Figure \ref{fig:BN relations}. 

Let $A= \Z[X]/(X^2)$ denote the Frobenius algebra with comultiplication, unit, and counit described below. 
\begin{align*}
\Delta \: A & \to A \o A & \eta \: \Z & \to A & \varepsilon \: A & \to \Z \\
1 & \mapsto X \o 1 + 1\o X & 1 &\mapsto 1 & 1 & \mapsto 0\\
X & \mapsto X\o X & & & X & \mapsto 1
\end{align*}
This is the Frobenius algebra underlying $\sl_2$ link homology \cite{Kh1}. Define the quantum grading $\qdeg$ on $A$ by setting 
\begin{equation}\label{eq:1 and X gradings}
\qdeg(X) = -1 \hskip2em \qdeg(1) = 1.
\end{equation}
Let $\Z \gmod$ denote the category of $\Z$-graded abelian groups and graded maps (of any degree) between them. 
The Frobenius algebra $A$ defines a $(1+1)$-dimensional TQFT, which descends to a functor 
\begin{equation}\label{eq:tqft F}
\FKh \colon \BN \to \Z\gmod.
\end{equation}
We will refer to $\FKh$ as the \emph{Khovanov TQFT}.  The construction is reviewed below. 

On an object $\Cs$ which is a collection of $n$ disjoint circles in the plane, $\FKh(\Cs) = A^{\o n}$, generated by elements of the form 
\begin{equation}
\label{eq:generators expanded as labels on each circle}
y=y_1\o \cdots \o y_n\in A^{\o n}
\end{equation}
where each $y_i \in \{1, X\}$; this corresponds to a label of each circle in $\Cs$ by $1$ or $X$.  We refer to such generators as \emph{Khovanov generators}.

Morphisms in $\BN$ are generated by so-called \emph{elementary cobordisms} which contain only a single non-degenerate critical point with respect to the height function $\R^2\times I \rightarrow I$.  Index 0 (cup) and 2 (cap) elementary cobordisms  are assigned the unit and co-unit maps of $A$ respectively; index 1 (saddle) elementary cobordisms are assigned multiplication or comultiplication maps according to whether the cobordism merged two circles into one, or split one circle into two.  The presence of a dot on a cobordism corresponds to multiplication by $X\in A$.

Now let $D$ be a planar diagram for an oriented link $L\subset S^3$.  We recall the construction of the formal Bar-Natan chain complex $[[D]]$ from \cite{BN2}.  First form the \textit{cube of resolutions}.  Label the crossings of the diagram by $1,\ldots, n$. Each crossing may be resolved in one of two ways, called the \textit{0-smoothing} and \textit{1-smoothing}, as in \eqref{fig:0 and 1 smoothings}. 
\begin{equation}\label{fig:0 and 1 smoothings}
\begin{aligned}
\includestandalone{two_smoothings}
\end{aligned}
\end{equation}
For each $u = (u_1,\ldots, u_n)\in \{0,1\}^n$, perform the $u_i$-smoothing at the $i$-th crossing. This results in a collection of disjoint simple closed curves in the plane, which we denote $D_u$. Identifying $\{0,1\}^n$ with the vertices of an $n$-dimensional cube, we decorate each vertex $u$ by its corresponding smoothing $D_u$. 

Let $u=(u_1,\ldots, u_n)$ and $v=(v_1,\ldots, v_n)$ be vertices which differ only in the $i$-th entry, where $u_i = 0$ and $v_i=1$. The diagrams $D_u$ and $D_v$ are identical outside of a small disk around the $i$-th crossing. There is a cobordism from $D_u$ to $D_v$, which is the obvious saddle near the $i$-th crossing and the identity (product cobordism) elsewhere. We will call this the \emph{saddle cobordism} from $D_u$ to $D_v$, and denote it by $d_{u,v}$. Decorate each edge of the $n$-dimensional cube by the corresponding saddle cobordism. The result is a commutative cube in the category $\BN$. There is an assignment $s_{u,v} \in \{0,1\}$ for each edge so that multiplying the edge map $d_{u,v}$ by $(-1)^{s_{u,v}}$ results in an anti-commutative cube (see \cite[Section 2.7]{BN2}, also \cite[Definition 4.5]{LS}).

For $u=(u_1,\ldots, u_n) \in \{0,1\}^n$, set $\vert u \vert = \sum_i u_i$. Now, form the chain complex $[[D]]$ by setting 
\[
[[D]]^i = \bigoplus_{\vert u \vert = i+ n_-} D_{u}\{i + n_+ -n_-\} 
\] 
where $n_-$, $n_+$ are the number of negative and positive crossings in $D$, and $\{-\}$ denotes the formal grading shift in $\BN$. The differential is defined on each summand by the edge map $(-1)^{s_{u,v}}d_{u,v}$. That $[[D]]$ is a chain complex follows from anti-commutativity of the cube.

\begin{theorem}\emph{(\cite[Theorem 1]{BN2})}
If diagrams $D$ and $D'$ are related by a Reidemeister move, then $[[D]]$ and $[[D']]$ are chain homotopy equivalent.
\end{theorem}

Let
\[
CKh(D) := \FKh([[D]])
\]
denote the chain complex obtained by applying the Khovanov TQFT to $[[D]]$; it is precisely the invariant defined in \cite[Section 7]{Kh1}.

\subsection{Annular Khovanov homology}\label{sec:annular Khovanov homology} 

Asaeda-Przytycki-Sikora \cite{APS} introduced a homology theory for links in $I$-bundles over surfaces. The case of the thickened annulus is known as annular Khovanov homology, sometimes called \emph{sutured} annular Khovanov homology. This section reviews its construction, following \cite{Ro, GLW}. 

\begin{figure}
\centering
\includegraphics{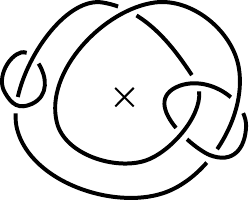}
\caption{An annular link diagram}\label{fig:annular link}
\end{figure}

Let $\A = S^1\times I$ denote the annulus. An \textit{annular link} is a link in the thickened annulus $\A\times I$, and its diagram is obtained by projecting onto the first factor $\A$. Embed $\A$ standardly in $\R^2$ via 
\[
\A = \{ x\in \R^2 \mid 1\leq \lr{x} \leq 2\},
\]
so that an annular link diagram and all of its smoothings are drawn in the punctured plane $\R^2\setminus (0,0)$. We represent the annulus by simply indicating the deleted point using the symbol $\times$. See Figure \ref{fig:annular link} for an example of an annular link diagram drawn using this convention. We distinguish two types of circles in $\A$: \emph{trivial} circles, which are contractible in $\A$, and \emph{essential} ones, which are not contractible. 

Let $\BN(\A)$ denote the Bar-Natan category of the annulus. Its objects are formal direct sums of formally $\Z^2$-graded collections of simple closed curves in $\A$. Morphisms are matrices with entries formal $\Z$-linear combinations of dotted cobordisms properly embedded in $\A \times I$, modulo isotopy relative to the boundary, and subject to the local Bar-Natan relations shown in Figure \ref{fig:BN relations}.

Let $\Z \ggmod$ denote the category of $\Z^2$-graded abelian groups and graded maps (of any bidegree) between them. We now describe the annular TQFT
\[
\Fa \colon \BN(\A) \to \Z\ggmod,
\]
following the construction in \cite{Ro}.

Let $\Cs \subset \A$ be a collection of $n$ trivial and $m$ essential circles. Viewing $\A$ as a subspace of $\R^2$, apply the Khovanov TQFT $\FKh$ from Section \ref{sec:Khovanov homology},
\[
\FKh(\Cs) = A^{\o n} \o A^{\o m}. 
\]
Introduce a second grading, called the \emph{annular} grading and denoted $\adeg$, on $\FKh(\Cs)$ as follows. Every tensor factor corresponding to a trivial circle is concentrated in annular degree zero. For a factor $A$ associated to an essential circle, introduce the notation
\[
\vm = X,\hskip2em \vp = 1
\]
to denote a basis for this copy of $A$, and define $\adeg$ by setting 
\[
\adeg(\vm) = -1 \hskip2em \adeg(\vp) = 1.
\]

The underlying abelian group of $\Fa(\Cs)$ is defined to be $\FKh(\Cs)$, with the bigrading given by $(\qdeg, \adeg)$. To define $\Fa$ on a cobordism $S\subset \A\times I$, first view $S$ as a surface in $\R^2 \times I$ and consider the map $\FKh(S)$ assigned to $S$ by the Khovanov TQFT. It was observed in \cite[Section 2]{Ro} that $\FKh(S)$ is of the form
\begin{equation}\label{eq:non alpha cob split}
\FKh(S) = \FKh(S)_0 + \FKh(S)_-
\end{equation}
where $\FKh(S)_0$ preserves $\adeg$ and $\FKh(S)_-$ lowers $\adeg$. Set 
\[
\Fa(S) := \FKh(S)_0
\]
to be the part of $\FKh(S)$ that preserves $\adeg$. It follows from \eqref{eq:non alpha cob split} that $\Fa$ is functorial with respect to composition of cobordisms, and it clearly factors through the relations in Figure \ref{fig:BN relations}. We will refer to $\Fa$ as the \emph{annular TQFT}. Note that if $S$ carries $d$ dots, then $\Fa(S)$ is a map of $(\qdeg, \adeg)$ bidegree
\[
(\chi(S) - 2d, 0).
\]

We distinguish the bigraded modules assigned to trivial and essential circles by writing 
\[
V = \Fa(C)
\]
if $C$ is an essential circle, with basis written as $\{\vm, \vp\}$. The notation $A$ is reserved for the module assigned to a trivial circle, with basis $\{1, X\}$. Then if $\Cs \subset \A$ consists of $n$ trivial and $m$ essential circles, the module assigned to $\Cs$ by $\Fa$ is written as 
\[
\Fa(\Cs) = A^{\o n} \o V^{\o m}.
\]

\begin{figure}
\begin{subfigure}[b]{.2\textwidth}
\begin{center}
\includegraphics{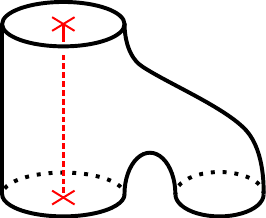}
\end{center}
\caption{Type (I)}\label{fig:type1}
\end{subfigure}
\begin{subfigure}[b]{.2\textwidth}
\begin{center}
\includegraphics{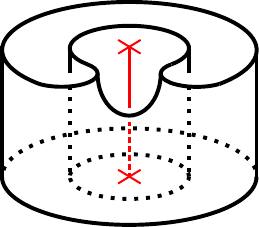}
\end{center}
\caption{Type (II)}\label{fig:type2}
\end{subfigure}
\begin{subfigure}[b]{.2\textwidth}
\begin{center}
\includegraphics{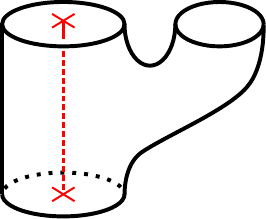}
\end{center}
\caption{Type (III)}\label{fig:type3}
\end{subfigure}
\begin{subfigure}[b]{.2\textwidth}
\begin{center}
\includegraphics{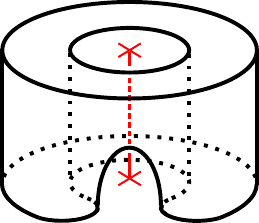}
\end{center}
\caption{Type (IV)}\label{fig:type4}
\end{subfigure}
\caption{Saddles involving essential circles.  The vertical red arc is the central axis of $\A \times I \subset \R^2\times I$.}\label{fig:el saddles}
\end{figure}

Cobordisms between trivial circles are assigned the same map by $\FKh$ and $\Fa$. We record here the maps assigned to the four elementary saddle cobordisms involving at least one essential circle, shown in Figure \ref{fig:el saddles}, using the $V$ notation.

\noindent\begin{minipage}{.5\linewidth}
\begin{align}\label{eq:formula1}
\begin{aligned}
   &V\o A  \rar{\hyperref[fig:type1]{\operatorname{(I)}}} V \\  
   &\vm \o 1 \mapsto \vm \\ 
   &\vp\o 1 \mapsto \vp \\
    &\vm \o X \mapsto 0 \\
    &\vp \o X  \mapsto 0
\end{aligned}
\end{align}
\end{minipage}%
\begin{minipage}{.35\linewidth}
\begin{align}\label{eq:formula2}
  \begin{aligned}
    &V \o V  \rar{\hyperref[fig:type2]{\operatorname{(II)}}} A \\  
    &\vm \o \vm \mapsto 0 \\
    &\vp \o \vm  \mapsto X \\
    &\vm \o \vp  \mapsto X \\
    &\vp \o \vp  \mapsto 0
 \end{aligned}
\end{align}
\end{minipage}

\vskip1em
 
\noindent\begin{minipage}{.49\linewidth}
\begin{align}\label{eq:formula3}
  \begin{aligned}
  &V  \rar{\hyperref[fig:type3]{\operatorname{(III)}}} V \o A \\
&\vm  \mapsto \vm \o X \\
&\vp  \mapsto \vp \o X
  \end{aligned}
\end{align}
\end{minipage} 
\begin{minipage}{.4\linewidth}
\begin{align}\label{eq:formula4}
  \begin{aligned}
    &A  \rar{\hyperref[fig:type4]{\operatorname{(IV)}}} V\o V \\
&1 \mapsto \vp \o \vm + \vm \o \vp   \\
&X  \mapsto 0
  \end{aligned}
\end{align}
\end{minipage}
\vskip1em

Finally, let $D$ be a diagram for an oriented annular link $L$. The construction of $[[D]]$ as described in Section \ref{sec:Khovanov homology} is completely local, and crossings are away from the puncture $\times$. Thus we may view $[[D]]$ as a chain complex over $\BN(\A)$, with the grading shifts $\{-\}$ in $\BN$ rewritten as bigrading shifts $\{-,0\}$ in $\BN(\A)$. Annular link diagrams representing isotopic annular links are related by Reidemeister moves away from the puncture. Therefore the chain homotopy class of $[[D]]$, viewed as a chain complex over $\BN(\A)$, is an invariant of $L$. Applying the annular TQFT $\Fa$, we obtain the annular Khovanov chain complex
\begin{equation}\label{eq:annular chain complex}
CKh_\A(D) := \Fa([[D]]).
\end{equation}
It is an invariant of $L$ up to chain homotopy equivalence.

\subsection{Arc diagrams}
\label{sec:Arc diagrams}
In our resolution diagrams throughout the Khovanov cube for a link diagram $D$, we will use additional decorations to record whether the $0$-smoothing or $1$-smoothing was used at a crossing. These form an enhancement of the \emph{resolution configurations} from \cite[Section 2]{LS}, and they will be crucial to our combinatorial analysis in Section \ref{sec:Path moduli spaces}.

\begin{definition}
\label{def:arc diagram} 
Let $D$ be an $n$-crossing link diagram and let $u\in \{0,1\}^n$ be a vertex. The \emph{arc diagram $D_u$} is obtained from the resolution $D_u$ as follows. Join the two strands at each smoothing by one of two types of arcs, a solid red arc, called a \emph{future arc}, if the $0$-resolution was used, or a dashed blue arc, called a \emph{past arc}, if the $1$-resolution was used. 
\end{definition}

The local picture is shown in Figure \ref{fig:smoothings with arcs}. To explain the terminology, a $0$-smoothing means that the saddle corresponding to that crossing is yet to be performed, while a $1$-smoothing means that the associated saddle has already been performed. For an arc $a$ in an arc diagram, we can obtain a new arc diagram by performing surgery along $a$ and then inserting the dual surgery arc to $a$ and switching its type, as shown in Figure \ref{fig:surgery along an arc}.  We will refer to this process as performing surgery along $a$, with the insertion of the dual arc of opposite type implied.  An edge $D_u \to D_v$ in the cube of resolutions of $D$ thus corresponds to performing surgery along a future arc in $D_u$.  More generally, given any pair of vertices $u,v\in\{0,1\}^n$, the arc diagram $D_v$ can be obtained from $D_u$ by a sequence of such surgeries along past and future arcs in $D_u$.

\begin{figure}
\hspace*{\fill}
\begin{subfigure}[t]{.4\linewidth}
\begin{center}
\includestandalone{two_smoothings_with_arcs}
\end{center}
\caption{Future ($0$-smoothing) and past ($1$-smoothing) arcs.}
\label{fig:smoothings with arcs}
\end{subfigure}
\hfill
\begin{subfigure}[t]{.4\linewidth}
\begin{center}
\includestandalone{surgery_along_an_arc}
\end{center}
\caption{Surgery along a future arc $a$ gives a new arc diagram with the dual to $a$ drawn as a past arc.}
\label{fig:surgery along an arc}
\end{subfigure}
\hspace*{\fill}
\caption{}
\end{figure}

\begin{figure}
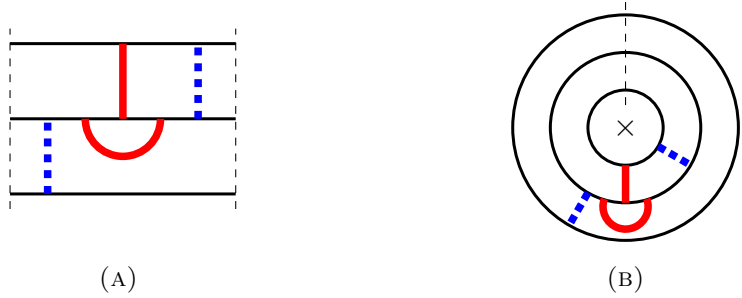

\begin{subfigure}[b]{.4\textwidth}
    \begin{center}
    \includestandalone{arc_diagram_convention_2}
    \end{center}
    \caption{}\label{fig:arc diagram}
\end{subfigure}
\begin{subfigure}[b]{.4\textwidth}
    \begin{center}
    \includestandalone{arc_diagram_convention_1}
    \end{center}
    \caption{}\label{fig:cut open arc diagram}
\end{subfigure}
\caption{A rectangular arc diagram and its annular closure}\label{fig:convention for arc diagrams}
\end{figure}

In what follows, we will depict annular arc diagrams by drawing planar tangles in the rectangle, with the understanding that the annular arc diagram is obtained by connecting the left endpoints to the right. Thus horizontal intervals are segments of distinct essential circles. See Figure \ref{fig:convention for arc diagrams} for an example. We will omit the dashed platforms in the future.

An arc diagram is \emph{connected} if it is connected as a subspace of the plane. Similarly, we will refer to the connected components of arc diagrams.  Note that, since surgery along an arc preserves connectedness, there is a natural correspondence between the connected components of $D_u$ and the connected components of $D_v$ for any two vertices $u,v\in\{0,1\}^n$.

Now let $u=(0,\cdots,0)$ be the starting vertex of the cube, and suppose we have a partition of $D_u$ into disjoint closed (not necessarily connected) components
\[D_u = \Cs_1 \sqcup \cdots \sqcup \Cs_k.\]
In a similar fashion, this partition determines such a partition of the arc diagram $D_v$ at every vertex $v\in\{0,1\}^n$.  Abusing notation slightly, we will use the same notation for these closed components at every vertex:
\[D_v = \Cs_1 \sqcup \cdots \sqcup \Cs_k.\]
Then for any Khovanov generator $x\in\FKh(D_v)$, we will write
\begin{equation}
\label{eq:restrictions of generators}
x=x|_{\Cs_1} \o \cdots \o x|_{\Cs_k} \in \FKh(D_v)
\end{equation}
where $x \vert_{\Cs_i}$ denotes the restriction of the labels in $x$ to the circles in $\Cs_i\subset D_v$.  Note that this is again a slight abuse of notation, since the use of the symbol $\otimes$ in Equation \eqref{eq:restrictions of generators} is, strictly speaking, incompatible with the use of the same symbol in Equation \eqref{eq:generators expanded as labels on each circle}, but this should not cause any confusion in our arguments.

In many places throughout the paper, we will be concerned with certain \emph{subcubes} within the large cube $\{0,1\}^n$.  A $k$-dimensional subcube, denoted $\subcube{k}{v}{u}$, is determined by a pair of vertices $u,v$ such that $u_i\leq v_i$ for all $i$, with precisely $k$ coordinates $i$ where $u_i=0$ and $v_i=1$.  When focusing on such a subcube, we will draw the arc diagram $D_u$ by including only the $k$ future arcs corresponding to these $k$ coordinates, and likewise for any intermediate vertices $u\leq w \leq v$ within the subcube.  See Figure \ref{fig:subcube arc diagram} for an example.

\begin{figure}
\hspace*{\fill}
\begin{subfigure}[b]{.4\textwidth}
    \begin{center}
    \includestandalone{subcube_arc_diagram1}
    \end{center}
    \caption{The arc diagram at the starting vertex $(0,0,0,0)$}\label{fig:subcube arc diagram1}
\end{subfigure}
\hfill
\begin{subfigure}[b]{.4\textwidth}
    \begin{center}
    \includestandalone{subcube_arc_diagram2}
    \end{center}
    \caption{An arc diagram for the subcube $(1,0,0,0) \leq_2 (1,0,1,1)$}\label{fig:subcube arc diagram2}
\end{subfigure}
\hspace*{\fill}
\caption{}\label{fig:subcube arc diagram}
\end{figure}

\subsection{Action of $\sl_2$ on annular Khovanov homology}\label{sec:sl2 action}

It is shown in \cite{GLW} that the annular TQFT $\Fa$ can be viewed as taking values in the category $\gRep(\sl_2)$ of graded $\sl_2$ representations. We first recall some basic definitions. Let $\sl_2$ be the $\Z$-span of the $2 \times 2$ matrices
\[
E = 
\begin{pmatrix}
0 & 1 \\
0 & 0
\end{pmatrix},\,
F = 
\begin{pmatrix}
0 & 0 \\
1 & 0
\end{pmatrix},\,
H = 
\begin{pmatrix}
1 & 0 \\
0 & -1
\end{pmatrix}.
\]
The usual commutator bracket, 
\[
[x,y] = xy-yx,
\]
makes $\sl_2$ into a Lie algebra, with the bracket given on generators by 
\begin{equation}\label{eq:sl2 relations}
[E,F] = H, \;  [H,E] = 2E, \; [H,F] = -2F.
\end{equation}
\begin{remark}
The Lie algebra $\sl_2$ in \cite{GLW} is defined over $\C$. For the purpose of this paper we continue working over $\Z$, since all desired results already hold integrally. 
\end{remark}

Recall that if $M$ is an $\sl_2$-module, then its linear dual 
\[
M^* = \Hom_\Z(M, \Z)
\]
is again an $\sl_2$-module via 
\[
(xf)(m) := -f(xm)
\]
for $x\in \sl_2$, $f\in M^*$, and $m\in M$. Likewise, if $N$ is another $\sl_2$-module, then the tensor product $M\o_\Z N$ inherits an $\sl_2$ action by setting 
\[
x(m\o n) := xm \o n + m \o xn
\]
for $x\in \sl_2$, $m\in M$, and $n\in N$. 

Let $V_1 = \Z \oplus \Z$ denote the fundamental representation of $\sl_2$ with standard basis vectors denoted $\bv_0= (0,1), \bv_1 = (1,0)$, and let $V_1^*$ denote its dual, with dual basis vectors $\bv_0^*, \bv_1^* \in V_1^*$. The action of $\sl_2$ on $V_1$ and $V_1^*$ is summarized in \eqref{eq:sl2 action}. 

\begin{align}\label{eq:sl2 action}
\begin{aligned}
E \bv_0 &= \bv_1 & F \bv_0 &= 0  & H \bv_0 &= - \bv_0\\
E \bv_1 &= 0 & F \bv_1 &= \bv_0  & H \bv_1 &= \bv_1\\
& & & & & \\
E \bv_0^* &= 0 & F \bv_0^* &= -\bv_1^* & H \bv_0^* &=  \bv_0^*\\
E \bv_1^* &= -\bv_0^* & F \bv_1^* &= 0 & H \bv_1^* &= - \bv_1^*
\end{aligned}
\end{align}

Let $\Cs \subset \A$ be a collection of $n$ trivial and $m$ essential circles. There is a natural ordering on the essential circles, starting from the innermost circle (closest to the puncture $\times$). Make
\[
\Fa(\Cs) = A^{\o n} \o V^{\o m}
\]
into an $\sl_2$-module as follows. Every factor $A$ corresponding to a trivial circle is assigned the trivial (zero) two-dimensional representation of $\sl_2$. For a factor $V$ corresponding to the $i$-th essential circle in $\Cs$, identify $V$ with $V_1$ if $i$ is odd via the $\Z$-linear isomorphism \eqref{eq:odd id} and with $V_1^*$ if $i$ is even via \eqref{eq:even id}. 
\begin{align}
V_1 & \to V,\, \bv_0 \mapsto \vm,\, \bv_1 \mapsto \vp, \label{eq:odd id} \\
V_1^* & \to V,\, \bv_0^* \mapsto \vp,\, \bv_1^* \mapsto \vm. \label{eq:even id}
\end{align}
Thus trivial circles are always assigned the trivial two-dimensional representation, and essential circles are assigned $V_1$ and $V_1^*$ in an alternating manner, with the convention that the innermost essential circle is assigned $V_1$. 

\begin{example}\label{ex:action on two essential circles}
The action of $E$ and $F$ on
two essential circles is recorded below.
\begin{align*}
E(\vm \o \vm) &= \vp \o \vm - \vm \o \vp & E(\vp\o \vm) &= - \vp \o \vp \\
E(\vm \o \vp) &= \vp \o \vp & E(\vp \o \vp) &= 0\\
& & & & & & \\
F(\vm \o \vm) &= 0 & F(\vp\o \vm) &=  \vm \o \vm \\
F(\vm \o \vp) &= -\vm \o \vm & F(\vp \o \vp) &= \vm \o \vp - \vp \o \vm
\end{align*}

\end{example}

\begin{remark}
The $\sl_2$ representations $V_1$ and $V_1^*$ are isomorphic via 
\[
\bv_0 \mapsto -\bv_1^*, \hskip1em \bv_1 \mapsto \bv_0^*.
\]
However, the map $V_1 \to V_1^*$ obtained by composing \eqref{eq:odd id} with the inverse of \eqref{eq:even id} is not $\sl_2$-linear; indeed, the identifications \eqref{eq:odd id} and \eqref{eq:even id} equip $V$ with two actions of $\sl_2$ which differ by a sign. 
\end{remark}

By \cite[Lemma 4]{GLW}, maps assigned to cobordisms in $\A \times I$ commute with the action of $\sl_2$. Therefore the annular TQFT $\Fa$ is upgraded to a functor 
\[
\Fa \colon \BN(\A) \to \gRep(\sl_2)
\]
landing in the category of $\Z$-graded $\sl_2$ representations. The annular grading corresponds to the weight space decomposition, and the quantum grading corresponds to the external grading. It follows that $\sl_2$ acts on the annular chain complex \eqref{eq:annular chain complex} and its homology.

Let $\Cs \subset \A$ be a collection of $n$ trivial and $m$ essential circles. We will always assume that the essential circles are ordered from innermost to outermost, so that the $i$-th tensor factor of $V$ in 
\[
\Fa(\Cs) = A^{\o n} \o V^{\o m}
\]
corresponds to the $i$-th essential circle. Any standard generator in $V^{\o m}$ can be written as $v_{s_1} \o \cdots \o v_{s_m}$, with each $s_i \in \{\pm \}$.  Thus any standard generator $x\in \Fa(\Cs)$ is of the form 
\begin{equation}\label{eq:standard generator}
x =y \o v_{s_1} \o \cdots \o v_{s_m}
\end{equation}
where $y\in A^{\o n}$ is a standard generator for the trivial circles in $\Cs$.

By construction, $\sl_2$ weights correspond to the annular grading, so that for a standard generator $x\in \Fa(\Cs)$, $H$ acts as multiplication by $\adeg(x)$,
\begin{equation}\label{eq:H map}
Hx = \adeg(x) x.
\end{equation}
Although the formula \eqref{eq:H map} is straightforward, we can refine it by writing
\begin{equation}\label{eq:H map refined}
Hx = \sum_{i=1}^m s_i y \o v_{s_1} \o \cdots \o v_{s_m}
\end{equation}
using the notation \eqref{eq:standard generator}. The action of $E$ and $F$ can similarly be written as
\begin{align}\label{eq:E and F explicit}
\begin{aligned}
Ex &= \sum_{i=1}^m y \o v_{s_1} \o \cdots\o  E v_{s_i} \o \cdots \o v_{s_m},
\\ Fx & = \sum_{i=1}^m  y \o v_{s_1} \o \cdots\o F v_{s_i} \o \cdots \o v_{s_m}.
\end{aligned}
\end{align}
In the formulas \eqref{eq:E and F explicit}, $E v_{s_i}$ and $F v_{s_i}$ are given by 
\[
E v_{s_i} = 
\begin{cases}
(-1)^{i+1} v_+ & \text{ if } s_i = - \\
0 & \text{ if } s_i = +
\end{cases}
\hskip2em
F v_{s_i} = 
\begin{cases}
(-1)^{i+1} v_- & \text{ if } s_i = + \\
0 & \text{ if } s_i = -
\end{cases}
\]

Letting $C$ denote the $i$-th essential circle in $\Cs$, will say that $E$ \emph{acts on $C$} to mean that we apply $E$ on the $i$-th factor of $V$ in $\Fa(\Cs)$,
\[
x \mapsto y \o v_{s_1} \o \cdots\o  E v_{s_i} \o \cdots \o v_{s_m},
\]
and likewise for $F$. If $J$ denotes one of $E$, $F$, or $H$, we will also say the \emph{$J$ map} to mean the endomorphism of the module assigned to a collection of circles
\[
J\: \Fa(\Cs) \to \Fa(\Cs), \; x\mapsto Jx,
\]
or on the annular chain complex 
\[
J \: CKh_\A(D) \to CKh_\A(D),
\]
which is built out of the $J$ map on each smoothing.

The following is a combination of \cite[Lemma 2]{GLW} and \cite[Lemma 5]{GLW}.

\begin{lemma}\label{lem:E and F duality}
Let $D$ be an annular link diagram. For each smoothing $D_u$, let $\Theta_u \: \Fa(D_u) \to \Fa(D_u)$ be the involution given by $v_{\pm} \mapsto v_{\mp}$ on essential circles and the identity map on the trivial circles. The maps $\Theta_u$ assemble into an isomorphism
\[
\Theta \: CKh_\A(D) \rar{\sim} CKh_\A(D)
\]
of chain complexes of abelian groups, which fits into the commutative diagram
\[
\begin{tikzcd}
CKh_\A(D) \ar[d, "E"'] \ar[r, "\Theta"] & CKh_\A(D) \ar[d, "F"] \\
CKh_\A(D) \ar[r, "\Theta"] & CKh_\A(D).
\end{tikzcd}
\]
\end{lemma}

Therefore the action of $F$ is related to $E$ simply by $F = \Theta E \Theta$. Moreover, $\Theta$ sends generators to generators bijectively, so in our combinatorial analysis we can focus on one of $E$ or $F$. Note however that $\Theta$ does not preserve $(\qdeg, \adeg)$-bidegree. A modified quantum grading is used in \cite{GLW}, there denoted as $j'$. In our notation, $j' = \qdeg - \adeg$.

\section{Khovanov spectra}

\subsection{Using framed flow categories to construct Khovanov spectra}
\label{sec:ffc review}
There is a general method due to Cohen-Jones-Segal \cite{CJS} for constructing a spectrum $\X$ with cellular cochain complex $C^*(\X)$ isomorphic to a given cochain complex $C^*$.  The essential idea is to build a cell complex whose cells correspond to the generators of $C^*$ and whose attaching maps correspond to the differential of $C^*$.  Of course, the various compositions of attaching maps must be homotopically coherent in order for this to make sense.  The data needed to guarantee this coherence is organized in a \emph{framed flow category}, which we will review here.  Most of this material is based upon the corresponding material in \cite{LS,LLS} where the authors successfully constructed such a category (and thus a spectrum) lifting the Khovanov chain complex of a link in $S^3$.

\begin{definition}
\label{def:flow cat}
(\cite[Definition 3.12]{LS})
A \emph{flow category} $\ffc$ is a category consisting of the following data.
\begin{itemize}
    \item A finite set of \emph{objects}, also denoted $\ffc$, together with an indexing map $h:\ffc\rightarrow\Z$.  In analogy with Morse theory, an object $x\in\ffc$ is thought of as a critical point of index $h(x)$ for some Morse function.
    \item For each $x\neq y\in\ffc$, a \emph{morphism space} (usually called a \emph{moduli space}) $\ModSp(y,x)$ which is a $(h(y)-h(x)-1)$-dimensional manifold with corners (see \cite[Section 3.1]{LS})); negative dimensional manifolds are all empty.  In analogy with Morse theory, $\ModSp(y,x)$ is often thought of as the moduli space of flow lines from $y$ down to $x$.  Furthermore, the codimension-$m$ boundary $\de_{[m]}\ModSp(y,x)$ of any such moduli space can be identified with the disjoint union
    \[ \de_{[m]}\ModSp(y,x) \cong \coprod_{(z_1,\cdots,z_m)\in\ffc^m} \ModSp(y,z_m) \times \ModSp(z_m,z_{m-1}) \times \cdots \times \ModSp(z_1,x).\]
    In analogy with Morse theory, $\de_{[m]}\ModSp(y,x)$ is thought of as the $m$-times broken flow lines from $y$ down to $x$.
\end{itemize}
\end{definition}

\begin{remark}
In \cite{LS,LLS}, the authors use a notion of an \emph{$\langle n \rangle$-manifold}, which is a manifold with corners $X$ whose boundary $\de X$ is decomposed in a particular manner into $n$ pieces, $\de X = \de_1 X \cup \cdots \cup \de_n X$. We use square brackets in $\de_{[m]} X$ to denote the codimension $m$ boundary of a manifold with corners in order to distinguish from the earlier notation. 
\end{remark}

\begin{definition}
\label{def:framed flow cat}
(\cite[Definition 3.20]{LS})
A \emph{framed flow category} is a flow category $\ffc$ whose moduli spaces come equipped with \emph{framed neat embeddings} into cornered Euclidean space (see \cite[Definition 3.7]{LLS}) in a manner consistent with the boundary identifications indicated in Definition \ref{def:flow cat}.
\end{definition}

\begin{definition}
\label{def:totalization and refinement}
A framed flow category $\ffc$ determines a cochain complex $\Tot^*(\ffc)$, called the \emph{totalization complex} of $\ffc$,
defined as follows.
\begin{itemize}
    \item The generators of $\Tot^*(\ffc)$ are in bijection with the objects of $\ffc$.
    \item For any $x\in\ffc$, $h(x)$ gives the homological degree of the corresponding generator $x\in \Tot^*(\ffc)$.
    \item For any $x\in\Tot^*(\ffc)$, we have
    \[dx=\sum_{\{y\in\ffc \mid h(y)=h(x)+1\}} \vert \ModSp(y,x) \vert y\]
    where $|\ModSp(y,x)|$ is the signed count of the framed points comprising the 0-dimensional manifold $\ModSp(y,x)$.
\end{itemize}
In analogy with Morse theory, $\Tot^*(\ffc)$ can be viewed as the Morse cochain complex associated to the corresponding Morse function.  Then we say $\ffc$ \emph{refines} a complex $C^*$ if $\Tot^*(\ffc)\cong C^*$.  
\end{definition}

\begin{remark}
Note that the totalization notation in \cite{LLS} is assigned to Burnside functors, while the term \emph{associated cochain complex} is used in \cite{LS} to describe the complex we are calling the totalization in Definition \ref{def:totalization and refinement}
\end{remark}

\begin{theorem}\emph{(\cite[Lemmas 3.24 -- 3.26]{LS})}
\label{thm:framed flow cat gives spectrum}
A framed flow category $\ffc$ gives rise to a cell complex $|\ffc|$ whose cellular cochain complex is isomorphic to $\Tot^*(\ffc)$ up to an overall homological shift.  The stable homotopy type of $|\ffc|$ is an invariant of the given framed flow category $\ffc$.
\end{theorem}

The definitions above are considered in much greater detail in \cite{LS, LLS}, together with the process by which one can build the cell complex $|\ffc|$.  In short, each object $x\in\ffc$ corresponds to a cell of dimension $h(x)+k$ for some fixed offset $k$.  The boundary of such a cell is made up of cornered Euclidean spaces, and the framed embeddings of the moduli spaces into such cornered Euclidean spaces provide homotopically coherent attaching instructions.

We will not need to delve into any greater detail here, because (following \cite{LS, LLS}) our framed flow category will actually be built upon an already existing one.

\begin{definition}
\label{def:cube flow cat}
For any fixed $n$, the \emph{cube flow category} of dimension $n$ is $\ffc_\cube$ whose objects are the vertices $u$ of an $n$-dimensional cube, viewed as elements of $\{0,1\}^n$, and whose non-empty moduli spaces $\cMod$ are permutohedra.  The grading $h\: \ffc \to \Z$ is given by $h(u) = \sum_i u_i$. For more precise details, see \cite[Definition 4.1]{LS} where $\ffc_\cube$ is also described via a specific Morse function on $\R^n$, or \cite[Section 3.4]{LLS} for a combinatorial description. 
\end{definition}
\begin{remark}
In \cite{LS}, the notation $\mathscr{C}$ is used for flow categories.  We choose to use $\ffc$ in order to reserve $\mathscr{C}$ for collections of circles in a diagram.
\end{remark}
\begin{proposition}\emph{(\cite[Section 4]{LS}, \cite[Section 3.6]{LLS})}
\label{prop:cube flow is framed}
The cube flow category $\ffc_\cube$ is a framed flow category.  In particular, its moduli spaces can be neatly embedded (with framing) into cornered Euclidean spaces in a coherent fashion.
\end{proposition}

The proof of Proposition \ref{prop:cube flow is framed} does involve a choice of \emph{sign assignment} for the edges of the cube which correspond to a choice of signs making a cube-shaped chain complex anti-commute so that $d^2=0$ (see \cite[Proposition 4.12]{LS}).  However, different choices of sign assignments give rise to stably equivalent spectra, in a lift of the corresponding argument for chain complexes \cite[Proposition 6.1]{LS}.

In this framework, \cite{LS} show how to refine the Khovanov chain complex of a link $L\subset S^3$ with diagram $D$ into a flow category $\ffc_{Kh}(D)$ that ``trivially covers" $\ffc_\cube$ in the following sense.
\begin{definition}
\label{def:cubical flow cat}
(\cite[Definition 3.21]{LLS})
A \emph{cubical flow category} is a flow category $\ffc$ which comes equipped with a graded functor $f:\ffc\rightarrow\ffc_\cube$ to some fixed cube category such that for all $x,y\in\ffc$, the map $\ModSp_\ffc (y,x) \rightarrow \ModSp_\cube (f(y),f(x))$ is a trivial covering map.
\end{definition}

\begin{proposition}
\label{prop:cubical flow is framed}
Any cubical flow category $\ffc$ can be upgraded to a framed flow category.
\end{proposition}
\begin{proof}
The covering maps $\ModSp_\ffc \rightarrow \ModSp_\cube$ can be composed with the embeddings of $\ModSp_\cube$ into cornered Euclidean space, and all of the framings are then inherited.  These framings also provide a consistent manner in which to separate the components of the trivial covers.  In \cite{LLS} such embeddings are called \emph{cubical neat embeddings}; see Section 3.6 in that reference for more details.
\end{proof}

\begin{remark}
Alternatively, one may follow the arguments in Sections 3.1 and 3.2 of \cite{LS} to upgrade a cubical flow category into a framed flow category.  Our compositions of covering maps and embeddings give rise to \emph{neat immersions} with coherent framings in the language employed in that reference, which can be perturbed to give \emph{neat embeddings}.  See also Section 3.4.1 in that reference.
\end{remark}

\begin{theorem}
\label{thm:Khovanov flow covers cube flow}
(\cite{LS,LLS})
Given a link $L\subset S^3$ with diagram $D$, there is a cubical flow category $\ffc_{Kh}(D)$ which refines the Khovanov complex $CKh^*(D)$.  The resulting suspension spectrum $\X_{Kh}(D)= \Sigma^{\infty} |\ffc_{Kh}(D)|$ is an invariant of the link $L$ up to stable homotopy equivalence, allowing the notation $\X_{Kh}(L)$ for the spectrum.
\end{theorem}

One crucial feature of the Khovanov chain complex used in the proof of Theorem \ref{thm:Khovanov flow covers cube flow} is that, before sign assignment, all signs in the differential are positive.  Thus the commutation of any square face in the Khovanov cube corresponds to a \emph{bijection} of 0-dimensional moduli spaces coming from the edges of the square.  Such bijections correspond to trivial covers of the 1-dimensional permutohedron (which is an interval), allowing $\ffc_{Kh}$ to satisfy the requirements of Definition \ref{def:cubical flow cat} for the 1-dimensional moduli spaces.  In our situation however, this feature will not be present.

\begin{remark}
\label{rmk:Odd Khovanov as signed cover}
It should be noted that the odd Khovanov differential used to define the odd Khovanov complex $CKh_o^*(D)$ in \cite{ORS} does include signs within the cube before the sign assignment.  In \cite{SSS}, the authors used the language of signed correspondences in a suitable Burnside category to resolve this issue.  In terms of our language here, $\ffc_{Kh^o}$ could be thought of as a ``signed" cubical flow category, where the various covering maps carry signs which indicate the need to reverse framing on certain components of the given cover.  Despite this need for signs, it is shown in \cite{SSS} that the square faces in the odd Khovanov cube still correspond to (signed) bijections, which give (signed) trivial covers of intervals as needed.
\end{remark}

The rest of the arguments in \cite{LS,LLS} used to construct $\ffc_{Kh}$ are less immediately relevant, so we simply summarize.  The moduli spaces are built inductively.  The 1-dimensional moduli spaces in $\ffc_{Kh}$ corresponding to square faces involve a choice of bijection between the compositions of $0$-dimensional moduli spaces coming from the two ways to traverse the edges in a square.  For all but one case, these two compositions are either both empty or both consist of one point, so the bijection is determined. The remaining nontrivial case is known as the \emph{ladybug configuration} \cite[Figure 5.1]{LS}, where the two $0$-dimensional moduli spaces each have two points. In this situation, the choice is the \emph{ladybug matching} (\cite[Section 5.4]{LS}).  The 2-dimensional $\ModSp_{Kh}$ are then bounded by some cover of $\de_{[1]}\ModSp_\cube$ (a hexagon, topologically a circle $S^1$).  The hexagon relation (\cite[Section 5.5]{LS}) is the check that this cover is indeed the trivial cover, so that the 2-dimensional $\ModSp_{Kh}$ can be chosen to be a trivial cover of $\ModSp_\cube$.  The higher dimensional $\de\ModSp_\cube$ are all topologically spheres of higher dimension, and so have only trivial covers, allowing the construction of all higher dimensional $\ModSp_{Kh}$ from here.

If we instead let $L\subset \A \times I$ be an annular link with diagram $D$, then the methods in \cite{LS, LLS} extend in a straightforward way to build a stable homotopy refinement of annular Khovanov homology. A formulation using the language of Burnside functors is given in \cite[Section 4.3]{AKW}. The correspondence between Burnside functors and cubical flow categories \cite[Section 4.3]{LLS} yields a flow category $\ffc_\A(D)$ refining $CKh_\A(D)$. The associated suspension spectrum
\begin{equation}
\label{eq:annular spectrum}
\X_\A(D) = \Sigma^\infty|\ffc_\A(D)|
\end{equation}
is an invariant of $L$ up to stable homotopy equivalence. An alternative construction can be found in \cite{LLSCK}.

\begin{remark}
In \cite{LS}, there is a global choice between two types of ladybug matchings, the so-called \emph{left pair} and \emph{right pair}. By \cite[Proposition 6.5]{LS}, the two spectra associated to these choices are stably equivalent.  To the authors' knowledge, it is not known if the annular spectrum is independent of the choice of ladybug matching; in particular the argument in \cite[Proposition 6.5]{LS} does not hold for links in $\A \times I$. For the remainder of this paper we fix one of the two choices. 
\end{remark}

\subsection{Constructing a map between Khovanov spectra}
\label{sec:constructing a map in general} 
Given an annular link diagram $D$, each of the generators $E, F$, and $H \in \sl_2$ yield chain maps $CKh_\A(D) \to CKh_\A(D)$. The goal of this paper is to lift these endomorphisms to maps of spectra. This section discusses our general strategy for lifting maps of cochain complexes to maps of spectra.

\begin{definition}
\label{def:grading shift ffc}
Let $\ffc$ be a flow category with grading $h \: \ffc \to \Z$. The \emph{suspension of $\ffc$}, denoted $\Sigma \ffc$, is the flow category whose underlying category is identical to $\ffc$ but whose grading $\Sigma h$ is shifted up by one, 
\[
\Sigma h (x) = h(x) +1.
\]
\end{definition}

It is straightforward to verify that if $\ffc$ is a framed flow category, then $\Sigma \ffc$ inherits a natural framing, and moreover there is a homeomorphism
\[
\vert \Sigma \ffc \vert \cong \Sigma \vert \ffc \vert.
\]

Consider a map $A^*\xrightarrow{f} B^*$ of cochain complexes, where $A^*$ and $B^*$ are refined by framed flow categories $\ffc(A)$ and $\ffc(B)$ respectively.  In order to lift $f$ to a map of cell complexes $|\ffc(B)|\xrightarrow{F}|\ffc(A)|$ satisfying $F^* = f$ on cohomology (up to a shift), the strategy is to refine the complex $\Cone~f$ into a framed flow category $\ffc(\Cone~f)$.  The key point is that, by construction, $\ffc(\Cone~f)$ will contain complementary downward and upward closed subcategories (\cite[Section 3.4]{LS}) corresponding to $\ffc(A)$ and $\Sigma \ffc(B)$ respectively. The downward closed subcategory yields a subcomplex $\vert \ffc(A) \vert \subset \vert \ffc(\Cone~f) \vert$ whose quotient complex is canonically identified with $\vert \Sigma \ffc(B) \vert = \Sigma\vert \ffc(B)\vert$. The Puppe sequence yields 
\[
\vert \ffc(A)\vert \hookrightarrow \vert \ffc(\Cone~f)\vert \twoheadrightarrow \vert \ffc(\Cone~f)\vert / \vert \ffc(A)\vert = \Sigma \vert \ffc(B) \vert \xrightarrow{\star}  \Sigma \vert \ffc(A) \vert \to \cdots 
\]
where $\star$ is the Puppe map. After shifting degrees, the Puppe map is the desired $F$.

Now we note that, if $A^*=B^*=CKh^*(D)$ for a link diagram $D$ having $n$ crossings and $f$ is assembled from maps $\Fa(D_u) \to \Fa(D_u)$ on each smoothing, then the complex $\Cone~f$ takes the shape of an $(n+1)$-dimensional cube $\cube$.  If in addition the map $f$ contains no signs, then one can seek to build a cubical flow category $\ffc_{Kh}(\Cone~f)$ covering $\ffc_{\cube}$ as before.  Assuming this can be done, the resulting Puppe map $F$ induces a map of spectra $\X(D)\xrightarrow{\Sigma^\infty F} \X(D)$ which lifts $f$ as desired.  All of this is illustrated by the exact triangle for a crossing (\cite[Theorem 2]{LS}) which is also used in \cite[Section 3.3]{LS2} to define the map induced by a saddle cobordism between links. 

In our situation, $D$ is an annular link and $f$ is one of $E, F$, or $H$. Moreover, as illustrated in Section \ref{sec:cubic illustration}, these maps contain signs leading to cancellations, which prevents $\ffc(\Cone~f)$ from being made into a cubical flow category. Nevertheless, we will refine $\Cone~f$ and lift $E, F$, and $H$ by using the above procedure.

\begin{remark}
\label{rmk:old Burnside version}
All of the constructions thus far can be converted into the language of functors from cube categories to Burnside categories.  This is the framework introduced in \cite{LLS} to investigate disjoint unions and connect sums of links, and it provides a convenient combinatorial approach to constructing stable homotopy types, as in \cite{SSS} for the odd Khovanov spectrum.  In this language, chain maps are lifted as natural transformations of Burnside functors.  In a previous paper on annular Khovanov spectra \cite{AKW}, we used (an equivariant version of) this approach to refine the quantum annular Khovanov homology constructed in \cite{BPW} and the action coming from the generators $K,K^{-1} \in U_q(\sl_2)$ on quantum annular Khovanov homology.  However, as we shall see, this viewpoint is not sufficient to refine the action of the remaining generators of $U_q(\sl_2).$  See \cite[Section 8]{AKW} for a related discussion. 
\end{remark}

\section{New features: cancellations and topology of moduli spaces}
\label{sec:cubic illustration}

This section considers a particular example to illustrate new complexities that appear in the analysis of the $\sl_2$ action. These features motivate the introduction of new techniques in follow-up sections. 

Consider the $2$-crossing annular knot diagram  and the associated surgery arc diagram (see Section \ref{sec:Arc diagrams}) in Figure \ref{fig:example}. The cube of resolutions for this link diagram (to be more precise, a square since there are two crossings) is given in Figure \ref{fig:example square}; the coordinate directions are labeled in accordance with the numbering of the surgery arcs in Figure \ref{fig: example arc diagram}.

\begin{figure}[ht]
\hspace*{\fill}
\begin{subfigure}[b]{.4\textwidth}
    \begin{center}
    \includestandalone{cubic_knot_diagram}
    \end{center}
    \caption{An annular link diagram} \label{fig:example diagram}
\end{subfigure}
\hfill
\begin{subfigure}[b]{.4\textwidth}
    \begin{center}
    \includestandalone{cubic_surgery_diagram}
    \end{center}
    \caption{The corresponding arc diagram} \label{fig: example arc diagram}
\end{subfigure}
\hspace*{\fill}
\caption{}\label{fig:example}
\end{figure}

\begin{figure}
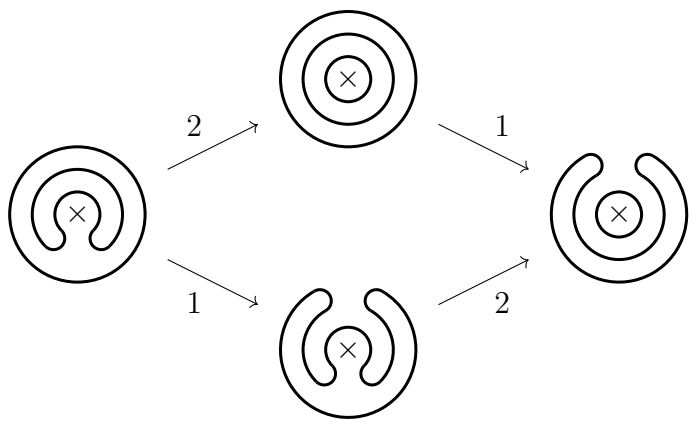

\begin{center}{\includestandalone{cubic_cube_of_resolutions}}
\end{center}
\caption{The cube of resolutions for the link in Figure \ref{fig:example diagram}}\label{fig:example square}
\end{figure}

\begin{figure}
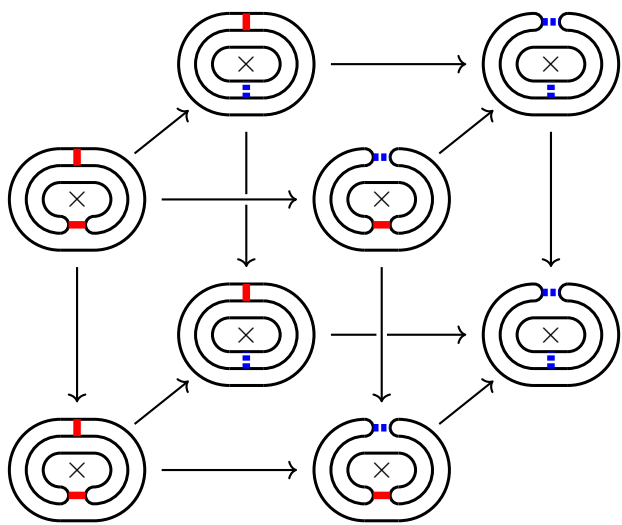

\begin{center}{\includestandalone{cubic_cube_of_resolutions_prime}}
\end{center}
\caption{The cube-shaped complex associated to $E$}\label{fig:example cube}
\end{figure}

The cube-shaped chain complex ${\rm Cone}(E)$ is shown in Figure \ref{fig:example cube}. This cube consists of two copies (pictured horizontally) of the cube of resolutions from Figure \ref{fig:example square}; the action of $E$ takes place along the vertical axis.  Each diagram in the cube is endowed with past/future arcs as in Definition \ref{def:arc diagram}.

Consider a particular generator, $1\otimes v_-$, at the upper left vertex of the cube.  We use the computations from Section \ref{sec:sl2 action} to arrive at the following diagram which summarizes the cube edge maps (we suppress the tensor product notation for the $v$ factors, for brevity).

\begin{equation} \label{cubic cube}
\begin{tikzcd}[row sep=scriptsize,column sep=scriptsize]
& v_+v_-v_-+ v_-v_+v_-  \arrow[rr, "1"] \arrow[dd, "E", pos=0.2] & & X\otimes v_-  \arrow[dd, "E"] \\ 
1\otimes v_- \arrow[rr, crossing over, "1", pos=0.7] \arrow[dd, "E"] \arrow[ur, "2"] & & v_- \arrow[ur, "2"]  \\
 & v_+v_-v_++v_-v_+v_+ \arrow[rr, "1", pos=0.8] & & X \otimes v_+ \\ 
 1\otimes v_+ \arrow[rr, "1"] \arrow[ur, "2"] & & v_+ \arrow[ur, "2"] \arrow[uu, crossing over, leftarrow, "E"', pos=0.8]\\
\end{tikzcd}
\end{equation}

A crucial feature in this example is {\em cancellation}. It is not seen in the result of the algebraic calculation in the cubical diagram above, and we will explain it now in more detail. The generator $1\otimes v_-$ is sent by the edge map labeled $2$ to $v_+v_-v_-+ v_-v_+v_-$. The value of the $E$ map on these two summands is given as follows:

\begin{equation} \label{1st eq} E(v_+v_-v_-)=  v_+v_-v_+ -v_+v_+v_-
\end{equation}
\begin{equation} \label{2nd eq} E(v_-v_+v_-)=  v_-v_+v_+ + v_+v_+v_- 
\end{equation}

Considering the composition $E\circ 2$, the terms $-v_+v_+v_-$, $v_+v_+v_-$ in equations (\ref{1st eq}), (\ref{2nd eq}) cancel; this is an instance of   cancellation that is the root of the complexity of this problem as we explain next.

\begin{figure}[ht]
\begin{center}{\includestandalone{cubic_2_morphism}}
\end{center}
\caption{}\label{fig:example 2 morphism}
\end{figure}

A relevant part of the moduli spaces is illustrated in Figure \ref{fig:example 2 morphism}. The vertical dashed lines are labeled by compositions of edges of length three in the cube, and the corresponding $1$-morphisms are indicated by thickened dots. For example, the effect of the edge path $E21$ on the generator $1\otimes v_-$ is
$$ 1\otimes v_-\overset{1}\longmapsto v_-\overset{2}\longmapsto X\otimes v_-\overset{E}\longmapsto X\otimes v_+$$
The resulting single flow line is represented by a dot above the label {$E21$} in Figure \ref{fig:example 2 morphism}. There is a single flow line associated to the edge path $E12$ as well, since the summand $v_+v_-v_-$ is sent to zero by the edge map labeled 1. The interesting case in this example is the edge path $1E2$: there are three flow lines
$$ 1\otimes v_-\longmapsto v_+v_-v_- \longmapsto v_+v_-v_+\longmapsto X\otimes v_+$$
$$ 1\otimes v_-\longmapsto v_+v_-v_-\longmapsto -v_+v_+v_-\longmapsto - X\otimes v_+$$
$$ 1\otimes v_-\longmapsto v_-v_+v_-\longmapsto v_+v_+v_-\longmapsto X\otimes v_+$$
represented by the three dots above the label $1E2$. (Note that $v_-v_+v_+$ is sent to zero by the map $1$, so this term does not contribute a flow line.) The rest of the edge paths give rise to a single flow line.

The vertical strips between dashed lines in Figure \ref{fig:example 2 morphism} represent $2$-morphisms, or in other words $1$-dimensional moduli spaces corresponding to square faces in the cube ${\rm Cone}(E)$.
Restricting to the square face in the back of the cube (\ref{cubic cube}), the flow line 
\[v_-v_+v_-\longmapsto X\otimes v_-\longmapsto X\otimes v_+\]
along the edge path $E1$ corresponds to the flow line 
\[v_-v_+v_-\longmapsto v_+v_+v_-\longmapsto X\otimes v_+\]
along $1E$. This explains the matching of the flow lines in the $2$-morphism $E12-1E2$, i.e. the interval that goes across that vertical strip in Figure \ref{fig:example 2 morphism}.

On the other hand, $v_+v_-v_-$ is sent to zero by the map $1$, and thus it does not contribute to a flow line along $E1$. The same generator $v_+v_-v_-$ is sent by the map $E$ to a pair of generators in equation (\ref{1st eq}) which then cancel under the map $1$. Thus the $2$-morphism $E12-1E2$ has to ``match'' the empty set at $E12$ to two canceling flow lines at $1E2$. It is natural to represent this pair of flow lines by oppositely framed points in the interval labeled by $1E2$, and to connected them by a framed ``turnback'' in the strip $E12-1E2$. 

A similar analysis along the square face on the left of the cube forces the matching of the flow lines under the $2$-morphism $1E2-12E$, creating the ``cubic'' shape of the $1$-dimensional moduli space indicated in Figure \ref{fig:example 2 morphism}. 

This example shows that bijections across square faces, used in \cite{LS} in the setting of Khovanov homology, are insufficient in the presence of cancellations. On the other hand, $2$-dimensional faces of the cube force the matching of flow lines associated to edge paths by framed $1$-dimensional moduli spaces. As indicated in the introduction, in general building higher-dimensional moduli spaces is a rather non-trivial problem. Our solution in this paper is based on representing the moduli spaces as codimension-$1$ submanifolds, as illustrated in Figure \ref{fig:example 2 morphism}. A subtlety in this approach is that such a codimension-1 position of the flow lines ($0$-dimensional moduli spaces) in the interval has to be prescribed for all edge paths in a given cube-shaped chain complex ${\rm Cone}(E)$. Indeed, these flow lines for the edge paths of arbitrary lengths are a part of the inductive construction of higher-dimensional moduli spaces, based on Definition \ref{def:flow cat}. Moreover, given any two generators at two arbitrary vertices $u,v$ in the cube, the positions in the interval for all edge paths connecting $u,v$ in the cube have to be consistent, in the sense that they can be connected by the prescribed framed $1$-dimensional moduli spaces in codimension-1 without intersections. The analysis that enables a solution to this problem is the content of
Sections \ref{0 dim sec} - \ref{sec:Path moduli spaces}. 

The construction of higher-dimensional moduli spaces in Section \ref{sec:m-1 to m}
proceeds by induction. For example, in the basic case considered above, the $2$-dimensional moduli space has to be constructed for the two chosen generators $1\otimes v_-$, $X\otimes v_+$ in the cube (\ref{cubic cube}). Identifying the leftmost and rightmost vertical intervals in Figure \ref{fig:example 2 morphism} gives a boundary condition for the $2$-dimensional moduli space in a hexagon times an interval, see Figure \ref{fig:Hexagon times I}. In general, the higher-dimensional inductive step amounts to finding a null-cobordism in a permutohedron times an interval $I$, as explained in Section \ref{sec:m-1 to m}.

\section{Constructing the framed flow category: The set-up and the base case} \label{0 dim sec}

We begin by fixing an $n$-crossing annular link diagram $D$, and we let $J$ denote one of the maps $E$, $F$, or $H$.  As discussed in Section \ref{sec:constructing a map in general}, the chain map $J$ on the annular Khovanov complex $CKh_\A(D)$ gives rise to $\Cone~J$, a complex in the shape of an $(n+1)$-dimensional cube $\cube$.  Our goal is to build a framed flow category $\ffcJ(D)$ (we will often omit the diagram $D$ from the notation) that refines $\Cone~J$ with the help of $\ffc_\cube$, the framed flow category for the cube $\cube$.

\subsection{Definitions and the overall strategy}
\label{sec:Definitions and overall strategy}
To begin with we set some notation and basic definitions.  The vertices of the cube $\cube$ will use labels $u,v,w$, with \[u=(u',i)\in\{0,1\}^n \times \{0,1\}.\]
Here we are viewing the first $n$ coordinates $u'$ as indicating where in $CKh_\A(D)$ the corresponding vertex lies, while the last coordinate $i$ indicates whether we are in the domain or range of the $J$ map.

\begin{definition}
\label{def:resolution diagram}
Given a vertex $u=(u',i)\in\cube$, the \emph{resolution diagram} $D_u$ is the resolution of our link diagram $D$ at the vertex $u'$ in the $CKh_\A(D)$.  We further define $\KG{u}$ to be the set of Khovanov generators for $D_u$; see Equation \eqref{eq:generators expanded as labels on each circle}.
\end{definition}

\begin{definition}
An \emph{ $m$-dimensional sub-cube } is a choice of two vertices 
\[\scube\] 
such that $u$ and $v$ match in all but $m$ indices where $u$ has entry zero and $v$ has entry one.  A $1$-dimensional sub-cube is also called an \emph{edge}.  We say that a subcube $\scube$ \emph{involves $J$} if $u=(u',0)$ and $v=(v',1)$ (so that some edges corresponding to the $J$-map are included).
\end{definition}

\begin{definition}
For any sub-cube $\scube$, the \emph{cube moduli space} $\cMod$ is by definition the $(m-1)$-dimensional permutohedron (see Definition \ref{def:cube flow cat}).  This is a manifold with corners (for example, if $m=3$ this is a hexagon) that is topologically equivalent to $D^{m-1}$, whose codimension $i$ boundary (for all $i=1,\dots,m-1)$ is identified with $i$-fold ``composition'' moduli spaces
\[\de_{[i]} \cMod \cong \coprod \cubeMod{v}{w_i} \times \cdots \cubeMod{w_1}{u},\]
where the disjoint union is taken over all sequences of vertices $u\leq w_1 \leq \cdots \leq w_i \leq v$.  For more details, see \cite[Definition 4.1]{LS} and \cite[Section 3.4]{LLS}.
\end{definition}

\begin{remark}
\label{rmk:direction reversal}
Note the reversal of direction. For vertices $u$ and $v$ with $u\leq_1 v$, there is an edge from $u$ to $v$ in the cube of resolutions. However, moduli spaces are directed from $v$ to $u$.
\end{remark}

As described in Section \ref{sec:ffc review}, the original construction \cite{LS} used trivial coverings of the moduli spaces in the cube flow category (see Definition \ref{def:cubical flow cat}).  The example of Section \ref{sec:cubic illustration} shows that this framework alone is insufficient for building our desired $\ffcJ$.  Instead, we will need to keep track of signs of the $J$ map with the help of framings.  First we give ourselves an extra dimension to allow for meaningful framings to occur.

\begin{definition}
\label{def:cube ambient space}
For any sub-cube $\scube$, the \emph{cube ambient space} $\AmbS$ is defined depending on whether or not the $J$ map was involved in $\scube$.
\begin{equation}\label{eq:amb space}
\AmbS := \begin{cases}
\cMod & \text{if $J$ is not involved}\\
\cMod \times (0,1) & \text{if $J$ is involved}
\end{cases}
\end{equation}
In either case, $\AmbS$ is a manifold with corners whose boundary structure is determined by that of $\cMod$.  In our figures throughout the paper, we will draw the $(0,1)$-direction horizontally.
\end{definition}

Next we introduce the types of maps that we will be using in place of the trivial covering maps utilized in \cite{LS} via Definition \ref{def:cubical flow cat}.

\begin{definition}
\label{def:thick embedding}
A \emph{thick embedding} is a framed map $f=e\circ t$ where $t$ is a trivial covering and $e$ is an embedding.  We will denote thick embeddings with the symbol
$\thickemb$.  If $t$ is the trivial identity cover (i.e. the identity map), we call $f$ a \emph{thin embedding} and denote with the usual symbol $\hookrightarrow$.  Thus a thin embedding is just a framed embedding.
\end{definition}

\begin{remark}
Note that according to the definition, different sheets in a thick embedding may have different framings. In codimension $1$, such framings are either identical or are reversed. As we show later in the paper, the former case occurs for $J=E, F$ and the latter case occurs for $J=H$.
\end{remark}

With all of this notation in place, we can state our main construction as a theorem.

\begin{theorem}
\label{thm:ffcJ} 
Fix an $n$-crossing annular link diagram $D$ and an $\sl_2$ generator $J=E,F,$ or $H$ giving rise to a chain map on annular Khovanov complexes with mapping cone $\Cone~J$.  Let $\cube$ denote the $(n+1)$-dimensional cube, with corresponding framed flow category $\ffc_\cube$.  Then there exists a flow category $\ffc_J=\ffc_J(D)$ together with a ``thick embedding functor" $f:\ffcJ\rightarrow \ffc_\cube$ satisfying the following properties.
\begin{enumerate}
\item The objects of $\ffc_J$ are in bijection with the Khovanov generators throughout the complex $\Cone~J$:
\[\Obj(\ffcJ):=\coprod_{u\in\cube} \KG{u},\]
with grading given by homological degree in $\Cone~J$, and graded functor $f:\Obj(\ffcJ)\rightarrow\Obj(\ffc_\cube)$ being the obvious covering map of sets
\[f(x)= u \quad \Leftrightarrow \quad x\in \KG{u}.\]

\item For any generators $x\in\KG{u},y\in\KG{v}$ such that $u\not\leq v$, the moduli space $\MspJxy$ is the empty set.

\item
\label{item:ffcj boundary definitions}
Given a sub-cube $\scube$ and a pair of generators $x\in \KG{u}, y\in \KG{v}$, the moduli space $\MspJxy$ is an $(m-1)$-dimensional manifold with corners, whose codimension $i$ boundary is identified with $i$-fold compositions
\[
\de_{[i]}\MspJxy \cong \coprod \MspJ{y}{z_i} \times \cdots \times \MspJ{z_1}{x}
\]
where the disjoint union is taken over all tuples of generators $(z_1,\dots,z_i)\in \KG{w_1}\times \cdots \times \KG{w_i}$ for sequences of vertices $u\leq w_1 \leq \cdots \leq w_i \leq v$ in the subcube $\scube$.

\item
\label{item:ffcj thick embedding condition}
Each moduli space $\MspJxy$ of item \eqref{item:ffcj boundary definitions} comes equipped with a thick embedding
\[f:\MspJxy\thickemb \AmbS\]
such that, for any tuple $(z_1,\dots,z_i)\in \KG{w_1}\times \cdots \times \KG{w_i}$ giving rise to a component of $\de_{[i]}\MspJxy$, the following diagram of thick embeddings (using the product framing as necessary) commutes

\begin{equation}
\begin{aligned}
\label{eq:respecting boundary structure}
\begin{tikzpicture}[x=1.5in,y=-.5in]

\node(A) at (0,0) {$\MspJxy$};

\node(B) at (0,1) {$\de_{[i]} \MspJxy$};

\node(C) at (0,2) {$\MspJ{y}{z_i} \times \cdots \times \MspJ{z_1}{x}$};

\node(D) at (2,0) {$\AmbS$};

\node(E) at (2,1) {$\de_{[i]} \AmbS$};

\node(F) at (2,2) {$\AmbSp{v}{w_i}\times\cdots\times \AmbSp{w_1}{u}$};

\xthickembtikz{A}{D}{f}
\thickembtikz{B}{E}
\xthickembtikz{C}{F}{f\times \cdots \times f}

\draw[right hook->] (C) -- (B);
\draw[right hook->] (B) -- (A);
\draw[right hook->] (F) -- (E);
\draw[right hook->] (E) -- (D);

\end{tikzpicture}
\end{aligned}
\end{equation}
 
\item
\label{item:differential condition}
In the case $m=1$, the $0$-dimensional moduli spaces correspond to the components of the differential in $\Cone~J$.
\end{enumerate}
\end{theorem}

Item \eqref{item:ffcj boundary definitions} ensures that the moduli space boundaries of $\ffcJ$ are arranged in the same way as those of $\ffc_\cube$, while item \eqref{item:ffcj thick embedding condition} ensures that this arrangement is respected by the thick embeddings as well.  Meanwhile item \eqref{item:differential condition} refers to the edges of the cube, where the moduli spaces are all thickly embedded points in either $D^0$ or $D^0\times(0,1)$ corresponding to either the annular Khovanov differential in $CKh_\A(D)$ or the $J$ map, respectively.

Recall the annular Khovanov spectrum $\X_\A(D)$ from Equation \eqref{eq:annular spectrum} and the surrounding discussion.  

\begin{corollary}
\label{cor:ffcj refines Cone J}
The flow category $\ffcJ(D)$ defined above can be upgraded to a framed flow category which refines the complex $\Cone~J$, giving rise to a map of spectra $\J:\X_\A(D) \rightarrow \X_\A(D)$ whose induced map $\J^*$ on cohomology is equal to the action of $J$. 
\end{corollary}
\begin{proof}
The cube flow category $\ffc_\cube$ is a framed flow category, so each $\cMod$ comes equipped with a framed neat embedding $\cMod\hookrightarrow \cEuc$ where $\cEuc$ denotes a suitable cornered Euclidean space (see \cite[Section 3.6]{LLS}).  As in Definition \ref{def:cube ambient space}, we set
\begin{equation}\label{eq:Euc space}
\EucS := \begin{cases}
\cEuc & \text{if $J$ is not involved}\\
\cEuc \times \R & \text{if $J$ is involved}
\end{cases}
\end{equation}
We then have embeddings $\AmbS\hookrightarrow\EucS$, where $(0,1)$ embeds into $\R$ in the standard way when $J$ is involved.  We then compose with our thick embeddings provided by Item \eqref{item:ffcj thick embedding condition} to give
\[\MspJxy \thickemb \AmbS \hookrightarrow \EucS.\]
Just as in Proposition \ref{prop:cubical flow is framed} for the cubical flow categories of \cite{LLS}, the framing inherited by this composition provides us with a consistent way to separate the various components of each trivial cover, giving us framed neat embeddings $\MspJxy \hookrightarrow \EucS$ which respect the necessary boundary conditions as required (note that crossing with $\R$ in Equation \eqref{eq:Euc space} does not affect the cornered structure of $\cEuc$).  Item \eqref{item:differential condition} then ensures our resulting framed flow category $\ffcJ(D)$ has refined $\Cone~J$.  This gives rise to the desired map as described in Section \ref{sec:constructing a map in general}.
\end{proof}

The rest of this section together with Sections \ref{sec:m=2}, \ref{sec:Path moduli spaces}, and \ref{sec:m-1 to m}
are devoted to proving Theorem \ref{thm:ffcJ}.  Items \eqref{item:ffcj boundary definitions} and \eqref{item:ffcj thick embedding condition} of the theorem are the most crucial aspects of the construction, and the various thickly embedded moduli spaces will be built inductively.  If the sub-cube $u\leq_m v$ does not involve the $J$ map, then it corresponds to a sub-cube in the cube of resolutions of $D$. In this case the thick embedding $f:\MspJxy \thickemb \AmbS=\cMod$ will be a trivial cover built exactly as in \cite{LS} (based on the annular differential rather than the traditional Khovanov differential). If $u\leq_m v$ involves the $J$ map, then $\MspJxy$ will be thickly embedded as a codimension $1$ submanifold of $\AmbS$.  We handle the base case below.

\subsection{0-dimensional moduli spaces for edges}
\label{sec:0-dim moduli spaces}
Fix $u\leq_1 v$, $x\in \KG{u}$, and $y\in \KG{v}$. We begin by defining the $0$-dimensional moduli space $\MspJxy$.

Suppose first that $u\leq_1 v$ is not the $J$ edge. Then it corresponds to an edge $u\to v$ in the cube of resolutions of $D$. Let 
\[
\psi_{u,v} \: \Fa(D_u) \to \Fa(D_v)
\]
denote the map assigned to the saddle cobordism $D_u \to D_v$, and write
\[
\psi_{u,v}(x) = \sum_{z\in \KG{v}} k_z z
\]
where each $k_z \in \{0,1\}$. Say that $y$ \emph{appears in $\psi_{u,v}(x)$} if the above coefficient $k_y$ is nonzero. If $y$ appears in $\psi_{u,v}(x)$, then $\MspJxy$ is defined to consist of a single element, which we denote
\[
x\to y.
\]
If $y$ does not appear in $\psi_{u,v}(x)$, set $\MspJxy = \varnothing$.

Now suppose that $u \leq_1 v$ is the $J$ map. Let $s$ denote the number of essential circles in $D_u = D_v$. We distinguish two cases. First, if $J$ is one of $E$ or $F$, again write 
\[
Jx = \sum_{z\in \KG{v}} \l_z z
\]
where now each $\l_z\in \{-1,0,1\}$, and say that $y$ appears in $Jx$ if $\l_y\neq 0$. If $y$ appears in $Jx$, then it is obtained from $x$ by the action of $J$ on an essential circle in $D_u$. Let $1\leq i \leq s$ denote the position of this essential circle, and define $\MspJxy$ to consist of a single element, which we denote 
\[
x \xrightarrow{i} y.
\]
If $y$ does not appear in $Jx$, then set $\MspJxy = \varnothing$.

Finally, if $J=H$, set $\ModSp_H(y,x)=\varnothing$ unless the labels ($v_\pm,1,X$) on each circle in $x$ and $y$ match, in which case  $\ModSp_H(y,x)$ contains a point for every essential circle in $D_u=D_v$. These $s$ points in $\ModSp_H(y,x)$ will be denoted 
\[
x \xrightarrow{i} y
\]
for $1\leq i \leq s$.

We now describe the thick embeddings, which again depend on the type of edge $u\leq_1 v$. If $u\leq_1 v$ is not the $J$ map, then $\AmbS$ is a point, so the map is unique. 

Suppose now that $u\leq_1 v$ corresponds to the $J$ map. Let $s$ denote the number of essential circles in $D_u = D_v$. Recall that $\AmbS = \cMod \times (0,1) \cong (0,1)$. We define our thick embedding $\MspJxy\xthickemb{f}\AmbS$ by setting\footnote{The important data here is the ordering of points in the interval; a specific formula is given just to fix the convention.}
\begin{equation}
\label{eq:m=1 J-embedding}
f(x\xrightarrow{i}y):=\frac{i}{s+1}\in(0,1),
\end{equation}
with framing determined as follows.  If $J$ is either $E$ or $F$, then the framing is given by the sign of $(-1)^{i+1}$.  If $J=H$, then the framing is $\pm 1$ if $x$ (and thus $y$) is labeled by $v_{\pm }$ on the $i$-th essential circle in $D_u$.

Note that all of these are actually thin embeddings when the subcube is $1$-dimensional.  It follows from the construction of the spaces and thin embeddings that item \eqref{item:differential condition} of Theorem \ref{thm:ffcJ} is clearly satisfied.

\begin{remark}
Recall from \eqref{eq:H map} that $Hx = \adeg(x) y$, where the circle labels in $y$ match those in $x$. It may seem natural to define the $0$-dimensional moduli space $\ModSp_H(y,x)$ to contain $\vert\adeg(x) \vert$ points, and define the thick embedding to be constant, with framing given by the sign of $\adeg(x)$. In the above definition of $\ModSp_H(y,x)$ we are motivated by the refined formula \eqref{eq:H map refined} in which $H$ acts on each essential circle. 
\end{remark}

\begin{example}
\label{ex:0-dim ex}
We illustrate the $0$-dimensional moduli spaces for an edge corresponding to various $J$ maps when $D_u = D_v$ consists of two essential circles. The thin horizontal segment is the interval $\AmbS \cong (0,1)$, the dot $\bullet$ is the image of $\MspJxy$, and the arrows indicate the framing. 

Moduli spaces $\ModSp_E(y,x)$ are depicted below.

\begin{equation*}
\includestandalone{moduli_space_ex_E}
\end{equation*}
For $J=F$ the situation is similar; see also Lemma \ref{lem:E and F duality} and the discussion surrounding it. 

Moduli spaces $\ModSp_H(y,x)$ are shown below. There are four generators of $\Fa(D_u)$, hence four nonempty moduli spaces, each containing two elements of the form $x \rar{i} y$ for $i=1, 2$. The numbers above the dots $\bullet$ in the following diagram indicate the image of the point $x \rar{i} y$. 

\begin{equation*}
\includestandalone{moduli_space_ex_H}
\end{equation*}

The moduli spaces for pairs $(x,y)$ which are not depicted above are all empty. 
\end{example}

\section{Constructing the 1-dimensional moduli spaces for square faces}
\label{sec:m=2}

We now define the $1$-dimensional moduli spaces and their thick embeddings, with an eye towards the inductive step in Section \ref{sec:m-1 to m}. Although the combinatorial analysis in this section can be simplified by the more general considerations in Section \ref{sec:Path moduli spaces}, we present a ``hands-on'' approach to illustrate the notions in Section \ref{sec:Definitions and overall strategy}. 

Fix $u\leq_2 v$, $x\in \KG{u}$, and $y\in \KG{v}$.  The moduli space $\MspJxy$ will be a disjoint union of closed intervals, so we need to specify which points in $\de \MspJxy$ bound an interval in $\MspJxy$.  Item \eqref{item:ffcj boundary definitions} of Theorem \ref{thm:ffcJ} requires that these boundary points come from products corresponding to composition along the edges of the subcube.  Letting $w_1, w_2$ denote the two vertices with $u\leq_1 w_1, w_2 \leq_1 v$, we introduce the notations
\begin{align}
\begin{split}
    \pathPts{w_1}{\MspJxy} & := \coprod_{z\in \KG{w_1}} \ModSp_J(y,z) \times \ModSp_J(z,x),\\
    \pathPts{w_2}{\MspJxy} & := \coprod_{z\in \KG{w_2}} \ModSp_J(y,z) \times \ModSp_J(z,x)
\end{split}
\label{eq:m=2 path moduli spaces}
\end{align}
for these two composition moduli spaces - one for the edge-path passing through $w_1$, and the other for the edge-path passing through $w_2$.  Thus we have
\[\de\MspJxy = \pathPts{w_1}{\MspJxy} \,\coprod\, \pathPts{w_2}{\MspJxy}.\]
 We write elements of each $\pathPts{w_i}{\MspJxy}$ as
\begin{equation}
\label{eq:elements in composition moduli spaces}
(x\rightarrow z \rightarrow y):= (z\rightarrow y , x\rightarrow z)\in \ModSp_J(y,z) \times \ModSp_J(z,x)
\end{equation}
for various $z\in\KG{w_i}$.  If one of the edges of $\scube$ corresponds to the $J$ map, we decorate the corresponding arrow in $(x\rightarrow z \rightarrow y)$ according to which essential circle was acted upon as in Section \ref{sec:0-dim moduli spaces}.  Each composition moduli space comes equipped with a thick product embedding
\begin{equation}
\label{eq:m=2 boundary embeddings}
\pathPts{w_i}{\MspJxy} \thickemb \AmbSp{v}{w_i}\times \AmbSp{w_i}{u}\cong
\begin{cases}
D^0 & \text{if $J$ is not involved}\\
(0,1) & \text{if $J$ is involved}
\end{cases}
\end{equation}
induced by the definitions in Section \ref{sec:0-dim moduli spaces}.  Our goal now is to construct a moduli space with thick embedding
\[\MspJxy \thickemb \AmbS \cong
\begin{cases}
D^1 & \text{if $J$ is not involved}\\
D^1\times (0,1) & \text{if $J$ is involved}
\end{cases}
\]
which `fills' the thick embeddings \eqref{eq:m=2 boundary embeddings} in the sense that Equation \eqref{eq:respecting boundary structure} is satisfied.  Naturally, we split into cases depending on whether or not $J$ was involved.

\subsection{Square faces involving only the annular Khovanov differential}
\label{sec:m=2 J not involved}
A subcube $u\leq_2 v$ that does not involve the $J$ edge corresponds to a square face in the cube of resolutions of $D$, and in this case we proceed exactly as in \cite{LS}. From \eqref{eq:m=2 path moduli spaces} and the analysis in \cite{LS}, we see that 
\[ \size{\pathPts{w_1}{\MspJxy}} = \size{\pathPts{w_2}{\MspJxy}} =: N \in \{0,1,2\}\]
where the number of points $N$ depends on the combinatorics of the square face $D_u\rightarrow D_v$.  

Following \cite{LS}, we define
\[
\ModSp_J(y,x) := \coprod_S \cMod 
\]
where $S$ is a set of $N$ elements.  If $N\leq 1$, then there is no ambiguity in the above definition. If $N=2$, then the ladybug matching \cite[Section 5.5]{LS} is used to determine the pairing of the two pairs of boundary points in $\de \ModSp_J(x,y)$. To define the thick embedding, recall that 
\[
\AmbS = \cMod 
\]
in this case, so we let
\[
\ModSp_J(y,x) \xthickemb{\coprod \id} \AmbS 
\]
be a disjoint union of the identity map. In particular, the above thick embedding is a trivial cover requiring no framing data, and this is exactly how the cubical flow category $\ffc_{Kh}$ is constructed in \cite{LS} (the framing data in \cite{LS} is entirely inherited by the neat embedding of the cube moduli spaces $\cMod\hookrightarrow \cEuc$; see Section \ref{sec:ffc review} and Definition 5.5 there).

\subsection{Square faces involving $E$ or $F$}
\label{sec:m=2 J=E/F involved}
Now we suppose that the subcube $u\leq_2 v$ involves the $J=E$ map; the case involving $J=F$ is nearly identical (see Lemma \ref{lem:E and F duality}).  Without loss of generality, let the edges $u \to w_1$ and $w_2\to v$ be saddle maps $S$, while $u\to w_2$ and $w_1\to v$ are the $E$ map. In particular, we have equality of resolutions $D_u = D_{w_2}$, $D_{w_1} = D_v$, and the saddles $D_u \to D_{w_1}$ and $D_{w_2} \to D_v$ are identical. 
\begin{equation*}
\label{eq:square with J}
\begin{tikzcd}
D_u \ar[d, "E"'] \ar[r, "S"] & D_{w_1} \ar[d, "E"] \\
D_{w_2} \ar[r, "S"] & D_v
\end{tikzcd}
\end{equation*}

\begin{figure}
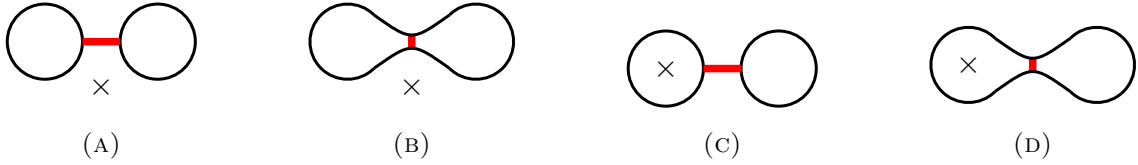

\begin{subfigure}[b]{.25\textwidth}
\begin{center}
\includestandalone{arc_diagram_multiplying_trivials}
\end{center}
\caption{}\label{fig:multiplying trivial circles}
\end{subfigure}
\begin{subfigure}[b]{.25\textwidth}
\begin{center}
\includestandalone{arc_diagram_splitting_trivial_into_trivials}
\end{center}
\caption{}\label{fig:splitting trivial into trivials}
\end{subfigure}
\begin{subfigure}[b]{.25\textwidth}
\begin{center}
\includestandalone{arc_diagram_multiplying_trivial_and_essential}
\end{center}
\caption{}\label{fig:multiplying trivial and essential}
\end{subfigure}
\begin{subfigure}[b]{.25\textwidth}
\begin{center}
\includestandalone{arc_diagram_splitting_trivial_from_essential}
\end{center}
\caption{}\label{fig:splitting trivial from essential}
\end{subfigure}
\caption{Connected arc diagrams where the number of essential circles is unchanged}\label{fig:saddles not changing number of ess circles}
\end{figure}

We focus first on connected arc diagrams $D_u$, and we begin by considering those for which the saddle $S$ does not change the number of essential circles.  There are four of these, as shown in Figure \ref{fig:saddles not changing number of ess circles}.

For Figures \ref{fig:multiplying trivial circles} and \ref{fig:splitting trivial into trivials}, the $E$ map is the zero map, and all of the moduli spaces are empty.

\begin{figure}
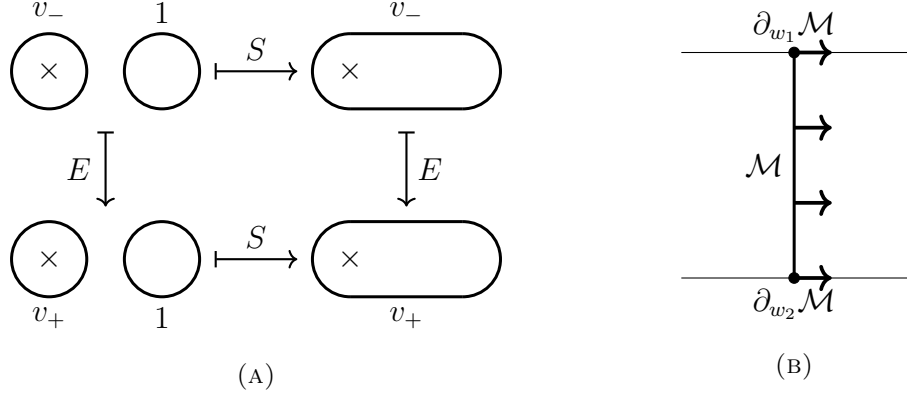

\hspace*{\fill}
\begin{subfigure}{.4\textwidth}
    \begin{center}
    \includestandalone{square_face_bijection_ex}
    \end{center}
    \caption{}\label{fig:mult ess and triv gens}
\end{subfigure}
\hfill
\begin{subfigure}{.2\textwidth}
    \begin{center}
    \includestandalone{framed_embedding_bijection_ex}
    \end{center}
    \caption{}\label{fig:mult ess and triv modsp}
\end{subfigure}
\hspace*{\fill}
\caption{The relevant generators in the subcube of resolutions for the case of Figure \ref{fig:multiplying trivial and essential}, together with the corresponding moduli space embedded and framed into $D^1\times (0,1)$.}
\end{figure}

In the case of Figure \ref{fig:multiplying trivial and essential}, we see that all of the moduli spaces are empty unless $x=1 \o v_-$ and $y=v_+$, in which case we have single point composition moduli spaces
\[
\pathPts{w_1}{\MspExy} = \{(1\o v_- \rightarrow v_- \xrightarrow{1} v_+)\},
\]
\[
\pathPts{w_2}{\MspExy} = \{(1 \o v_- \xrightarrow{1} 1 \o v_+  \rightarrow v_+)\}
\]
as illustrated in Figure \ref{fig:mult ess and triv gens}.

Figure \ref{fig:splitting trivial from essential} is dual to Figure \ref{fig:multiplying trivial and essential}, so the analysis is similar.  We have empty moduli spaces unless $x=v_-$ and $y=X\o v_+$, in which case we have
\begin{align*}
\pathPts{w_1}{\MspExy} & = \{(v_-\rightarrow X \o v_- \xrightarrow{1} X \o v_+)\},\\
\pathPts{w_2}{\MspExy} & = \{(v_- \xrightarrow{1} v_+ \rightarrow X \o v_+)\}
\end{align*}

In all of these cases where $S$ does not change the number of essential circles,  we see that there is a bijection of one-point composition moduli spaces that commutes with the thin embeddings
\begin{equation}
\begin{aligned}
\label{eq:1dim bij for trivial S}
\begin{tikzpicture}[x=1in,y=-.5in]

\node(A) at (0,0) {$\pathPts{w_1}{\MspExy}$};

\node(B) at (0,2) {$\pathPts{w_2}{\MspExy}$};

\node(C) at (2,0) {$\AmbSp{v}{w_1}\times\AmbSp{w_1}{u}$};

\node(D) at (2,2) {$\AmbSp{v}{w_2}\times\AmbSp{w_2}{u}$};

\node(E) at (2,1) {$(0,1)$};

\draw[->] (A)-- node[left]{$\cong$} (B);
\draw[right hook->] (A) -- (C);
\draw[right hook->] (B) -- (D);
\path (C)-- node[rotate=-90]{$\cong$} (E);
\path (E)-- node[rotate=-90]{$\cong$} (D);

\end{tikzpicture}.
\end{aligned}
\end{equation}

This commuting bijection allows us to define $\MspExy$ as a single closed interval $D^1$ which is thinly embedded as an identity
\[\MspExy:=D^1\xhookrightarrow{\id \times\{p_1\}} D^1\times (0,1)=\AmbS\]
with the constant framing. See Figure \ref{fig:mult ess and triv modsp} for a depiction of the case where $S$ is the saddle in Figure \ref{fig:multiplying trivial and essential}.

\begin{figure}
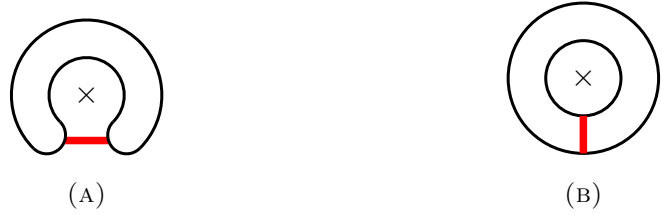

\begin{subfigure}[t]{.4\textwidth}
\begin{center}
\includestandalone{arc_diagram_splitting_trivial_into_essentials}
\end{center}
\caption{}\label{fig:splitting trivial into essentials}
\end{subfigure}
\begin{subfigure}[t]{.4\textwidth}
\begin{center}
\includestandalone{arc_diagram_multiplying_essentials}
\end{center}
\caption{}\label{fig:multiplying essentials}
\end{subfigure}
\caption{Connected arc diagrams where the number of essential circles is changed.}\label{fig:saddles changing number of ess circles}
\end{figure}

Next, we consider the case where $S$ changes the number of essential circles.  For this there are two possible connected arc diagrams, illustrated in Figure \ref{fig:saddles changing number of ess circles}.

In Figure \ref{fig:splitting trivial into essentials}, we see that all moduli spaces are empty unless $x=1$ and $y=v_+\otimes v_+$, in which case we have $\pathPts{w_1}{\MspExy}$ containing two points while $\pathPts{w_2}{\MspExy}$ is empty (see Figure \ref{fig:comult triv gens}).

\[\pathPts{w_1}{\MspExy} = \{(1\rightarrow v_-\otimes v_+ \xrightarrow{1} v_+\otimes v_+), (1\rightarrow v_+\otimes v_- \xrightarrow{2} v_+\otimes v_+)\},\]
\[\pathPts{w_2}{\MspExy} = \varnothing.\]

\begin{figure}
\hspace*{\fill}
\begin{subfigure}{.4\textwidth}
    \begin{center}
    \includestandalone[scale=.8]{square_face_turnback_ex}
    \end{center}
    \caption{}\label{fig:comult triv gens}
\end{subfigure}
\hfill
\begin{subfigure}{.2\textwidth}
    \begin{center}
    \includestandalone{framed_embedding_turnback_ex}
    \end{center}
    \caption{}\label{fig:comult triv modsp}
\end{subfigure}
\hspace*{\fill}
\caption{The relevant generators in the subcube of resolutions for the case of Figure \ref{fig:splitting trivial into essentials}, together with the corresponding moduli space embedded and framed into $\AmbS \cong D^1\times (0,1)$.}
\end{figure}

Following the embedding instructions from Section \ref{sec:0-dim moduli spaces}, we see that our two points in $\pathPts{w_1}{\MspExy}$ are embedded in $\AmbSp{v}{w_1}\times\AmbSp{w_1}{u}\cong (0,1)$ with framing shown below.

\begin{equation*}
\begin{aligned}

\begin{tikzpicture}
\draw (0,0) -- (3,0);

\filldraw[black] (.5,0) circle (2pt) node[draw=none] {};
\filldraw[black] (2.5,0) circle (2pt) node[draw=none] {};

\draw[->,line width=.5mm] (.5,0) -- (1,0);
\draw[->,line width=.5mm] (2.5,0) -- (2,0);

\node[draw=none, anchor = south] at (.5,0) {$p_1$};
\node[draw=none, anchor = south] at (2.5,0) {$p_2$};

\end{tikzpicture}

\end{aligned}
\end{equation*}

We then define $\MspExy$ in this setting to be a single interval which is thinly embedded into $\AmbS\cong D^1\times (0,1)$ as a framed turnback as indicated in Figure \ref{fig:comult triv modsp}.

The last connected case is dual to the previous one, where the moduli spaces are empty unless $x=v_-\otimes v_-$ and $y=X$, and the roles of $\pathPts{w_1}{\MspExy}$ and $\pathPts{w_2}{\MspExy}$ are reversed.  The details here are left to the reader.

Finally, we consider the situation when the diagram $D_u$ is disconnected.  We let $\ConnComp\subset D_u$ denote the connected component of the saddle arc $S$ within $D_u$, while $\ConnComp'$ denotes the complement of $\ConnComp$.  As in Section \ref{sec:Arc diagrams}, these components determine connected components in the resolutions at every vertex, which we also denote by $\ConnComp$ and $\ConnComp'$. Recall from Section \ref{sec:Arc diagrams} the notion of restriction of generators. This allows us to write generators as $x=x|_\ConnComp \otimes x|_{\ConnComp'}$ and $y=y|_\ConnComp\otimes y|_{\ConnComp'}$, and we use this breakdown to aid in our analysis.  

\begin{lemma}
\label{lem:m=2 conn comp and annular deg}
Suppose some $\pathPts{w_i}{\MspExy}$ is nonempty.  Then precisely one of the two following statements are true.
\begin{enumerate}
    \item \label{it:m=2 E acts on C}
    $\adeg(y|_\ConnComp)=\adeg(x|_\ConnComp)+2$ and $\adeg(y|_{\ConnComp'}) = \adeg(x|_{\ConnComp'})$.  This corresponds to $E$ acting on an essential circle in $\ConnComp$.
    \item \label{it:m=2 E acts on C'}
    $\adeg(y|_\ConnComp)=\adeg(x|_\ConnComp)$ and $\adeg(y|_{\ConnComp'}) = \adeg(x|_{\ConnComp'})+2$.  This corresponds to $E$ acting on an essential circle in $\ConnComp'$.
\end{enumerate}
When $J=F$, the value $2$ above is replaced by $-2$.
\end{lemma}
\begin{proof}
Saddle maps preserve annular degree, but the $E$ map increases annular degree by two.  Similarly, the $F$ map decreases annular degree by two.
\end{proof}

If $\pathPts{w_1}{\MspExy}=\pathPts{w_2}{\MspExy}=\varnothing$, we define $\MspExy:=\varnothing$ as well.  Otherwise, we use Lemma \ref{lem:m=2 conn comp and annular deg} to consider the two cases.

\begin{figure}
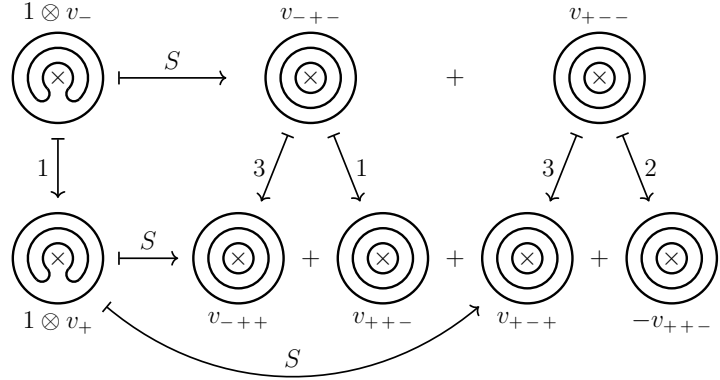
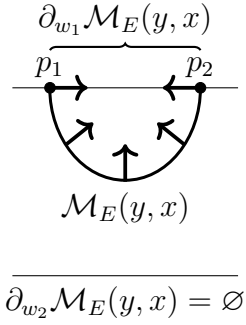
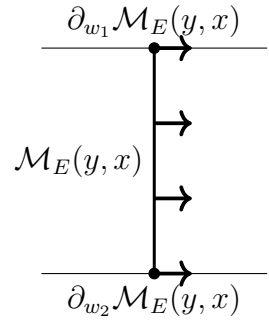

\begin{subfigure}{\textwidth}
    \begin{center}
    \includestandalone[scale=.8]{square_face_turnback_ex_disconnectedv2}
    \end{center}
    \caption{The generators and differentials within the subcube $\subcube{2}{v}{u}$.  The numbers on maps indicate which essential circle is being acted on by $E$.  Certain tensor products are written with subscripts to avoid clutter.}
\end{subfigure}\\
\hspace*{\fill}
\begin{subfigure}{.4\textwidth}
    \begin{center}
    \includestandalone{framed_embedding_turnback_ex_disconnected}
    \end{center}
    \caption{The moduli space for $x=1\o v_-$ and $y=v_+ \o v_+ \o v_-$, embedded and framed in $\AmbS$.}
    \label{fig:framed embedding turnback ex disconnected}
\end{subfigure}
\hfill
\begin{subfigure}{.4\textwidth}
    \begin{center}
    \includestandalone{framed_embedding_bijection_ex_disconnected_1}
    \end{center}
    \caption{The moduli space for $x=1\o v_-$ and either $y=v_-\o v_+\o v_+$ or $y=v_+\o v_-\o v_+$, embedded and framed in $\AmbS$.}
    \label{fig:framed embedding bijection ex disconnected}
\end{subfigure}
\hspace*{\fill}
\caption{Some examples of $\MspExy\hookrightarrow\AmbS\cong D^1 \times(0,1)$ for a subcube $\subcube{2}{v}{u}$ involving $E$ for a disconnected arc diagram.  For the fixed $x=1\o v_-$, different values of $y$ correspond to different cases of Lemma \ref{lem:m=2 conn comp and annular deg}. }
\end{figure}

In the case \eqref{it:m=2 E acts on C} where $E$ was acting upon the component $\ConnComp$ containing $S$, we must have the labels of $y|_{\ConnComp'}$ matching the labels of $x|_{\ConnComp'}$, in which case these components are irrelevant for the analysis of the moduli spaces and we may define
\[\MspExy:=\MspE{y|_\ConnComp}{x|_\ConnComp}\]
as defined and embedded earlier.  This case is illustrated for a specfic example in Figure \ref{fig:framed embedding turnback ex disconnected}.

Meanwhile, in the case \eqref{it:m=2 E acts on C'} where $E$ was acting on some essential circle in $\ConnComp'$ disconnected from the saddle $S$, we see that the generator at $w_1$ must be $y|_\CC \o x|_{\CC'}$ while the generator at $w_2$ must be the `opposite' $x|_\CC \o y|_{\CC'}$.  This allows us to show that the maps $S$ and $E$ (and their corresponding edge moduli spaces) commute in the following sense:

\begin{equation}
\label{eq:edges commute}
\begin{aligned}
\pathPts{w_1}{\MspExy} &=
\MspE{y|_\CC \o y|_{\CC'}} {y|_\CC \o x|_{\CC'}} \times \MspE{y|_\CC \o x|_{\CC'}} {x|_\CC \o x|_{\CC'}}\\
& \cong \MspE{y|_{\CC'}}{x|_{\CC'}} \times \MspE{y|_\CC}{x|_\CC}\\
& \cong \MspE{y|_\CC}{x|_\CC} \times \MspE{y|_{\CC'}}{x|_{\CC'}}\\
& \cong \MspE{y|_\CC \o y|_{\CC'}} {x|_\CC \o y|_{\CC'}} \times \MspE{x|_\CC \o y|_{\CC'}} {x|_\CC \o x|_{\CC'}}\\
& = \pathPts{w_2}{\MspExy}.
\end{aligned}
\end{equation}

The framed embedding data for such moduli spaces is determined by which circle is being acted on by $E$, and this is maintained throughout all of the equivalences above.  Therefore we again have a (one point) bijection $\pathPts{w_1}{\MspExy}\cong\pathPts{w_2}{\MspExy}$ that commutes with the thin embeddings as in Equation \eqref{eq:1dim bij for trivial S}, allowing us to once again define $\MspExy$ as
\[\MspExy:=D^1\xhookrightarrow{I\times\{p_1\}} D^1\times (0,1)=\AmbS.\]
This case is illustrated for a specific example in Figure \ref{fig:framed embedding bijection ex disconnected}.

\subsection{Square faces involving $H$}
\label{sec:m=2 J=H involved}

Finally, we turn to the case where $J=H$. As before, we distinguish between two cases. 

Suppose first that the saddle $S$ does not change the number of essential circles in $D_u$. Let $s$ denote the number of essential circles in $D_u$ (equivalently, in any of $D_u, D_{w_1}, D_{w_2}$, or $D_v$). If $s=0$, then $\de_{w_1} \MspHxy = \de_{w_2} \MspHxy = \varnothing$, and we define $\MspHxy = \varnothing$. Otherwise, there are unique generators in  $z \in \KG{w_1}$ and $z'\in \KG{w_2}$ such that 
\begin{align*}
\de_{w_1}\MspHxy &= \MspH{y}{z} \times \MspH{z}{x} \\ \de_{w_2} \MspHxy &= \MspH{y}{z'} \times \MspH{z'}{x}.
\end{align*}
In fact, $z=y$ and $z'=x$ as generators of $D_u=D_{w_2}$ and $D_v=D_{w_1}$, respectively. The edge moduli spaces $\de_{w_1}\MspHxy$ and $\de_{w_2}\MspHxy$ each have $s$ elements, thinly embedded in $(0,1)$ in order according to Equation \eqref{eq:m=1 J-embedding}.  From the formulas \eqref{eq:formula1} and \eqref{eq:formula3}, we see that the essential circles in $D_u = D_{w_2}$ and $D_v=D_{w_1}$ have the same labels in $x=z'$ and $y=z$, respectively, so the $i$-th point in $\de_{w_1}\MspHxy$ has the same framing as the $i$-th point in $\de_{w_2}\MspHxy$.  Thus we have a bijection of $s$-point composition moduli spaces that commutes with the thin embeddings precisely as in Equation \eqref{eq:1dim bij for trivial S}.  We therefore define $\MspHxy$ to be a disjoint union of $s$ intervals, embedded in $\AmbS$
as a framed identity cobordism from $\de_{w_1}\MspHxy$ to $\de_{w_2}\MspHxy$ connecting the $i$-th point of $\de_{w_1}\MspHxy $ to the $i$-th point of $\de_{w_2} \MspHxy$.

\begin{figure}
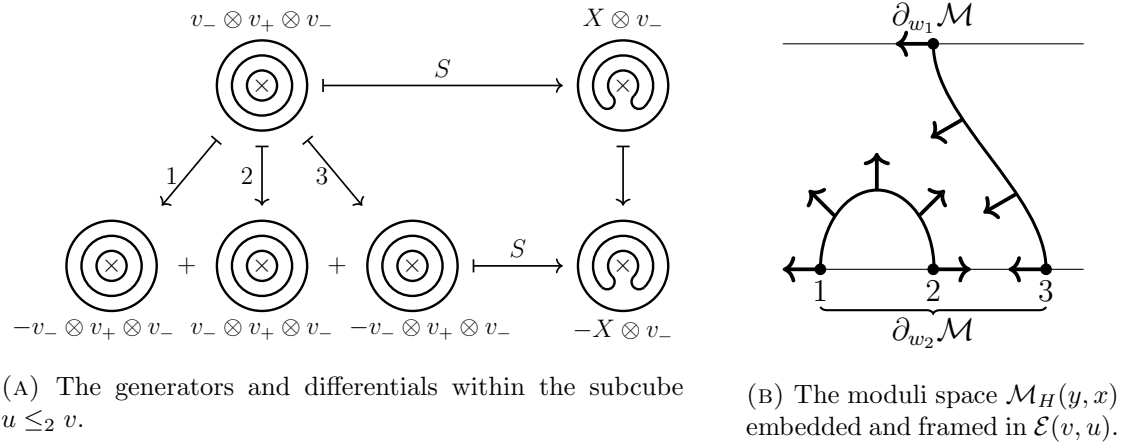

\hspace*{\fill}
\begin{subfigure}{.55\textwidth}
    \begin{center}
    \includestandalone[scale=.8]{square_face_turnback_ex_disconnected_H}
    \end{center}
    \caption{The generators and differentials within the subcube $\subcube{2}{v}{u}$.}
\end{subfigure}
\hfill
\begin{subfigure}{.3\textwidth}
    \begin{center}
    \includestandalone{framed_embedding_turnback_ex_disconnected_H}
    \end{center}
    \caption{The moduli space $\MspHxy$ embedded and framed in $\AmbS$.}
\end{subfigure}
\hspace*{\fill}
\caption{An example of $\MspHxy\hookrightarrow\AmbS\cong D^1 \times(0,1)$ for a subcube $\subcube{2}{v}{u}$ involving an $H$ edge and a saddle $S$ which decreases the number of essential circles by two.  The two merged circles must be consecutive, allowing a framed turnback to `cancel' the corresponding points in $\de_{w_2}\MspHxy$.}\label{fig:m=2 framed embedding for H ex}
\end{figure}

Suppose now that $S$ changes the number of circles. Suppose that $D_v=D_{w_1}$ has $s$ circles and $D_u = D_{w_2}$ has $s+2$, so that $S$ is the Type II saddle from Figure \ref{fig:type2}. If both $\de_{w_1}\MspHxy = \de_{w_2}\MspHxy = \varnothing$, then define $\MspHxy = \varnothing$. Otherwise, there are generators $z\in \KG{w_1}$, $z' \in \KG{w_2}$ such that 
\begin{align*}
\de_{w_1}\MspHxy &= \MspH{y}{z} \times \MspH{z}{x} \\ \de_{w_2} \MspHxy &= \MspH{y}{z'} \times \MspH{z'}{x}.
\end{align*}
The edge moduli spaces $\MspH{z}{x}$ and $\MspH{y}{z'}$ each have one element, while $\MspH{y}{z}$ and $\MspH{z'}{x}$ have $s$ and $s+2$ elements, respectively, corresponding to the number of essential circles in $D_{w_1}$ and $D_{w_2}$. Note that the two essential circles merged by $S$ are consecutive in $D_{w_2}$, and their labels in $z'$ are different. Define $\MspHxy$ to consist of $s+1$ intervals, with $s$ of them thinly embedded in $\AmbS$ matching the $s$ points in $\MspH{y}{z}$ with the $s$ points in $\MspH{z'}{x}$ that correspond to essential circles not participating in the saddle. The remaining interval in $\MspHxy$ is thinly embedded as a turnback with boundary the two remaining unmatched points in $\de_{w_2} \MspHxy$ corresponding to the essential circles created by $S$. See Figure \ref{fig:m=2 framed embedding for H ex} for an example. Finally, the case when $S$ increases the number of essential circles by two is dual to the above discussion. 

\section{Path moduli spaces}  
\label{sec:Path moduli spaces}

This section contains a detailed combinatorial analysis that will be crucial for completing our inductive construction in Section \ref{sec:m-1 to m}. For an $m$-dimensional subcube $u\leq_m v$ and generators $x\in \KG{u}, y\in \KG{v}$, Item \eqref{item:ffcj boundary definitions} of Theorem \ref{thm:ffcJ}  states that the codimension $m-1$ boundary of $\MspJxy$ is built out of products of the moduli spaces associated to a sequence of edges from $u$ to $v$. Such a sequence of edges is a \emph{path} from $u$ to $v$, and each edge moduli space, along with its thick embedding, was defined in Section \ref{sec:0-dim moduli spaces}. We analyze these \emph{path  moduli spaces} and their thick embeddings, obtained as the product of the thick embeddings for each edge moduli space. While each edge moduli space is thinly embedded, taking the disjoint union over intermediate generators as in Theorem \ref{thm:ffcJ} \eqref{item:ffcj boundary definitions} may produce a thick embedding.

We focus on connected arc diagrams in this section. The key results are Proposition \ref{prop:E/F conn is thin} and Proposition \ref{prop:conn H}. In Proposition \ref{prop:E/F conn is thin} we show that for $J=E, F$ all path moduli spaces involving $J$ are in fact thinly embedded. For $J=H$, path moduli spaces are thinly embedded except for one type of arc diagram, for which an explicit description is provided, see Proposition \ref{prop:conn H}. The outline is as follows. In Section \ref{sec:basic defs and lemmas} we introduce path moduli spaces, their thick embeddings, and present a helpful lemma which bounds the number of points in some path moduli spaces. Section \ref{sec:equiv path moduli spaces} establishes two key lemmas. First, Lemma \ref{lem:biject to all ess} allows us to assume that a path has the special property that the arc diagram just before the $J$ edge (the so-called \emph{$J$ resolution}) consists entirely of essential circles. Second, in Lemma \ref{lem:forbidden configs} we give restrictions on the type of arc configurations that may appear in such a diagram.  

To establish thinness we focus on the $J$-resolution arc diagram. While the two cases $J=E,F$ and $J=H$ require different arguments, the overall strategy is the following. Connectivity and the formulas \eqref{eq:formula2} and \eqref{eq:formula4} imply that picking an essential circle for $J$ to act on and a label $v_\pm$ on that circle determines labels on all the other circles. Lemma \ref{lem:genus considerations} ensures thinness ``before and after'' the $J$ edge. With further analysis, depending on $J=E, F$ or $J=H$, we arrive at our desired results. Note that connectivity of the arc diagram is critical: disconnected arc diagrams may yield non-thinly embedded path moduli spaces. See also the discussion at the beginning of Section \ref{sec:m-1 to m J=E disconn}.

\subsection{Basic definitions and lemmas}
\label{sec:basic defs and lemmas}

\begin{definition}
\label{def:path}
A \emph{path} $P$ through a subcube $\scube$ is a sequence of edges, denoted as
\[P = (u\rightarrow w_1 \rightarrow \cdots \rightarrow w_{m-1} \rightarrow v).\]
Fixing generators $x\in\KG{u},y\in\KG{v}$, any such path determines a \emph{path moduli space}
\[\pathP := \coprod \MspJ{y}{z_{m-1}} \times \MspJ{z_{m-1}}{z_{m-2}} \times \cdots \times \MspJ{z_1}{x},\]
where the disjoint union is taken over all sequences $(z_1,\ldots,z_{m-1})$ with $z_i\in\KG{w_i}$ for each $1\leq i \leq m-1$.  As in Section \ref{sec:m=2}, points in $\pathP$ will be denoted as
\[(x\rightarrow z_1 \rightarrow\cdots \rightarrow z_{m-1} \rightarrow y)\in\pathP,\]
with a label on an arrow corresponding to the $J$-edge indicating which essential circle is being acted on by $J$ if necessary (see Equation \eqref{path eq} for an example of such a label).
\end{definition}

Recall from Definition \ref{def:flow cat} that we use the notation $\de_{[i]}\ModSp$ to indicate the codimension-$i$ portion of the boundary of a moduli space.  In order to satisfy the demands of Theorem \ref{thm:ffcJ}, we define
\begin{equation}\label{disjoint union paths}
\de_{[m-1]}\MspJxy \cong \coprod_{\text{paths $P$}} \pathP
\end{equation}
where the disjoint union is taken over paths $P$ from $u$ to $v$; to make the notation less cumbersome we omit reference to $u, v$.

Furthermore, every $\pathP$ inherits a thick embedding
\begin{equation}
\label{eq:product cube ambient space}
\pathP \xthickemb{\varphi_P} \mathcal{E}(P) := \AmbSp{v}{w_{m-1}} \times \cdots \times \AmbSp{w_1}{u}
\end{equation}
via the disjoint union of products of the thin embeddings
\[
\MspJ{z_{m-i}}{z_{m-i-1}} \hookrightarrow \AmbSp{w_{m-i}}{w_{m-i-1}}
\]
defined in Section \ref{sec:0-dim moduli spaces}.  If the subcube $\scube$ does not involve $J$, then $\pAmb$ is a point and the thick embedding is constant.  Otherwise, $J$ is involved and there is a natural identification $\pAmb \cong (0,1)$.   Then the thick embedding \eqref{eq:product cube ambient space} sends a point
\begin{equation} \label{path eq}
(x \to z_1 \to \cdots \to z_i \rar{j} z_{i+1} \to \cdots \to z_{m-1} \to y) \in \pathPts{P}{\MspJxy},
\end{equation}
to the fixed $j$-th point in $(0,1)$ with framing determined by the sign $(-1)^{j+1}$ if $J=E, F$ and otherwise given by the label of the $j$-th circle in $z_i$ if $J=H$, as described in Section \ref{sec:0-dim moduli spaces}. 

We may also consider \emph{sub-paths} of a path $P$.   If $w_k$ is any vertex along a path
\[P=(u\rightarrow w_1\rightarrow \cdots \rightarrow w_k \rightarrow \cdots \rightarrow v),\]
we may define $P_{\leq w_k} :=(u\rightarrow\cdots \rightarrow w_k)$ and $P_{\geq w_k} :=(w_k\rightarrow \cdots \rightarrow v)$.
This then allows us to split the path moduli space as a product
\[\pathP = \coprod_{z\in \KG{w_k}} \pathPts{P_{\geq w_k}}{\MspJ{y}{z}} \times \pathPts{P_{\leq w_k}}{\MspJ{z}{x}}\]
embedded as a product with the product framing.

\begin{definition}
For a fixed path $P$ through a sub-cube $\scube$ that involves $J$, we single out the vertices along $P$ directly before and after the $J$ map has been applied:

\begin{itemize}[label={$\bullet$}]
\item $e_0=e_0(P)$ denotes the source vertex of the $J$ edge in $P$, and
\item $e_1=e_1(P)$ denotes the target vertex of the $J$ edge in $P$.
\end{itemize}

The arc diagram $D_e:=D_{e_0}=D_{e_1}$ being acted on by $J$ will be referred to as the \emph{$J$-resolution of $P$}.\footnote{To relate this notation to equation (\ref{path eq}), $e_0=w_i$, $e_1=w_{i+1}$.}
\end{definition}

Note that the $J$-resolution of a path $P$ depends on the path, although we do not include $P$ in the notation to avoid clutter.

Now if $u\leq_m v$ does not involve the $J$ map, then any path $P$ from $u$ to $v$ represents a path in the cube of resolutions of $D$.  For each such path there is an associated cobordism $S \: D_u \to D_v$.

\begin{lemma}\label{lem:genus considerations}
Let $P$ be a path from $u$ to $v$ which does not involve the $J$ map, and let $S \: D_u \to D_v$ denote its corresponding surface cobordism. Let $g$ denote the maximal genus among all components of $S$. 

\begin{enumerate}
\item If $g=0$, then for every $x\in \KG{u}$ and $y\in \KG{v}$, the moduli space $\pathPts{P}{\MspJxy}$ contains at most one element. 
\item If $g \geq 1$ and one of $D_u$ or $D_v$ consists only of essential circles, then $\Fa(S) = 0$. In particular in this case, for every $x\in \KG{u}$ and $y\in \KG{v}$, the moduli space $\pathPts{P}{\MspJxy}$ is empty. 
\end{enumerate}
\end{lemma}
\begin{proof}
For (1), first view $S$ as a cobordism in $\R^2\times I$ using the inclusion $\A \times I \hookrightarrow \R^2\times I$. Since $\Fa$ is the $\adeg$-preserving part of $\FKh$, the cardinality of $\pathPts{P}\MspJxy$ is bounded above by the coefficient of $y$ in $\FKh(S)(x)$. Represent the generators $x$ and $y$ using (dotted) cup cobordisms, denoted $\Sigma_x$ and $\Sigma_y$ respectively. Let $\Sigma_y^*$ denote the dual cap cobordism to $\Sigma_y$, so that the closed surface 
\[
\Sigma_y^* \circ \Sigma_y
\]
evaluates to $1$ under $\FKh$. The coefficient of $y$ in $\FKh(S)(x)$ is then the value of $\FKh$ on the closed surface
\[
S':= \Sigma_y^* \circ S \circ \Sigma_x.
\]
Since $S$ has genus zero, every component of $S'$ is a (possibly dotted) sphere. Then $\FKh(S') \leq 1$, which completes the proof of statement (1). 

We now address (2). Let $\overline S\subset S$ be a component of positive genus. It suffices to show $\Fa(\overline S)=0$. By neck-cutting, Figure \ref{fig:neck-cutting}, we can write $\overline S$ as a sum 
\[
S_1 + \cdots + S_k,
\]
where each component of each $S_i$ is either an annulus with essential boundary or a disk, possibly carrying dots. We may assume that each component carries at most one dot. By assumption, each $S_i$ contains an annulus. We will show that each $S_i$ has a dotted annulus, which implies $\Fa(\overline S) = 0$ since $X$ acts by zero on essential circles, see formula \eqref{eq:formula1}.

Fix $1\leq i \leq k$. Recall that the Bar-Natan relations are homogeneous, so $\deg(\overline S) = \deg(S_i)$. Let $n$ denote the number of (possibly dotted) disks in $S_i$. If none of the annuli in $S_i$ are dotted, then 
\[
\deg(S_i) \geq -n. 
\]
On other other hand, considering the connected surface $\overline{S}$,
\[
\deg(\overline S) = \chi(\overline S) < -n,
\]
since $g(\overline S)\geq 1$ and the number of boundary components of $\overline S$ is strictly greater than $n$. Therefore $S_i$ contains a dotted annulus, which completes the proof. 
\end{proof}

\subsection{Equivalent path moduli spaces and a non-emptiness condition}
\label{sec:equiv path moduli spaces}
In this section we analyze the manner in which two paths $P,P'$ can give rise to equivalent path moduli spaces $\pathPts{P}{\MspJxy}\cong\pathPts{P'}{\MspJxy}$.  This analysis will allow us to choose certain preferred paths through our subcubes.  It will also give rise to an important non-emptiness condition on path moduli spaces.

\begin{definition}
Let $P$ and $P'$ be two paths in a subcube $u\leq_m v$. Write
\begin{align*}
P &= \left( u = w_0 \to w_1 \to \cdots \to w_{m} = v \right) \\
P' &= \left( u = w_0' \to w_1' \to \cdots \to w_{m}' = v \right).
\end{align*}
We say that $P$ and $P'$ are \emph{related by a square face} if for some $1 \leq j \leq m-1$, we have $w_j\neq w_j'$ while $w_i = w_i'$ for all indices $i\neq j$. We say the square face involves $J$ if one of the edges in the $2$-dimensional subcube $w_{j-1} \leq_2 w_{j+1}$ is the $J$ map. In this case, two of the opposite edges are $J$, and the remaining two opposite edges are a saddle. 
\end{definition}

\begin{lemma}\label{lem:square face bijections}
Let $P$ and $P'$ be two paths in a subcube $u\leq_m v$ that are related by a square face $w_{j-1} \leq_2 w_{j+1}$ as above, and let $x\in \KG{u}$, $y\in \KG{v}$. If either
\begin{enumerate}
\item the square face does not involve the $J$ map, or
\item the square face involves the $J$ map, and the number of essential circles in the resolutions appearing in $w_{j-1} \leq_2 w_{j+1}$ is the same,
\end{enumerate} 
then there is a bijection 
\[
\pathP \cong \pathPts{P'}{\MspJxy}
\]
such that the diagram \eqref{eq:path bijection} commutes.
\begin{equation}
\label{eq:path bijection}
\begin{tikzpicture}[x=1.2in,y=-.8in,baseline=(current  bounding  box.center)]
\node(A) at (0,0){$\pathP$};
\node(B) at (0,1){$\pathPts{P'}{\MspJxy}$};
\node(C) at (1,0){$\mathcal{E}(P)$};
\node(D) at (1,1){$\mathcal{E}(P')$};

\thickembtikz{A}{C}
\thickembtikz{B}{D}
\draw[->] (A) -- node[left]{$\cong$} (B);
\draw[->] (C) -- node[left]{$\cong$} (D);

\end{tikzpicture}
\end{equation}

\end{lemma}

\begin{proof}
The desired bijection is assembled from the bijections in Section \ref{sec:m=2} on the factors of $\pathP$ and  $\pathPts{P'}{\MspJxy}$ involving the square face and the identity on all other factors. Commutativity of the diagram is immediate for item (1), since $\mathcal{E}(P) = D^0 = \mathcal{E}(P')$ in this case. For item (2) it follows from commutativity of the diagram \eqref{eq:1dim bij for trivial S} and the discussion below Lemma \ref{lem:m=2 conn comp and annular deg}.
\end{proof}

In the case that the number of essential circles in the $J$-resolution has changed, we instead see turnbacks as in Figure \ref{fig:comult triv modsp}, but this fact does not concern us at the moment.

\begin{corollary}
\label{cor:fixed J edge}
Let $P$ and $P'$ be paths from $u$ to $v$ such that the $J$-resolution of $P$ and $P'$ are at the same vertex, $e_0(P) = e_0(P')$. Then there is a bijection 
\[
\pathP \cong \pathPts{P'}{\MspJxy}
\]
such that the diagram \eqref{eq:path bijection} commutes.
\end{corollary}

\begin{proof}
The paths $P$ and $P'$ can be related by a sequence of square faces that do not involve the $J$ map. The claim follows from Lemma \ref{lem:square face bijections}. 
\end{proof}

Now we consider how the various bijections of Lemma \ref{lem:square face bijections} can be used to replace a given path through a subcube with one that is easier to analyze.

\begin{lemma}\label{lem:biject to all ess}
Consider a path $P$ through the sub-cube $\scube$, specifying a $J$-resolution arc diagram $D_e$, and fix generators $x\in\KG{u},y\in\KG{v}$. Then there exists a path $P'$ through $\scube$, with $J$-resolution $D_{e'}$, satisfying the following properties.
\begin{itemize}
\item Each connected component of $D_{e'}$ contains the same number of essential circles as the corresponding component of $D_e$.
\item Any connected component of $D_{e'}$ which contains at least one essential circle also contains no trivial circles.
\item There is a bijection $\pathPts{P}{\MspJxy} \cong \pathPts{P'}{\MspJxy}$ such that the diagram \eqref{eq:path bijection} commutes.
\end{itemize}
\end{lemma}

\begin{proof}
We work one connected component $\Cs$ of $D_e$ at a time.  If $\Cs$ contains only essential circles or only trivial circles, then we are done.  Otherwise, $\Cs$ must contain a trivial circle that is joined by an arc (either past or future) to an essential circle. This arc specifies a square face in the subcube $\scube$, and the swap across this square face yields a path $P''$ in which the corresponding component $\Cs''$ in the $J$-resolution of $P''$ contains one fewer trivial circle and the same number of essential circles as $\Cs$. By Lemma \ref{lem:square face bijections} (2), the moduli spaces for $P$ and $P''$ are in bijection with diagram \eqref{eq:path bijection} commuting. We continue this process until there are no trivial circles left in any of the connected components with essential circles.   
\end{proof}

Lemma \ref{lem:biject to all ess} ensures that it is enough to understand $\pathPts{P}{\MspJxy}$ for paths $P$ having the $J$-resolution $D_e$ consisting of only essential circles, possibly together with a disconnected component consisting of trivial circles.  However, in the following sections, we will focus largely on connected arc diagrams; the disconnected cases will follow from the connected ones in a manner very similar to the arguments of Section \ref{sec:m=2 J=E/F involved}.

We end this section with an important consequence of Corollary \ref{cor:fixed J edge}.   Recall from Section \ref{sec:Arc diagrams} that a subcube $\scube$ gives rise to its own arc diagrams at each vertex including only those arcs that correspond to saddle maps within the subcube.  Note further that, for any path $P$, the arc diagrams $D_{e_0}$ and $D_{e_1}$ directly before and after the $J$ map are the same.  The following lemma provides an important restriction on the types of arc diagrams that can appear for non-empty path moduli spaces.

\begin{lemma}\label{lem:forbidden configs}
Let $P$ be a path through $u\leq_m v$ and let $w$ be a vertex in $P$. Suppose that any of the configurations in Figure \ref{fig:forbidden configs}, as well as their vertical and horizontal reflections, appear for some pair of arcs in the arc diagram $D_w$. Then $\pathPts{P}\MspJxy = \varnothing$ for every $x\in \KG{u}, y\in \KG{v}$. 

\begin{figure}
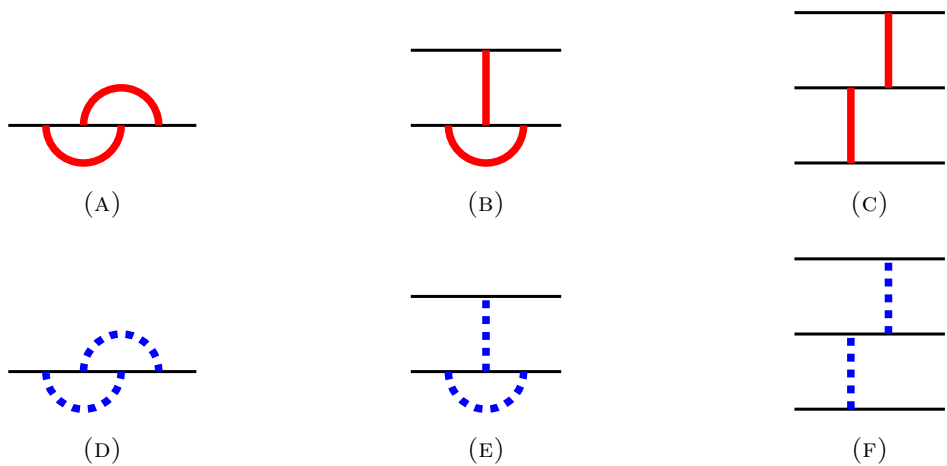

\begin{subfigure}[b]{.3\textwidth}
\begin{center}
\includestandalone{forbidden1}
\end{center}
\caption{}\label{fig:forbidden1}
\end{subfigure} 
\begin{subfigure}[b]{.3\textwidth}
\begin{center}
\includestandalone{forbidden3}
\end{center}
\caption{}\label{fig:forbidden3}
\end{subfigure} 
\begin{subfigure}[b]{.3\textwidth}
\begin{center}
\includestandalone{forbidden5}
\end{center}
\caption{}\label{fig:forbidden5}
\end{subfigure} \\
\vskip1em
\begin{subfigure}[b]{.3\textwidth}
\begin{center}
\includestandalone{forbidden2}
\end{center}
\caption{}\label{fig:forbidden2}
\end{subfigure} 
\begin{subfigure}[b]{.3\textwidth}
\begin{center}
\includestandalone{forbidden4}
\end{center}
\caption{}\label{fig:forbidden4}
\end{subfigure} 
\begin{subfigure}[b]{.3\textwidth}
\begin{center}
\includestandalone{forbidden6}
\end{center}
\caption{}\label{fig:forbidden6}
\end{subfigure} 
\caption{Forbidden configurations; only the two relevant arcs are pictured, but there could be an arbitrary number of additional past and future arcs anywhere within such a diagram.}\label{fig:forbidden configs}
\end{figure}

\end{lemma}

\begin{proof}
Surgery arcs can be performed in any order by Corollary \ref{cor:fixed J edge}. We assume that the depicted pair of arcs correspond to consecutive edges in $P$.  Direct calculations verify that they all yield the zero map.
\end{proof}

\subsection{The $E$ and $F$ maps for connected arc diagrams}
\label{sec:The $E$ and $F$ maps for connected arc diagrams}

In this section we analyze the path moduli spaces $\pathPts{P}\MspJxy$ for connected arc diagrams when $J=E,F$ and show that the thick embeddings \eqref{eq:product cube ambient space} in such cases are in fact thin. (We will see in Section \ref{sec:m-1 to m J=E disconn} how genuinely thick embeddings appear when the arc diagrams are disconnected.) The analysis is stated only for the $E$ map, but the $F$ map is nearly identical, see Lemma \ref{lem:E and F duality}; in particular all results here hold with $F$ and $E$ interchanged, though the proofs may need slight modifications.

Fix an annular link diagram $D$ with $n$-crossings, so that together with the $E$ map we have an $(n+1)$-dimensional cube $\cube$. We also fix a sub-cube $u\leq_m v$. For the remainder of this section, we assume that the sub-cube $u\leq_m v$ contains the $E$ edge and that the arc diagram $D_u$ is connected.  Our first lemma rules out a very specific arrangement of arcs for such a diagram.

\begin{lemma}\label{lem:forbidden E configuration}
Suppose the arc diagram $D_u$ is connected. Let $P$ be a path through $u\leq_m v$ such that the $E$-resolution $D_e$ contains only essential circles. If the configuration shown in \eqref{eq:forbidden E configuration} appears for some pair of arcs in the arc diagram $D_{e_0}$, where the two horizontal intervals are segments of distinct essential circles, then $\pathPts{P}\MspExy = \varnothing$ for every $x\in \KG{u}, y\in \KG{v}$. 

\begin{equation}
\label{eq:forbidden E configuration}
\begin{aligned}
\includestandalone{forbidden_E_configuration}
\end{aligned}
\end{equation}
\end{lemma}

\begin{proof}
We consider two cases. If there are no other essential circles in $D_e$, then by Corollary \ref{cor:fixed J edge} we may assume that the past arc occurs directly before the $E$ map, and the future arc occurs directly after. It is straightforward to verify that the configuration results in empty moduli spaces. If there are other essential circles, then since the arc diagram $D_e$ is connected, at least one of the two depicted circles must be connected to a third circle. Then one of Figure \ref{fig:forbidden5} or Figure \ref{fig:forbidden6} appears in the arc diagram, so we are done by earlier considerations. 
\end{proof}

We now present the main technical result of this section which concerns the thinness of path moduli spaces for connected arc diagrams.

\begin{proposition}\label{prop:E/F conn is thin}
Suppose $P$ is a path through a subcube $\scube$ having connected arc diagram $D_u$.  Then for any fixed generators $x\in\KG{u},y\in\KG{v},$ the thick embedding $\varphi_P$ from \eqref{eq:product cube ambient space} is in fact a thin embedding
\[\pathPts{P}{\MspExy}\xhookrightarrow{\varphi_P} (0,1)\cong \AmbS.\]
\end{proposition}

\begin{proof}
If the path moduli space $\pathPts{P}\MspExy$ is empty, then the statement is trivially true.  Otherwise, by Lemma \ref{lem:biject to all ess}, we can assume that our $E$-resolution consists entirely of essential circles. We fix two points $z,z'\in \pathPts{P} \MspExy$ and seek out to prove that, if $\varphi_P(z)=\varphi_P(z'),$ then in fact we must have had $z=z'$.

We write $z,z'$ as
\begin{align*}
z &= (x \to z_1 \to \cdots \to  z_i \rar{j} z_{i+1} \to \cdots \to z_{m-1} \to y)\\
z' &= (x \to z'_1 \to \cdots \to  z'_i \rar{j'} z'_{i+1} \to \cdots \to z_{m-1}' \to y),
\end{align*}
where $z_i, z_i'\in \KG{e_0}$ are the generators that $E$ acts on.  Then by definition of $\varphi_P$, we must have $j=j'$.  Letting $C$ denote the $j$-th essential circle in $D_e$, this implies that the label on $C$ must be $v_-$ in $z_i$ and $v_+$ in $z_{i+1}$. Based on the assumption that the arc diagram is connected and using Lemmas \ref{lem:forbidden configs}, \ref{lem:forbidden E configuration}, one observes that the curves are connected by past/future arcs in an alternating manner as illustrated in the following figure (also see Figure (\ref{eq:alternating arcs})).
The labels discussed above then determine the labels on the circles adjacent to $C$ as shown below.

\begin{equation*}
\includestandalone{propagation}
\end{equation*}
By propagating this argument we see that the labels on all essential circles in $D_e$ are determined by $j$, and thus for $j=j'$ we must have $z_i=z'_i$ and $z_{i+1}=z'_{i+1}$ as well.

From this we can conclude that both $z$ and $z'$ belong to the same composition moduli space
\[z,z'\in \pathPts{P_{\geq e}}{\MspE{y}{z_{i+1}}} \times \MspE{z_{i+1}}{z_i} \times \pathPts{P_{\leq e}}{\MspE{z_i}{x}}.\]
Then since $D_e$ contains only essential circles, Lemma \ref{lem:genus considerations} ensures that both $\pathPts{P_{\leq e}}{\MspE{z_i}{x}}$ and $\pathPts{P_{\geq e}}{\MspE{y}{z_{i+1}}}$ are at most single point moduli spaces, while $\MspE{z_{i+1}}{z_i}$ is an edge moduli space and so contains at most a single point as indicated in Section \ref{sec:0-dim moduli spaces}.  Thus we have $z=z'$ as desired, and the embedding $\varphi_P$ must in fact be thin.
\end{proof}

Although Proposition \ref{prop:E/F conn is thin} suffices to continue our construction, we can actually prove a stronger statement about our path moduli spaces which makes them easier to understand in concrete examples.  We begin with a special case.

\begin{lemma}\label{lem:E/F conn is one point per circle preliminary}
Suppose $P$ is a path through a subcube $u\leq_m v$ with (connected) $E$-resolution arc diagram $D_e$ consisting entirely of $s$ essential circles, with precisely one arc connecting each consecutive pair.  Fix $x\in \KG{u}, y\in \KG{v}$ such that the moduli space $\pathPts{P}{\MspExy}$ is nonempty.  Then for each $1\leq j \leq s$, there is at least one point of the form
\[ (z^j):= (x\rightarrow\cdots\rightarrow z_i \xrightarrow{j} z_{i+1} \rightarrow \cdots \rightarrow y)\in\pathPts{P}{\MspExy}.\]
\end{lemma}

\begin{proof}
Let $S_{\leq e} \: D_u \to D_{e}$, $S_{\geq e} \: D_e\to D_v$ denote the cobordisms obtained by composing the saddle maps in $P$ before the $E$ map and after the $E$ map, respectively. 

Fixing $1\leq j \leq s$, we first note that the $j$-th essential circle in $D_e$ must be labeled by $\vm$ in $z_i$. Using the argument in Proposition \ref{prop:E/F conn is thin}, there is a unique choice of generators $z_i,z_{i+1}$ satisfying this condition. To construct $(z^j)$, we will show that $\F_\A(S_{\geq e})(z_{i+1})$ is a generator of $\Fa(D_v)$, which we denote $y'$. Then we show there exists a unique $x'\in \KG{u}$ such that $z_i$ appears in $\F_\A(S_{\leq e})(x')$. Finally we will show that $x'=x$ and $y'=y$, completing the proof.

By Lemma \ref{lem:forbidden configs}, Figures \ref{fig:forbidden5} and \ref{fig:forbidden6}, the arc diagram for $D_e$ consists of $s$ essential circles connected by past and future arcs in an alternating manner, 
\begin{equation}\label{eq:alternating arcs}
\begin{aligned}
\includestandalone{alternating_arcs}
\end{aligned}
\end{equation}

Let us first show that $z_{i+1}$ is not killed by applying the future arcs in $D_e$. Note that each circle in $D_e$ is incident to at most one future arc. Pick circles $C_1, C_2$ in $D_e$ that are connected by a future arc $a$. Surgery along $a$ results in a trivial circle $C'$, and by construction of $z_{i+1}$ the labels on $C_1$ and $C_2$ are different. Therefore the map assigned to $a$ is nonzero on $z_{i+1}$, and moreover $C'$ is labeled by $X$ in the image of $z_{i+1}$. Repeating this for all future arcs, we see that $y':=\F_\A(S_{\geq e})(z_{i+1})$ is a generator of $D_v$. Each trivial circle in $D_v$ is labeled $X$ in $y'$. Essential circles in $D_v$ are those with no future arc incident on them in $D_e$, and they have the same label as in $z_{i+1}$. 

Let us now define $x'\in \KG{u}$ such that $z_i$ appears in $\F_\A(S_{\leq e})(x')$. As before, each circle in $D_e$ is incident to at most one past arc. If $a'$ is a past surgery arc joining the essential circles $C_1'$ and $C_2'$ in $D_e$, then undoing $a'$ we see a trivial circle $T$ splitting into the two essential circles $C_1' \cup C_2'$, 

\begin{equation*}
\begin{tikzpicture}[yscale=.75]
   \coordinate (p1) at (0,0);
   \coordinate (p2) at (.5,0);
   \coordinate (p3) at (0,1);
   \coordinate (p4) at (.5,1);
   
   \coordinate (p1') at (2,0);
   \coordinate (p2') at (2.5,0);
   \coordinate (p3') at (2,1);
   \coordinate (p4') at (2.5,1);
   
   \draw[line width=.4mm] (p1) -- (p2);
   \draw[line width=.4mm] (p3) -- (p4);
   \draw[line width=.4mm] (p2) arc[start angle=-90, end angle=90, radius=.5];
   
    \draw[line width=.4mm] (p1') -- (p2');
   \draw[line width=.4mm] (p3') -- (p4');
      \draw[line width=.4mm] (p3') arc[start angle=90, end angle=270, radius=.5];
     
    \coordinate (l1) at (4.5,0);
    \coordinate (r1) at (7,0);
    \coordinate (l2) at (4.5,1);
    \coordinate (r2) at (7,1);
    
    \draw[line width=.4mm] (l1) -- (r1);
   \draw[line width=.4mm] (l2) -- (r2);

\node[draw=none] at (3.5,.5)  {$\longrightarrow$};
   
       \coordinate (l1') at (2.25,2);
    \coordinate (r1') at (4.75,2);
    \coordinate (l2') at (2.25,3);
    \coordinate (r2') at (4.75,3);
    \coordinate (m1) at (3.5,2);
    \coordinate (m2) at (3.5,3);
    
       \draw[line width=.4mm] (l1') -- (r1');
   \draw[line width=.4mm] (l2') -- (r2');
   \draw[line width=1mm, dashed, blue] (m1) -- (m2);
   
   \node[draw=none, anchor = north] at (m1)  {$a'$};
   \node[draw=none, anchor = west] at (r1)  {$C_2'$};   
   \node[draw=none, anchor = west] at (r2)  {$C_1'$};
   
   \node[draw=none, anchor = east] at (p1) {$T$};
\end{tikzpicture}
\end{equation*}
We must label $T$ by the generator $1$. By construction, $C_1'$ and $C_2'$ have opposite labels in $z_i$. The saddle map corresponding to $a'$ sends $1\in \Fa(T)$ to a sum of two generators in $\Fa(C_1' \cup C_2')$, one of which matches the labels on $C_1'\cup C_2'$ in $z_i$. Repeat this for all past arcs in $D_e$ to obtain the desired generator $x'\in \KG{u}$. Note that trivial circles in $D_u$ are labeled by $1$ in $x'$. Essential circles in $D_u$ correspond to essential circles in $D_e$ with no past arc incident on them, and they have the same label in $x'$ and $z_i$.  

It remains to show that $x'=x$ and $y'=y$. By the above analysis, $x'$ and $x$ are labeled $1$ on all trivial circles in $D_u$, and likewise $y'$ and $y$ are labeled $X$ on all trivial circles in $D_v$. We will argue that labels on essential circles in $D_u$ and $D_v$ are determined by the arc diagram $D_e$.

Let $r$ and $t$ denote the number of essential circles in $D_u$ and $D_v$, respectively. Since the arc diagram $D_e$ is connected, Lemma \ref{lem:forbidden configs}, Figures \ref{fig:forbidden5} and \ref{fig:forbidden6} imply that all but the first and last essential circles in $D_e$ are connected to adjacent essential circles by both a past and future arc. It follows that the pair $(r,t)$ is one of $(2,0)$, $(1,1)$, or $(0,2)$, which can be read off from the type of arcs incident on the first and last circle in $D_e$, 

\begin{equation*}
\includestandalone{first_and_last_circles}
\end{equation*}

Finally, note that if $x'' \in \KG{u}$, $y''\in \KG{v}$ such that $\pathPts{P}{\MspE{y''}{x''}} \neq \varnothing$, then $\adeg(y'') = \adeg(x'') + 2$. For each of the three possibilities for $(r,t)$, we see that the labels on essential circles in $D_u$ and $D_v$ are always uniquely determined whenever the moduli space is nonempty. Therefore $x=x'$ and $y=y'$. 
\end{proof}

With Lemma \ref{lem:E/F conn is one point per circle preliminary} in place, we can state and prove the following proposition which characterizes all non-empty path moduli spaces in any subcube for a connected diagram which involves the $E$ map.

\begin{proposition}\label{prop:E/F conn is one point per circle}
Suppose $P$ is a path through a subcube $\scube$ having connected arc diagram $D_u$, such that the $E$-resolution $D_e$ of $P$ contains $s$ essential circles.  Then for any fixed generators $x\in\KG{u},y\in\KG{v}$, if the path moduli space $\pathPts{P}{\MspExy}$ is non-empty, then the thick embedding $\varphi_P$ from \eqref{eq:product cube ambient space} is a thin embedding mapping $\pathPts{P}{\MspExy}$ bijectively onto $\{ \frac{1}{s+1}, \ldots, \frac{s}{s+1} \} \subset (0,1)$.
\end{proposition}

\begin{proof}
Thinness of the embedding $\pathPts{P}{\MspExy}\xhookrightarrow{\varphi_P} (0,1)$ was established in Proposition \ref{prop:E/F conn is thin}.  By Lemma \ref{lem:biject to all ess}, we can assume that our $E$-resolution consists entirely of $s$ essential circles. It thus remains to show that each essential circle in $D_e$ contributes at least one point to $\pathPts{P}{\MspExy}$.

Fix $1\leq j \leq s$. By Corollary \ref{cor:fixed J edge}, the order in which the surgery arcs in the sub-cubes $u\leq e_0$ and $e_1 \leq v$ are performed is irrelevant, and we use this implicitly. 

Pick future surgery arcs $a_1, \ldots, a_f$ and past surgery arcs $b_1, \ldots, b_p$ in $D_e$, such that consecutive circles in $D_e$ are connected by exactly one of the chosen arcs. We assume that in the path $P$, the $f$ edges directly after the $E$ edge correspond to performing surgery on $a_1, \ldots, a_f$, and the $p$ edges directly before the $E$ edge correspond to surgery on the dual arcs to $b_1, \ldots, b_p$. In other words, $P$ is of the form
\[
u \to \cdots \to u_0 \overset{b_1}{\to} u_1 \overset{b_2}{\to} \cdots \overset{b_p}{\to} e_0 \overset{E}{\to} e_1 \overset{a_1}{\to} v_1 \overset{a_2}{\to} \cdots \overset{a_f}{\to} v_f \to \cdots \to v.
\]

Consider the sub-cube $u_0\leq v_f$ and the path $P^2$ through this sub-cube, consisting of all edges in $P$ from $u_0$ to $v_f$. Let $D_{u_0}'$ denote the diagram obtained from $D_{u_0}$ by including only arcs in $P^2$. The $E$-resolution of $P^2$ is connected by construction, so $D_{u_0}'$ is connected as well. Let $P^1$ and $P^3$ denote the sub-paths of $P$ consisting of edges from $u$ to $u_0$ and $v_f$ to $v$, respectively. 

\[
\lefteqn{\overbrace{\phantom{u \to \cdots \to u_0}}^{P^1}}u \to \cdots \to
\underbrace{u_0 \overset{b_1}{\to} u_1 \overset{b_2}{\to} \cdots \overset{b_p}{\to} e_0 \overset{E}{\to} e_1 \overset{a_1}{\to} v_1 \overset{a_2}{\to} \cdots \overset{a_f}{\to} 
\lefteqn{\overbrace{\phantom{ v_f \to \cdots \to v}}^{P^3}}v_f}_{P^2} \to \cdots \to v
\]

Since $\pathPts{P}{\MspExy}$ is nonempty, there exist generators $x'\in \KG{u_0}$, $y'\in \KG{v_f}$ such that the moduli spaces
\[
\pathPts{P^1}\MspE{x'}{x},\; \pathPts{P^2}{\MspE{y'}{x'}},\; \pathPts{P^3} \MspE{y}{y'}
\]
are all nonempty.  By Lemma \ref{lem:E/F conn is one point per circle preliminary} applied to $D_{u_0}'$ and $P^2$, we can find $z^2\in  \pathPts{P^2}{\MspE{y'}{x'}}$ of the form
\[z^2 = (x' \rightarrow \cdots \rightarrow z_i \xrightarrow{j} z_{i+1} \rightarrow \cdots \rightarrow y')\in \pathPts{P^2}{\MspE{y'}{x'}}.\] 

Then choosing any $z^1 \in \pathPts{P^1}\MspE{x'}{x}$ and $z^3 \in \pathPts{P^3} \MspE{y}{y'}$ gives rise to at least one point in the composition moduli space
\[z=( z^3,z^2,z^1 ) \in \pathPts{P^3}{\MspE{y}{y'}} \times \pathPts{P^2}{\MspE{y'}{x'}} \times \pathPts{P^1}{\MspE{x'}{x}} \subset \pathPts{P}{\MspExy}\]
which $\varphi_P$ maps to the $j$-th point in $(0,1)$ as desired.

\end{proof}

\subsection{The $H$ map for connected arc diagrams}

As in the previous section, we continue to fix a subcube $u\leq_m v$ involving $J=H$ such that the initial arc diagram $D_u$ is connected (and hence all intermediate diagrams $D_w$ for $u\leq w \leq v$ are connected as well). The main technical result is a characterization of the path moduli spaces in $u\leq_m v$, Proposition \ref{prop:conn H}. We begin with several lemmas.

\begin{lemma}
\label{lem:conn is thin H}
Let $P$ be a path through $u\leq_m v$ such that the $H$-resolution $D_e$ consists entirely of $s$ essential circles. Fix generators $x\in \KG{u}$, $y\in \KG{v}$ such that $\pathPts{P}{\ModSp_H(y,x)}$ is nonempty. If every pair of consecutive circles in $D_e$ is joined by exactly one type of arc (either past or future), then the thick embedding $\varphi_P$ from \eqref{eq:product cube ambient space} is in fact thin.
\end{lemma}
\begin{proof}
Let $z,z'\in \pathPts{P}{\MspHxy}$ with 
\begin{align*}
z &= (x\rightarrow\cdots\rightarrow z_i \xrightarrow{j} z_{i+1} \rightarrow \cdots \to y),\\
z' &= (x\rightarrow\cdots\rightarrow z_i' \xrightarrow{j} z_{i+1}' \rightarrow \cdots \to y),
\end{align*}
for some $1\leq j \leq s$. Let $C$ denote the $j$-th circle in $D_e$. We will first show that $z_i = z_i'$ (and thus $z_{i+1} = z_{i+1}'$ as well).

Choosing a label of $v_{\pm}$ on $C$ uniquely determines the generator $z_i$ (and thus $z_{i+1}$), as in the proof of Proposition \ref{prop:E/F conn is thin}.  Therefore, in order to have $z_i\neq z_i'$, we must have the labels on $C$ in $z_i$ and $z_i'$ be opposites.  It follows that $z_i$ is obtained from $z_i'$ by swapping the labels on all circles in $D_e$, and likewise for $z_{i+1}$ and $z'_{i+1}$.

On the other hand, Lemma \ref{lem:forbidden configs} (Figures \ref{fig:forbidden5} and \ref{fig:forbidden6}) implies that the circles in $D_e$ are connected by past and future arcs in an alternating manner as in \eqref{eq:alternating arcs},
where there may be more than one past or future arc between any two adjacent circles.  Since the arc diagram $D_e$ is connected, all but the first and last essential circles in $D_e$ are connected to adjacent essential circles by both a past and future arc. Therefore at least one of $D_u$ or $D_v$ has an essential circle. Moreover, the labels on the essential circle(s) in $x$ (resp. $y$) are the same as the labels on the corresponding essential circle(s) in $z_i$ (resp. $z_{i+1}$).  Since $x$ and $y$ are fixed, at least one of these prevents the circles in $D_e$ from having opposite labels for either $z_i$ and $z_i'$, or $z_{i+1}$ and $z_{i+1}'$.

Thus we must have that $z_i=z_i'$ and $z_{i+1}=z'_{i+1}$, and the lemma follows by the same argument as in the proof of Proposition \ref{prop:E/F conn is thin}.

\end{proof}

Lemma \ref{lem:conn is thin H} does not handle all cases where $D_e$ consists of only essential circles.  After ruling out the forbidden configurations Figure \ref{fig:forbidden5} and Figure \ref{fig:forbidden6} from Lemma \ref{lem:forbidden configs}, it remains to analyze the situation where $D_e$ consists of exactly two essential circles with both a past and future arc connecting them, as in \eqref{eq:forbidden E configuration}. The following lemma gives an explicit description of such a path moduli space and its thick embedding.

\begin{lemma}
\label{lem:conn is thick H outlier}
Let $P$ be a path through $u\leq_m v$ such that the $H$-resolution $D_e$ consists entirely of two essential circles, with at least one past and one future arc connecting them. Fix generators $x\in \KG{u}$, $y\in \KG{v}$ such that $\pathPts{P}{\ModSp_H(y,x)}$ is nonempty. Then $\pathPts{P}{\MspHxy}$ consists of four points, and the thick embedding from \eqref{eq:product cube ambient space} is a two-to-one map onto the points $\left\{\frac 13,\frac 23\right\} \subset (0,1) \cong \AmbS$.  Moreover, any two points in $\pathPts{P}{\MspHxy}$ that map to the same point in $\AmbS$ are oppositely framed.
\end{lemma}

\begin{proof}

Let $a_p$ and $a_f$ denote the saddles corresponding to one past and one future surgery arc in $D_e$, respectively.  By Corollary \ref{cor:fixed J edge}, we may assume that $a_p$ and $a_f$ are performed directly before and after the $H$ edge, respectively.  Let $P'$ denote the 3-edge subpath of $P$ corresponding to $a_p$, $H$, and $a_f$, and note that the diagrams at the start and end of this path both consist of a single trivial circle.  Let $P_{\leq p}$ (resp. $P_{\geq f}$) denote the subpath consisting of the edges before $a_p$ (resp. after $a_f$).

It is a straightforward calculation to check that, in order for $\pathPts{P'}{\MspH{z_f}{z_p}}$ to be non-empty, the generator $z_p$ (resp. $z_f$) must consist of a label $1$ (resp. $X$) on the trivial circle at the start (resp. end) of the path $P'$.  In this unique case one further computes that $\pathPts{P'}{\MspH{z_f}{z_p}}$ consists of four points thickly embedding into $(0,1)$ as in the statement of the lemma.

Then for the full path moduli space we have
\[\pathPts{P}{\MspHxy} \cong \pathPts{P_{\geq f}}{\MspH{y}{z_f}} \times \pathPts{P'}{\MspH{z_f}{z_p}} \times
\pathPts{P_{\leq p}}{\MspH{z_p}{x}}.\]
The cobordism corresponding to the path $P_{\leq p}$ (resp. $P_{\geq f}$) must have genus zero, otherwise Lemma \ref{lem:genus considerations} would cause the path moduli space to be empty after composing with $a_p$ (resp. $a_f$) and connecting to $D_e$ which consists of all essential circles.  Thus both $\pathPts{P_{\leq p}}{\MspH{z_p}{x}}$ and $\pathPts{P_{\geq f}}{\MspH{y}{z_f}}$ consist of a single trivially embedded point and we are done.

\end{proof}

\begin{remark}
The two past and future arcs in Lemma \ref{lem:conn is thick H outlier} form a ladybug configuration (\cite[Figure 5.1]{LS}). Lemma \ref{lem:forbidden E configuration} disallows this configuration for $J=E, F$, and it is the only case among connected arc diagrams for which the thick embedding \eqref{eq:product cube ambient space} is not thin. The configuration in Lemma \ref{lem:forbidden configs}, Figure \ref{fig:forbidden1} is also a ladybug, which is disallowed for all $J$. 
\end{remark}

The following important lemma concerns the planar topology of an arc diagram $D_e$ as in the hypothesis of Lemma \ref{lem:conn is thick H outlier}.

\begin{lemma}
\label{lem:ray}
Let $P$ be a path through $u\leq_m v$ such that the $H$-resolution $D_e$ consists entirely of two essential circles, with at least one past and one future arc connecting them. Let $x\in \KG{u}$, $y\in \KG{v}$ be generators such that $\pathPts{P}{\ModSp_H(y,x)}$ is nonempty.  Then there exists a ray from the puncture $\times$ to the point at infinity which is disjoint from all the arcs in $D_e$ and which intersects the circles in $D_e$ in exactly two points. 
\end{lemma}

\begin{proof}
Let $a_p$ and $a_f$ denote two past and future arcs, respectively, which join the two essential circles in $D_e$. We assume that $a_p$ and $a_f$ are adjacent, in the sense that they can be connected by an interval whose interior is disjoint from the arc diagram $D_e$. We will show that there exists such a ray passing between them, as in Figure \ref{fig:cutting H configuration}. Indeed, the only possible obstruction to such a ray would be an arc, or a sequence of interlacing arcs, along one of the essential circles with endpoints `outside' of the arcs $a_p,a_f$ (see Figure \ref{fig:forbidden H configuration} with single arc case and Figure \ref{fig:forbidden linking configuration} for the interlacing arc case).  By considering the forbidden configurations of Figures \ref{fig:forbidden1} and \ref{fig:forbidden2}, we see that the arcs on one side of the essential circle must be either all past or all future.  In either case, we will arrive at another forbidden configuration as in Figure \ref{fig:forbidden3} or Figure \ref{fig:forbidden4}.  Thus no such arrangement may exist and our ray can be drawn.

\begin{figure}
\hspace*{\fill}
\begin{subfigure}[t]{.2\textwidth}
\begin{center}
\includestandalone{cutting_H_configuration}
\end{center}
\caption{The desired ray}\label{fig:cutting H configuration}
\end{subfigure}
\hfill
\begin{subfigure}[t]{.2\textwidth}
\begin{center}
\includestandalone{forbidden_configurations_for_H}
\end{center}
\caption{One obstructing arc}\label{fig:forbidden H configuration}
\end{subfigure}
\hfill
\begin{subfigure}[t]{.5\textwidth}
\begin{center}
\includestandalone{forbidden_linking_configuration}
\end{center}
\caption{An obstructing sequence of interlaced arcs}\label{fig:forbidden linking configuration}
\end{subfigure}
\hspace*{\fill}
\caption{}\label{fig:ray obstruction}
\end{figure}
\end{proof}

We are now ready for the main result in this section.

\begin{proposition}
\label{prop:conn H}
Consider a subcube $u\leq_m v$ such that the arc diagram $D_u$ is connected, and fix generators $x\in \KG{u}, y\in \KG{v}$. Then one of the following two holds.

\begin{enumerate}
\item 
\label{item:H path is thin}
For any path $P$ from $u$ to $v$ such that $\pathPts{P}\MspHxy \neq \varnothing$, the thick embedding from \eqref{eq:product cube ambient space} is in fact thin.

\item
\label{item:H path is thick}
For any path $P$ from $u$ to $v$ such that $\pathPts{P}\MspHxy \neq \varnothing$, the moduli space $\pathPts{P}{\MspHxy}$ consists of four points, and the thick embedding from \eqref{eq:product cube ambient space} is a two-to-one map onto the points $\left\{\frac 13,\frac 23\right\} \subset (0,1) \cong \AmbS$.  Moreover, any two points in $\pathPts{P}{\MspHxy}$ that map to the same point in $\AmbS$ are oppositely framed.
\end{enumerate}
\end{proposition}

\begin{proof}
If all paths yield empty moduli spaces then there is nothing to show. Take a path $P$ from $u$ to $v$ with $\pathPts{P}\MspHxy \neq \varnothing$. By Lemma \ref{lem:biject to all ess}, we may assume that the $H$-resolution $D_e$ of $P$ consists entirely of essential circles. We proceed by considering two cases for this fixed $P$. 

Suppose first that $D_e$ contains a pair of essential circles which are joined by both a past and future arc. We will argue that case (2) of the lemma holds. By Lemma \ref{lem:forbidden configs}, Figures \ref{fig:forbidden5} and \ref{fig:forbidden6}, we know that $D_e$ contains exactly two essential circles. It was shown in Lemma \ref{lem:conn is thick H outlier} that the conclusion of case (2) holds for $\pathPts{P} \MspHxy$, so it remains to show that the same holds for every path from $u$ to $v$. To that end, let $P'$ be another path from $u$ to $v$ with nonempty moduli space, and denote its $H$-resolution by $D_{e'}$. As usual, it suffices to assume that $D_{e'}$ contains only essential circles. The arc diagram $D_{e'}$ can be obtained from $D_e$ by performing a sequence of surgeries along the arcs in $D_e$. Then Lemma \ref{lem:ray} and the assumption that $\pathPts{P'}\MspHxy \neq \varnothing$ guarantees that $D_{e'}$ contains two essential circles. Moreover, the ray constructed in the proof of Lemma \ref{lem:ray}, along with the forbidden configurations of Lemma \ref{fig:forbidden configs}, Figures \ref{fig:forbidden1} and \ref{fig:forbidden2}, imply that the two circles in $D_{e'}$ are joined by both a past and a future arc.  Thus (2) holds by Lemma \ref{lem:conn is thick H outlier}.

Now suppose that every pair of consecutive circles in $D_e$ are joined by exactly one type (either past or future) of arc.  The above argument guarantees that the $H$-resolution of every path with nonempty moduli space also has the property that each pair of consecutive circles are joined by exactly one type of arc. Then we conclude that case (1) holds by Lemma \ref{lem:conn is thin H}. 
\end{proof}

Proposition \ref{prop:conn H} suffices to continue our construction. However, as in Section \ref{sec:The $E$ and $F$ maps for connected arc diagrams}, we can completely describe the thickly embedded path moduli spaces for $J=H$.  In Case \eqref{item:H path is thick} above, the moduli spaces are already described completely.  Meanwhile, for Case \eqref{item:H path is thin}, we have the following proposition; compare with Proposition \ref{prop:E/F conn is one point per circle}.

\begin{proposition}
\label{prop:H conn bijection}
Let $x\in \KG{u}$, $y\in \KG{v}$ such that we are in Case \eqref{item:H path is thin} of Proposition \ref{prop:conn H}.  Then for any path $P$ with $\de_P\MspHxy\neq\varnothing$, the thin embedding $\varphi_P$ from \eqref{eq:product cube ambient space} surjects onto $\left\{\frac{1}{s+1}, \ldots, \frac{s}{s+1}\right\}$ where $s$ denotes the number of essential circles in $D_e$.
\end{proposition}

\begin{proof}
Pick an element
\[
(x\rightarrow\cdots\rightarrow z_i \xrightarrow{j} z_{i+1} \rightarrow \cdots \to y) \in \pathPts{P}\MspHxy
\]
for some $1\leq j \leq s$. Since $H$ does not change the labels on essential circles, there must be an element of the form
\[
(x\rightarrow\cdots\rightarrow z_i \xrightarrow{k} z_{i+1} \rightarrow \cdots \to y) \in \pathPts{P}\MspHxy
\]
for every $1\leq k \leq s$. 
\end{proof}

\section{Constructing the higher dimensional moduli spaces (the inductive step)}
\label{sec:m-1 to m}

In order to complete the proof of Theorem \ref{thm:ffcJ}, we induct on the dimension of the subcube being considered.  That is to say, we fix $m>2$ and inductively assume that, for all subcubes $\subcube{i}{b}{a}$ with $i\leq m-1$ and generators $r\in\KG{a},s\in\KG{b}$, we have constructed moduli spaces $\MspJ{s}{r}$ together with thick embeddings
\[\MspJ{s}{r}\thickemb \AmbSp{b}{a}\]
satisfying all relevant boundary conditions.  For any $\subcube{i}{b}{a}$ that does not involve $J$, this thick embedding is just
\[\MspJ{s}{r} \cong \coprod \cubeMod{b}{a} \xthickemb{\coprod \id} \cubeMod{b}{a} = \AmbSp{b}{a}.\]
Meanwhile, for any $\subcube{i}{b}{a}$ that does involve $J$, the thick embedding is of codimension $1$ due to the extra factor of $(0,1)$ in the definition of $\AmbSp{b}{a}$.

Now we fix an $m$-dimensional subcube $\scube$ and seek to build new moduli spaces $\MspJxy$ thickly embedding into $\AmbS$ for various $x\in\KG{u},y\in\KG{v}$.  In the case that $\scube$ does not involve $J$, none of the subfaces of $\scube$ involve $J$ either and the entire construction follows from \cite{LS}.  In this case, the ``trivial cover" language of \cite{LS} gives a thickened identity embedding
\[\MspJxy \cong \coprod \cMod \xthickemb{\coprod \id} \cMod = \AmbS\]
and we are done.  Otherwise, the subcube $\scube$ involves $J$ and we split into cases.  

In each case, the essential strategy is the same.  The boundary of $\MspJxy$ is made up of pieces of various codimensional pieces $\de_{[i]}\MspJxy$, each of which consists of products of lower dimensional moduli spaces, which we inductively assume that we have already constructed and thickly embedded into the relevant ambient spaces.  When we assemble these pieces together, we have a thick embedding of $\de\MspJxy$ into $\de\AmbS$.  Completing the inductive step is then equivalent to proving that this thick embedding can be `thickly filled' by a moduli space embedding $\MspJxy\thickemb\AmbS$ that satisfies Equation \eqref{eq:respecting boundary structure}.

\subsection{Subcubes involving $E$ or $F$ with connected arc diagrams}
\label{sec:m-1 to m J=E conn}
In the case that $J=E$ (the case $J=F$ is again dual to this case and left to the reader), we begin by considering subcubes $\scube$ whose corresponding arc diagram $D_u$ is connected.  We fix generators $x\in\KG{u},y\in\KG{v}$ and consider the codimension $(m-1)$ (i.e. dimension $0$) part of the boundary
\[\de_{[m-1]}\MspExy = \coprod_{\text{paths $P$}} \pathPts{P}{\MspExy}.
\]

If all of these path moduli spaces are empty, then we define $\MspExy:=\varnothing$.  Otherwise, Proposition \ref{prop:E/F conn is thin} shows that any non-empty $\pathPts{P}{\MspExy}$ is \emph{thinly} embedded into $\pAmb\cong (0,1)$.

Now consider $\de_{[m-2]}\MspExy$, which Theorem \ref{thm:ffcJ} demands is of the form
\[
\de_{[m-2]} \MspExy = \bigcup \left( \coprod_{\substack{r\in\KG{a} \\s\in\KG{b}}}  \pathPts{P^b}{\MspE{y}{s}} \times \MspE{s}{r} \times \pathPts{P^a}{\MspE{r}{x}} \right),
\]
where the outer union is taken over
\begin{itemize}
    \item all square sub-faces $\subcube{2}{b}{a}$ in the subcube $\scube$, and
    \item all sub-paths $P^a$ from $u$ to $a$ and $P^b$ from $b$ to $v$.
\end{itemize}

We will denote the parenthetical term by $\sqMspExy$.  Then our (inductively defined) lower dimensional moduli space embeddings determine a thick product embedding
\begin{equation}\label{eq:E/F square face embedding}
\sqMspExy \thickemb \pathAmb{P^b} \times \AmbSp{b}{a} \times \pathAmb{P^a} \cong [0,1]\times(0,1)
\end{equation}
with boundary conditions based on the embeddings of the two path moduli spaces corresponding to composing $P^a$ and $P^b$ with the two paths around the edges of the square face $\subcube{2}{b}{a}$.

\begin{equation}
\begin{aligned}
\begin{tikzpicture}[commutative diagrams/every diagram, row sep = 1.4ex]

\matrix[matrix of math nodes, name=m, commutative diagrams/every cell] 
{  &        &   & c &   &       & \\
 u & \cdots & a &   & b &\cdots & v\\
   &        &   & c' &  &       & \\
};

\path[commutative diagrams/.cd,every arrow,every label]

(m-2-3) edge (m-1-4)
(m-2-3) edge (m-3-4)
(m-1-4) edge (m-2-5)
(m-3-4) edge (m-2-5)
(m-2-1) edge (m-2-2) (m-2-2) edge (m-2-3)  (m-2-5) edge (m-2-6) (m-2-6) edge (m-2-7)
;

\draw [
    thick,
    decoration={
        brace,
        mirror,
        raise=0.5cm
    },
    decorate
] (-5.1,.9) -- (-1.5,.9) 
node [pos=0.5,anchor=north,yshift=-0.55cm] {$P^a$};

\draw [
    thick,
    decoration={
        brace,
        mirror,
        raise=0.5cm
    },
    decorate
] (1.5,.9) -- (5.1,.9) 
node [pos=0.5,anchor=north,yshift=-0.55cm] {$P^b$};
\end{tikzpicture}
\end{aligned}
\end{equation}

The key point is that, because the path moduli spaces on either boundary of $\sqMspExy$ are \emph{thinly} embedded for our connected $D_u$, we can conclude that our composition moduli space embeddings \eqref{eq:E/F square face embedding} are thin as well.  Thus we see that all of $\de_{[m-2]}\MspExy$ is thinly embedded, and we move on to consider $\de_{[m-3]}\MspExy$.  This is again built from products of lower dimensional moduli spaces that we inductively assume have been constructed and thickly embedded in such a way that their boundaries are built out of parts of the thinly embedded $\de_{[m-2]}\MspExy$ - and thus $\de_{[m-3]}\MspExy$ must actually have been thinly embedded also.

\begin{figure}[ht]
\includegraphics[height=5cm]{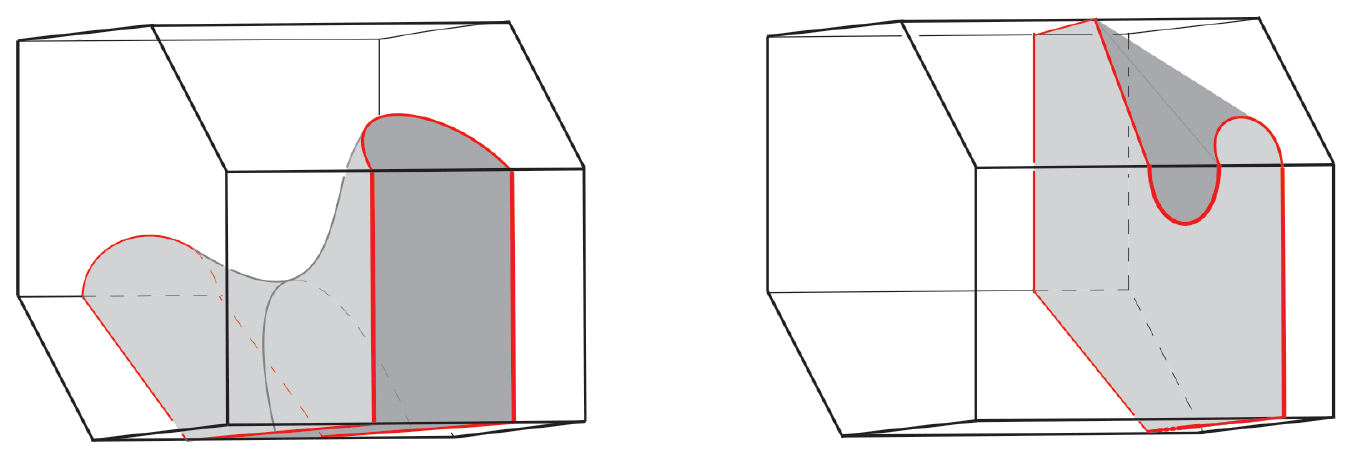}
\caption{Examples of 2-dimensional $\MspExy\thickemb\AmbS\cong\ModSp_\cube(v,u)\times(0,1)$, where $\ModSp_\cube(v,u)$ is the 2-dimensional permutohedron (a hexagon, topologically $D^2$).  The Pontryagin-Thom construction ensures that our framed embedded $\de\MspExy$ (drawn in red) can be filled as shown, regardless of the presence of turnbacks (framing data is omitted to avoid clutter).}
\label{fig:Hexagon times I}
\end{figure}

This reasoning continues all the way until we consider the thinly embedded $\de_{[1]}\MspExy=\de\MspExy$.  This is an $(m-2)$-dimensional manifold with corners $\de\MspExy$ thinly embedded into the $(m-1)$-dimensional ambient space boundary $\de\AmbS\cong S^{m-2}\times (0,1)$.  We may view this as a codimension-1 framed embedding $\de\MspExy\hookrightarrow S^{m-1}$ that misses the poles. Using the Pontryagin-Thom construction, this framed submanifold corresponds to a map $S^{m-1}\longrightarrow S^1$. Since $m>2$, the map is null-homotopic so it extends to a map $D^m\longrightarrow S^1$. Therefore $\de\MspExy\hookrightarrow S^{m-1}$ bounds a framed embedding of an $(m-1)$-dimensional manifold into the interior of $D^m$, which is equivalent to the interior of $D^{m-1}\times (0,1)\cong\AmbS$. It is worth noting that this is the point in the proof that uses (and propagates) the crucial codimension $1$ inductive assumption.  This is our construction of $\MspExy$ for connected $D_u$.  See Figure \ref{fig:Hexagon times I}.

\subsection{Subcubes involving $E$ or $F$ with disconnected arc diagrams}
\label{sec:m-1 to m J=E disconn}
Next we consider the situation when $J=E$ and the diagram $D_u$ is disconnected.  Unlike in Section \ref{sec:m-1 to m J=E conn}, path moduli spaces for subcubes with disconnected arc diagrams need not be thinly embedded. For example, adding a disjoint, contractible ladybug configuration (\cite[Figure 5.1]{LS}) has the effect of turning a nonempty thinly embedded path moduli space into one where the thick embedding is not injective. We begin with a simple lemma generalizing Lemma \ref{lem:m=2 conn comp and annular deg} from the $m=2$ construction.

\begin{lemma}\label{lem:general E acts on single conn comp}
Suppose there is some path $P$ for which $\pathPts{P}{\MspExy}\neq\varnothing$.  Then there is a single connected component $\Cs$ of the arc diagram $D_u$ such that
\[ \adeg(y|_\Cs)=\adeg(x|_\Cs)+2,\quad \adeg(y|_{\Cs'})=\adeg(x|_{\Cs'}),\]
where $\Cs'$ is the complement of $\Cs$.  We say that $E$ is acting on the component $\Cs$.  When $J=F$, the value $2$ above is replaced by $-2$.
\end{lemma}
\begin{proof}
As in the proof of Lemma \ref{lem:m=2 conn comp and annular deg}, saddles maps preserve annular degree, but the $E$ map increases annular degree by two.  
\end{proof}

As before, we consider
\[\MspExy=\MspE{y|_\Cs \otimes y|_{\Cs'}}{x|_\Cs \otimes x|_{\Cs'}}.\]

We let $\MspECCxy$ denote the moduli space that would be built for the diagram $\Cs$ for generators $x|_\Cs,y|_\Cs$ if $\Cs'$ were not present.  We also include the $E$ edge in the cube used to construct $\MspECCxy$, since $E$ is acting on $\Cs$.  Because $\Cs$ is connected, we can build $\MspECCxy$ using the techniques of Sections \ref{sec:0-dim moduli spaces}, \ref{sec:m=2 J=E/F involved}, and \ref{sec:m-1 to m J=E conn}.  Furthermore, if $\Cs'$ contains no saddles, then $y|_{\Cs'}=x|_{\Cs'}$ and we can define
\[\MspExy\cong\MspECCxy\]
as thickly embedded moduli spaces, and so we are done.

Therefore we may assume that $\Cs'$ contains at least one saddle.  In this case we let $\MspCCpxy$ denote the moduli space that would be built for the diagram $\Cs'$ for generators $x|_{\Cs'},y|_{\Cs'}$ if $\Cs$ were not present (not including the $E$ edge, since $E$ is not acting on $\Cs'$).  Note that, because $E$ is not involved, we can build $\MspCCpxy$ as a trivial cover of a suitable permutohedron using the techniques in \cite{LS}.

Now consider that, in the large subcube $\scube$, there are many $k$-dimensional subfaces $\subcube{k}{b}{a}$ whose edges consist of all of the saddles in $\Cs'$.  Indeed for any path $P$ starting from $u$ which consists of edges coming only from saddles in $\Cs$ (or the $E$ edge), the ending vertex of $P$ can be taken as the starting vertex $a$ of such a subface.

\begin{lemma}\label{lem:J=E subcube for the complement}
Let $k>0$ denote the number of saddles in $\Cs'$ and let $\subcube{k}{b}{a}$ denote \emph{any} subface of the subcube $\scube$ consisting of only edges corresponding to these $k$ saddles.  Then for any choice of generator $z$ on $\Cs\subset D_a$, we have a commuting diagram

\begin{equation}
\begin{aligned}
\label{eq:J=E framed bij for the complement}
\begin{tikzpicture}[x=1in,y=-.5in]

\node(A) at (0,0) {$\MspE{z\otimes (y|_{\Cs'})}{z\otimes (x|_{\Cs'})}$};

\node(B) at (0,2) {$\coprod \cubeMod{b}{a}$};

\node(C) at (2,0) {$\AmbSp{b}{a}$};

\node(D) at (2,2) {$\cubeMod{b}{a}$};

\node(E) at (0,1) {$\MspCCpxy$};

\path (A)-- node[rotate=-90]{$\cong$} (E);
\path (E)-- node[rotate=-90]{$\cong$} (B);
\path (C)-- node[rotate=-90]{$=$} (D);
\thickembtikz{A}{C}
\xthickembtikz{B}{D}{\coprod \id}

\end{tikzpicture}
\end{aligned}
\end{equation}

\end{lemma}

\begin{proof}
Because the arbitrary generator $z$ is kept constant throughout such a subcube $\subcube{k}{b}{a}$, the definitions make this clear.
\end{proof}

Implicit in Lemma \ref{lem:J=E subcube for the complement} is the fact that the diagram \eqref{eq:J=E framed bij for the complement} commutes for any portion of the boundaries of the moduli spaces in question.  Furthermore, this trivial cover remains consistent across all choices of $\subcube{k}{b}{a}$ and $z$.

The idea now is to assign a \emph{label} to each component of this trivial cover \[\MspCCpxy\xthickemb{\coprod \id} \cubeMod{b}{a}.\]  Then because $\Cs$ is disjoint from $\Cs'$, any such label will be maintained across the entire subcube $\scube$ and we are free to build the large moduli space $\MspExy$ one label at a time.

For each fixed label, we can regard $\MspCCpxy$ as a single copy of $\cubeMod{b}{a}$ embedding via the identity map.  Thus in any product used to build a composition moduli space, any factors coming from subfaces in $\Cs'$ contribute only trivial identity embeddings.  In this way we are free to build $\MspExy$ in the same manner as in the previous sections, constructing the boundary one codimension at a time.  At each stage of this process, products with the trivial identity embeddings maintain thinness and we eventually arrive at a thinly embedded $\de\MspExy$ of codimension one which can be filled using the Pontryagin-Thom construction as before.  Since this procedure is identical for each fixed label of component in $\MspCCpxy$, we conclude that our eventual thick embedding
\[\MspExy \thickemb \AmbS\]
contains a covering map of the same order as the trivial cover $\MspCCpxy\xthickemb{\coprod \id} \cubeMod{b}{a}$ and we are done.

\subsection{Subcubes involving $H$ with connected arc diagrams}
\label{sec:m-1 to m J=H conn}

In the case that $J=H$ and the subcube $\scube$ has corresponding arc diagram $D_u$ that is connected, Proposition \ref{prop:conn H} splits us into two further cases.  In case (1) of the proposition, we again have that all of our path moduli spaces are thin and the construction of $\MspHxy$ proceeds in precisely the same manner as in Section \ref{sec:m-1 to m J=E conn} when $J=E$.

Therefore we focus on case (2) of Proposition \ref{prop:conn H}, where for fixed $x\in\KG{u},y\in\KG{y}$, every non-empty path moduli space $\pathPts{P}\MspHxy$ consists of four points thickly embedded in two-to-one fashion onto $\left\{\frac 13 , \frac 23 \right\}\subset (0,1)\cong \AmbS$.  The proposition also shows that the two points embedded at any single $\frac{i}{3}$ are oppositely framed.  Thus for any path $P$ we can split $\pathPts{P}\MspHxy$ into two disjoint sets
\[\pathPts{P}\MspHxy=\pathPts{P}^+\MspHxy\sqcup\pathPts{P}^-\MspHxy,\]
where each $\pathPts{P}^\pm\MspHxy$ consists of two points, one embedded at $\frac{1}{3}$ with framing $\pm 1$, and the other embedded at $\frac{2}{3}$ with framing $\mp 1$.  In particular, each $\pathPts{P}^\pm\MspHxy$ is \emph{thinly} embedded into $\AmbS$.

The main idea then is to construct $\MspHxy$ in two `signed' pieces, which we denote by $\ModSp_H^\pm(y,x)$.  Each of these is built and thinly embedded into $\AmbS$ precisely as in Section \ref{sec:m-1 to m J=E conn} starting from the thinly embedded 
\[\de_{[m-1]}\ModSp_H^\pm(y,x):=\coprod_{\text{paths $P$}} \pathPts{P}^\pm\MspHxy,\]
and ending with a codimension one Pontryagin-Thom argument.  Note that the computations in Section \ref{sec:m=2 J=H involved} that are used to build $\de_{[m-2]}\MspHxy$ respect this decomposition into signed pieces, since any turnback described there must connect two oppositely framed points at $\frac 13$ and $\frac 23$.  The processes for building $\ModSp_H^+(y,x)$ and $\ModSp_H^-(y,x)$ are thus identical except that all framings are reversed, and we can combine them into a two-to-one thick embedding
\[\MspHxy=\ModSp_H^+(y,x)\sqcup\ModSp_H^-(y,x)\thickemb\AmbS.\]

\begin{remark}
The process above can be compared to Section \ref{sec:m-1 to m J=E disconn}, where identical thin constructions could be done one label at a time leading to thickly embedded moduli spaces.  The difference here is that our labels $+$ and $-$ give rise to oppositely framed moduli spaces throughout.
\end{remark}

\subsection{Subcubes involving $H$ with disconnected arc diagrams}
\label{sec:m-1 to m J=H disconn}

Finally we consider subcubes $\scube$ involving the $J=H$ edge, where the corresponding arc diagram $D_u$ is disconnected.  The construction in this case involves a mixture of the ideas in Sections \ref{sec:m-1 to m J=E disconn} and \ref{sec:m-1 to m J=H conn}, which we will refer to freely without reiterating any details.  As always, we begin by fixing $x\in\KG{u}$ and $y\in\KG{v}$ with some path $P$ through $\scube$ having non-empty $\pathPts{P}\MspHxy$.

Suppose the arc diagram $D_u$ has $k$ homologically essential connected components, which we denote by $\Cs_1,\dots,\Cs_k$, ordered by their nesting in the annulus.

The first idea is to build $\MspHxy$ one homologically essential connected component at a time, in the spirit of Section \ref{sec:m-1 to m J=H conn}.  To this end, note that each path moduli space $\pathPts{P}\MspHxy$ can be decomposed as
\[\pathPts{P}\MspHxy = \pathPts{P}^{\Cs_1}\MspHxy \sqcup \cdots \sqcup \pathPts{P}^{\Cs_k}\MspHxy,\]
where each $\pathPts{P}^{\Cs_i}\MspHxy$ consists of only the points corresponding to $H$ acting on essential circles in $\Cs_i$.  This decomposition is clearly respected throughout the entire cube, allowing us to decompose our desired moduli space $\MspHxy$ as
\[\MspHxy = \ModSp_H^{\Cs_1}(y,x) \sqcup \cdots \sqcup \ModSp_H^{\Cs_k}(y,x).\]
Then fixing a connected component $\Cs_i$, we  seek to build and thickly embed the space $\ModSp_H^{\Cs_i}(y,x)$ starting from its codimension $(m-1)$ boundary
\[\de_{[m-1]}\ModSp_H^{\Cs_i}(y,x) :=
\coprod_{\text{paths $P$}} \pathPts{P}^{\Cs_i}\MspHxy\]
as usual.

In order to do this, we appeal to the arguments in Section \ref{sec:m-1 to m J=E disconn} letting $\Cs_i$ take the role of $\Cs$ as notated there.  If the component $\Cs_i$ corresponds to case (1) of Proposition \ref{prop:conn H}, then this construction proceeds precisely as before and we produce a trivial cover according to the usual Lipshitz-Sarkar construction for generators on $\Cs_i'$, the complement of $\Cs_i$.  If the component $\Cs_i$ corresponds to case (2) of Proposition \ref{prop:conn H}, then the procedures of Section \ref{sec:m-1 to m J=E disconn} are done one signed label at a time, leading to a trivial cover with twice the order of the covering coming from Lipshitz-Sarkar's construction for generators on $\Cs_i'$.

\begin{figure}
\includegraphics[height=5cm]{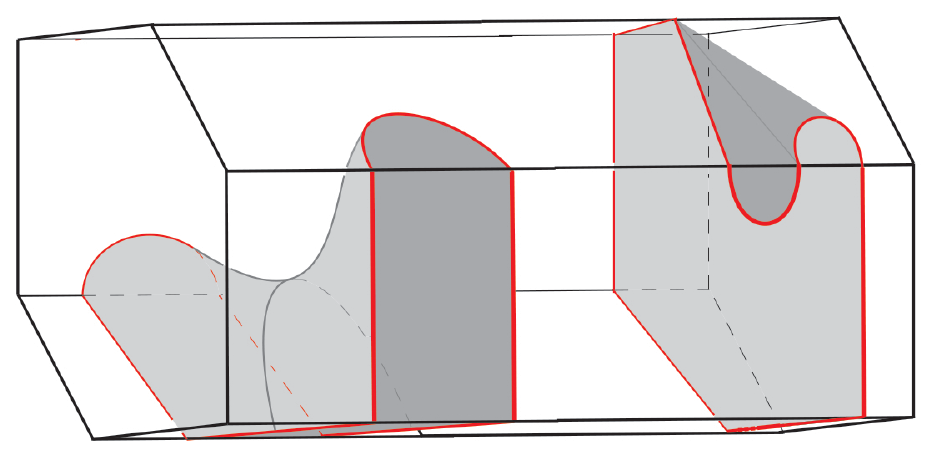}
\caption{A possible embedding for $\MspHxy$ for a disconnected diagram.  The two components are constructed and embedded individually for $H$ acting on different connected components of the arc diagram $D_u$.  They are then horizontally `stacked' together as shown.}
\label{fig:Hexagon times I stacked}
\end{figure}

In any case, each moduli space $\ModSp_H^{\Cs_i}(y,x)$ can be built and thickly embedded into a copy of $\AmbS$ which we write as
\[\ModSp_H^{\Cs_i}(y,x)\thickemb\AmbS\cong \ModSp_\cube \times (i-1,i).\]  We then `stack' all of these thick embeddings together as
\[\MspHxy = \coprod_{i=1}^k \ModSp_H^{\Cs_i}(y,x) \thickemb \coprod_{i=1}^k \ModSp_\cube \times (i-1,i) \hookrightarrow \ModSp_\cube \times (0,k) \cong \AmbS\]
to complete our construction.  See Figure \ref{fig:Hexagon times I stacked}.

\section{Invariance and naturality with respect to link cobordisms}
\label{sec:invariance}

In order to complete the proof of Theorem \ref{main thm}, we must show that our construction of the map $\J:\X_\A(D)\rightarrow\X_\A(D)$, or equivalently our spectrum $\X_J(D)$ lifting $\Cone(J)$, is invariant under the choices involved.  The majority of these choices are considered in \cite{LS} where invariance of the similarly constructed $\X_{Kh}(D)$ is shown, and the proofs there apply to our construction in the same manner.  However, in Section \ref{sec:m-1 to m} we made an additional choice of null-cobordism while appealing to the Pontryagin-Thom construction.  We show below that $\X_J(D)$ is invariant under this choice as well.

Finally, we will also show that $\J$ commutes, up to homotopy, with the maps on spectra assigned to link cobordisms.  Proposition \ref{prop:J commutes with cobs} below will concern elementary link cobordisms, including Reidemeister moves, and will thus conclude the proofs of both Theorem \ref{main thm} and Theorem \ref{thm:J commutes with cobs}.

\subsection{Invariance under choices of cobordisms}
The inductive step in the construction of higher-dimensional moduli spaces (see Section \ref{sec:m-1 to m}) relied on choosing a codimension-1 framed null-cobordism for the boundary. We will now show that different choices are in turn related by a framed codimension-1 cobordism.  (Even though we work with codimension-1 orientable submanifolds where the framing can always be assumed to exist, we choose to explicitly specify the framing since it is central to the construction in the context of framed flow categories.) 

\begin{lemma} \label{lem:cob choices linked by higher cob}
Suppose $\MspJxy,\ModSp'_J(y,x)$ are two moduli spaces constructed as in Section \ref{sec:m-1 to m} using some choices of null-cobordisms.  Then there exists a codimension-1 cobordism, that is a thick embedding
\[\overline{\ModSp}_J(y,x)\thickemb \AmbS\times[0,1]\]
whose boundary image in $\AmbS\times\{0\}$ (respectively $\AmbS\times\{1\}$) recovers the thick embedding $\MspJxy\thickemb\AmbS$ (respectively $\ModSp'_J(y,x)\thickemb\AmbS$).
\end{lemma}

\begin{proof} 
The construction of $\overline{\ModSp}_J(y,x)$ is done inductively in parallel with that of the moduli spaces. There is no choice involved in the construction of $0$- and $1$-dimensional moduli spaces, so in these cases $\overline{\mathcal{M}}$ is defined to be the product $\mathcal M \times I$.  For the inductive step we focus first on the case when $J=E$ (or $F$) for connected diagrams as in Section \ref{sec:m-1 to m J=E conn}.

Recall that we begin with the $(m-2)$-manifold with corners $\de\MspExy$ (for $m>2$) thinly embedded into $\de\AmbS$, viewed as $S^{m-1}$ minus the poles. Consider $\de\MspExy\subset S^{m-1}\times\{0\}$, $\de\mathcal{M}'_E(x,y)\subset S^{m-1}\times\{1\}$.  Just as in the arguments of Section \ref{sec:m-1 to m J=E conn}, the inductive assumption of this lemma provides the existence of a manifold with corners (denote it $\partial \overline{\mathcal{M}}_E(x,y)$) thinly embedded as a framed cobordism for these boundaries in $S^{m-1}\times[0,1]$.   Assembling this data: \[\MspExy\cup \partial \overline{\mathcal{M}}_E(x,y)\cup \mathcal{M}'_E(x,y),\]
we have a codimension-1 framed submanifold of $\partial(D^m\times I)$. The corresponding map $\partial(D^m\times I)\longrightarrow S^1$ extends to $D^m\times I$, giving a desired cobordism in $\AmbS\times[0,1]$ between $\MspExy, \mathcal{M}'_E(x,y)$. Note that it extends the cobordism $\partial \overline{\mathcal{M}}_E(x,y)$  in $\partial \AmbS\times[0,1]$ that was given by the inductive assumption.

The extension to the other cases, where the embeddings are indeed thick, is done in direct analogy to the constructions in Sections \ref{sec:m-1 to m J=E disconn}, \ref{sec:m-1 to m J=H conn}, and \ref{sec:m-1 to m J=H disconn}, where all components of the trivial cover can be considered one at a time.
\end{proof}

Next we will use Lemma \ref{lem:cob choices linked by higher cob} to show that stable homotopy types constructed using different choices of moduli spaces are homotopy equivalent. 

\begin{proposition}
\label{prop:invariance under cobs}
Fix an annular link diagram $D$, and let $\X_J(D),\X_J'(D)$ denote two stable homotopy types built via Theorem \ref{thm:ffcJ} using different choices of null-cobordisms to construct the moduli spaces as described above.  Then $\X_J(D)\simeq\X_J'(D)$.
\end{proposition}

\begin{proof}
We work in the setting of \cite[Section 3.3]{LS}. An analogue of Lemma 3.16 in that reference applied to our construction in Lemma \ref{lem:cob choices linked by higher cob}, together with the ability to `separate sheets of the cover' as in the proof of Proposition \ref{prop:cubical flow is framed}, gives a collection of embeddings $\overline{\mathcal M}_J(y,x)\subset \mathbb{E}_d[{\rm gr}(y):{\rm gr}(x)]\times[0,1]$ which extend the neat embeddings $\iota_{y,x}, \iota'_{y,x}$ of ${\mathcal M}_J(y,x), {\mathcal M}'_J(y,x)$ into $\mathbb{E}_d[{\rm gr}(y):{\rm gr}(x)]\times\{ 0\}$, respectively $\times \{ 1\}$. Informally, this may be thought of as a $1$-parameter family of neat embeddings of flow categories, with singularities at times corresponding to critical points of the cobordisms $\overline{\mathcal M}$.

Consider the construction of the cell complex corresponding to a framed flow category \cite[Definition 3.24]{LS}.  To each object $y$ in the flow category, one associates a cell $\mathcal{C}(y)$.  To define the attaching map of $\mathcal{C}(y)$ to a lower dimensional cell $\mathcal{C}(x)$, one uses the neat embedding of the moduli space $\ModSp(y,x)$ to identify a subset $\mathcal{C}_x(y)\cong\ModSp(y,x)\times \mathcal{C}(x)\subset \de\mathcal{C}(y)$.  The map is then defined on $\mathcal{C}_x(y)$ to be the projection to $\mathcal{C}(x)$, while sending all of $\de {\mathcal C}(y)\smallsetminus \cup_x {\mathcal C}_x(y)$ to the basepoint.

In our setting, the neat embedding of $\overline{\ModSp}_J(y,x)$ into $\mathbb{E}_d[{\rm gr}(y):{\rm gr}(x)]\times[0,1]$ provides an identification of $\overline{\ModSp}_J(y,x)\times\mathcal{C}(x)\subset(\de\mathcal{C}(y))\times[0,1]$, and the corresponding projection $\overline{\ModSp}_J(y,x)\times\mathcal{C}(x) \longrightarrow \mathcal{C}(x)$ gives the desired homotopy between attaching maps corresponding to the neat embeddings of $\ModSp_J(y,x)$ and $\ModSp'_J(y,x)$.  The constructions of the thickly embedded $\overline{\ModSp}_J(y,x)$ in Lemma \ref{lem:cob choices linked by higher cob} give rise to homotopies of this sort which are coherent across all attaching maps of all cells, showing that our resulting spectra $\X_J(D)$ and $\X'_J(D)$ are stably homotopy equivalent.
\end{proof}

\subsection{Link cobordisms and Reidemeister invariance}

Our goal in this section is to show that the $\J$ map constructed in Corollary \ref{cor:ffcj refines Cone J} commutes up to homotopy with the maps assigned to link cobordisms in \cite[Section 3]{LS2}.  This will include the case of link isotopies, indicating that $\J$ commutes with the stable equivalences assigned to Reidemeister moves in \cite[Section 6]{LS}.  We begin with some general notions about certain subcategories of flow categories that are used to construct the desired maps.

Let $\ffc$ be a framed flow category. Recall from \cite[Definition 3.29]{LS} the notion of a downwards (resp. upwards) closed subcategory $\subcat \subset \ffc$, which is a certain full subcategory whose geometric realization $\lr{\subcat}$ is naturally a subcomplex (resp. quotient complex) of $\lr{\ffc}$. Note that any full subcategory is specified by its objects. Thus if $\subcat \subset \ffc$ is downward closed, then it has a complementary upward closed subcategory $\subcat'$ whose objects are $\Obj(\ffc)\setminus \Obj(\subcat)$. Moreover, upon geometric realization, $\lr{\subcat}$ is a subcomplex of $\lr{\ffc}$ with quotient complex $\lr{\subcat'}$, allowing us to write the cofibration sequence
\begin{equation}
\label{eq:complementary cofib seq}
\lr{\subcat}\hookrightarrow \lr{\ffc} \twoheadrightarrow \lr{\subcat'}.    
\end{equation}
We will say a full subcategory $\subcat\subset \ffc$ is \emph{closed} if it is either upwards or downwards closed. 

We now specialize to the situation at hand.  Let $D$ be an annular link diagram, and let $\ffcJ(D)$ be the flow category refining $\Cone~J$ (see Theorem \ref{thm:ffcJ}).  By definition, the objects of $\ffcJ(D)$ consist of two copies of the objects of $\ffc_\A(D)$, 
\[
\Obj(\ffcJ(D)) = \Obj(\ffc_\A(D))_0 \coprod \Obj(\ffc_\A(D))_1,
\]
where we use the subscript $0$ and $1$ to indicate that some objects are before the $J$ map (subscript $0$), and some are after (subscript $1$). The gradings of objects in $\Obj(\ffc_\A(D))_1$ are shifted up by $1$ from the gradings in $\Obj(\ffc_\A(D))_0 $.

\begin{definition}
For a full subcategory $\subcat \subset \ffc_\A(D)$, let $\overline{\subcat}$ denote the full subcategory of $\ffcJ(D)$ whose objects are two copies of the objects of $\subcat$, 
\[
\Obj(\overline{\subcat}) = \Obj({\subcat})_0 \coprod \Obj({\subcat})_1,
\]
where $\Obj({\subcat})_i \subset \Obj(\ffc_\A(D))_i$. We call $\overline{\subcat}$ the \emph{double} of $\subcat$, and we let $\subcat_0,\subcat_1$ denote the two full subcategories of $\overline{\subcat}\subset\ffcJ(D)$ defined by
\[
\Obj(\subcat_i) = \Obj(\subcat)_i.
\]
\end{definition}

As with $\ffcJ(D)$, we think of $\subcat_0$ and $\subcat_1$ as two parallel copies of $\subcat$ within the double $\overline{\subcat}$, separated by the $J$ edge direction.  In this way, $\subcat_0$ and $\subcat_1$ are complementary downward and upward closed subcategories of $\overline{\subcat}$, so that upon geometric realization, there is an inclusion and projection 
\[
\lr{\subcat_0} \hookrightarrow \lr{\overline{\subcat}} \twoheadrightarrow \lr{\subcat_1}.
\]
Combining this with the canonical\footnote{The subcategories $\subcat_0, \subcat_1$ do not involve the $J$ edge, so they are built in exactly the same way as $\subcat$.} identifications $\subcat_0=\subcat, \subcat_1 = \Sigma\subcat$ (the $\Sigma(-)$ denotes a grading shift as in Definition \ref{def:grading shift ffc}) allows us to write, for any full subcategory $\subcat\subset\ffc_\mathbb{A}(D)$, the cofibration sequence
\begin{equation}
\label{eq:doubles cofib seq}
    \lr{\subcat} \hookrightarrow \lr{\overline{\subcat}} \twoheadrightarrow \Sigma\lr{\subcat}.
\end{equation}

Meanwhile, if $\subcat'$ denotes the complementary full subcategory of $\subcat$ in $\ffc_\A(D)$, then the complement of the double is the double of the complement, $\overline{\subcat}'=\overline{\subcat'}$.
However, if $\subcat$ is closed in $\ffc_\A(D)$, then its double $\overline{\subcat}$ need not be closed in $\ffcJ(D)$, since the $J$ map may take certain generators corresponding to objects of $\subcat$ into $\subcat'$, or vice versa.

\begin{lemma}
\label{lem:combined cofib seqs}
Suppose that $\subcat$ is a closed subcategory of $\ffc_\mathbb{A}(D)$ whose double $\overline{\subcat}$ is also a closed subcategory of $\ffcJ(D)$.  Then the cofibration sequences \eqref{eq:complementary cofib seq} for complementary pairs $\subcat,\subcat'$ and $\overline{\subcat},\overline{\subcat}'$ commute with the cofibration sequences \eqref{eq:doubles cofib seq} for complementary copies within a double, as illustrated below for the case that $\subcat$ is downward closed.
\begin{equation}
\label{eq:commuting cofib seqs}
\begin{tikzpicture}[x=1in,y=-.6in,baseline=(current  bounding  box.center)]

\node(AA) at (0,0) {$\lr{\subcat}$};
\node(AB) at (1,0) {$\lr{\overline{\subcat}}$};
\node(AC) at (2,0) {$\Sigma\lr{\subcat}$};

\node(BA) at (0,1) {$\lr{\ffc_\mathbb{A}(D)}$};
\node(BB) at (1,1) {$\lr{\ffcJ(D)}$};
\node(BC) at (2,1) {$\Sigma\lr{\ffc_\mathbb{A}(D)}$};

\node(CA) at (0,2) {$\lr{\subcat'}$};
\node(CB) at (1,2) {$\lr{\overline{\subcat}'}$};
\node(CC) at (2,2) {$\Sigma\lr{\subcat'}$};

\draw[right hook->] (AA)--(AB);
\draw[right hook->] (BA)--(BB);
\draw[right hook->] (CA)--(CB);
\draw[right hook->] (AA)--(BA);
\draw[right hook->] (AB)--(BB);
\draw[right hook->] (AC)--(BC);

\draw[->>] (AB)--(AC);
\draw[->>] (BB)--(BC);
\draw[->>] (CB)--(CC);
\draw[->>] (BA)--(CA);
\draw[->>] (BB)--(CB);
\draw[->>] (BC)--(CC);

\end{tikzpicture}
\end{equation}

\end{lemma}

\begin{proof}
The inclusion maps, the quotient maps, and the canonical identifications involved are all defined by identifying objects of the categories with cells in the geometric realizations, so the definitions make the commutation clear.
\end{proof}

Note that, whenever we do have a closed subcategory $\subcat\subset\ffc_\mathbb{A}(D)$ whose double $\overline{\subcat}$ remains closed in $\ffcJ(D)$ as in Lemma \ref{lem:combined cofib seqs}, each row in \eqref{eq:commuting cofib seqs} gives rise to a `horizontal' Puppe map which can be identified as the $\J$ map restricted to the setting of its row.  Similarly, each column in \eqref{eq:commuting cofib seqs} also gives rise to a `vertical' Puppe map.

An annular link cobordism can be decomposed into a sequence of so-called \emph{elementary cobordisms} corresponding to either Morse moves or Reidemeister moves.  To prove Theorem \ref{thm:J commutes with cobs}, it suffices to show that the map $\J$ commutes up to homotopy with all maps assigned to such elementary cobordisms.

\begin{proposition}
\label{prop:J commutes with cobs}
Let $\phi:\X_\mathbb{A}(D)\rightarrow\X_\mathbb{A}(D')$ be the map on annular Khovanov spectra assigned to an elementary annular link cobordism from $D'$ to $D$ as in \cite[Section 3]{LS2}.  Then the map $\J$ commutes with $\phi$ up to homotopy:
\[
\begin{tikzcd}
\X_\A(D) \ar[r, "\J"] \ar[d, "\phi"'] & \X_\A(D) \ar[d, "\phi"] \\
\X_\A(D') \ar[r, "\J"] & \X_\A(D')
\end{tikzcd}
\]
In particular, if $D,D'$ are two diagrams for the same annular link, then the stable equivalence $\X_\A(D) \cong \X_\A(D')$ commutes with $\J$ up to homotopy.
\end{proposition}

\begin{proof}
Each elementary cobordism map arises from identifying a closed subcategory $\subcat$ in the framed flow category $\ffc_\mathbb{A}(\td{D})$ for a suitable diagram $\td{D}$, leading to a cofibration sequence of the form \eqref{eq:complementary cofib seq} (or in the case of Reidemeister II and III, a finite sequence of such closed subcategories).  We claim that each such closed subcategory $\subcat$ used in this way gives rise to a closed double $\overline{\subcat}\subset \ffcJ(\td{D})$.  Indeed the saddle map arises from a diagram $\td{D}$ having a crossing at the place of the saddle, so that $\subcat$ and its complement $\subcat'$ correspond to 0- and 1-resolutions of this crossing, which clearly remain closed upon doubling.  Meanwhile, all of the other elementary cobordisms (cups, caps, and Reidemeister moves) identify closed subcategories $\subcat$ by fixing labels on \emph{trivial} circles in their corresponding $\td{D}$.  Since the $J$ map does not affect the labels on trivial circles, any such closed subcategory remains closed upon doubling.

Thus any closed subcategory $\subcat$ used to define an elementary cobordism map $\phi$ gives rise to a commuting diagram of cofibrations as in \eqref{eq:commuting cofib seqs}.  In all such cases, the map being built will be identified with one of the vertical maps in the first and/or last column (or perhaps the corresponding `vertical' Puppe map), while the $\J$ maps will be identified with the `horizontal' Puppe maps, and these two maps will commute due to the naturality of the Puppe construction.
\end{proof}

\begin{remark}
Conjecturally the maps $\phi$ assigned to link cobordisms in \cite[Section 3]{LS2} are well-defined and do not depend (up to homotopy) on the decomposition into elementary cobordisms. If this is the case, the statement of Theorem \ref{thm:J commutes with cobs} can be similarly strengthened as well. 
\end{remark}

\bibliographystyle{alpha}
\bibliography{main}

\end{document}